\documentclass[10pt]{article}
\newcommand{\hatw}{\hat}
\usepackage{geometry}
\usepackage[utf8]{inputenc}
\usepackage{amsfonts, amsmath, amssymb, amsthm}
\usepackage{physics}
\usepackage{graphicx}
\usepackage{wrapfig}
\usepackage[hidelinks]{hyperref}
\usepackage{makecell}
\usepackage{float}
\usepackage[title]{appendix}
\usepackage[justification=centering]{caption}

\usepackage{tikz}
\usetikzlibrary{arrows.meta}

\usepackage{fancyhdr}
\usepackage{datetime}
\usepackage[backend=bibtex,style=numeric,maxbibnames=99,sorting=none]{biblatex}
\renewbibmacro{in:}{}
\usepackage{mathtools}

\hypersetup{
    colorlinks=true,
    linkcolor=blue,
    filecolor=magenta,
    urlcolor=blue,
}

\mathtoolsset{showonlyrefs}

\renewcommand{\norm}[1]{\left\lVert#1\right\rVert}

\DeclareMathOperator{\Div}{div}
\DeclareMathOperator{\diam}{diam}

\DeclareMathOperator{\sgrad}{sgrad}
\DeclareMathOperator{\Int}{int}

\renewcommand{\dv}[2]{\frac{d#1}{d#2}}

\def\sloppy{%
  \hfuzz 3pt
  \vfuzz=\hfuzz}
\newtheorem{theorem}{Theorem} %
\newtheorem{corollary}[theorem]{Corollary}
\newtheorem{lemma}[theorem]{Lemma}
\newtheorem{remark}[theorem]{Remark}
\newtheorem{proposition}[theorem]{Proposition}
\newtheorem{definition}[theorem]{Definition}
\numberwithin{equation}{section}
\numberwithin{theorem}{section}
\author{Anatoly Neishtadt and Alexey Okunev}
\date{}
\title{Averaging and passage through resonances in two-frequency systems near separatrices\thanks{The work was supported by the Leverhulme Trust (Grant No. RPG-2018-143).}}
\begin{document}
\maketitle
\abstract{
The averaging method is a classical powerful tool in perturbation theory of dynamical systems. There are two major obstacles to applying the averaging method, resonances
and separatrices. In this paper we obtain realistic asymptotic estimates that justify the use of averaging method in a generic situation where both these obstacles are present at the same time, passage through a separatrix for time-periodic perturbations of one-frequency Hamiltonian systems. As a general phenomenon, resonances accumulate at separatrices. The Hamiltonian depends on a parameter that slowly changes for the perturbed system (so slow-fast Hamiltonian systems with two and a half degrees of freedom are included in our class).
Our results can also be applied to perturbations of generic two-frequency integrable systems near separatrices, as they can be reduced to periodic perturbations of one-frequency systems.
}

\tableofcontents

\section{Introduction}

Small perturbations of integrable Hamiltonian systems is an important class of dynamical systems encountered in various applications.
Far from the separatrices of the unperturbed system one can use the action-angle variables of the unperturbed system. The values of action variables slowly change for the perturbed system and their evolution can be approximately described by the averaged system that is obtained by averaging the rate of change of the action variables over all values of the angle variables.
For perturbations of one-frequency systems, the averaged system describes the evolution of the action variables for all initial data with accuracy $O(\varepsilon)$ over times $\sim \varepsilon^{-1}$ (\cite{fatou1928mouvement}, \cite{bogolyubov61}). Here $\varepsilon$ is the small parameter of perturbed system.

For two-frequency systems, the averaged system describes the evolution of most initial data (except a set with measure $O(\sqrt{\varepsilon})$)
with accuracy $O(\sqrt{\varepsilon} |\ln \varepsilon|)$ over times $\sim \varepsilon^{-1}$
(\cite{neishtadt1975passage}, see also review~\cite{neishtadt2014} and references therein and earlier work~\cite{arnold65}, where two-frequency systems were studied under a condition that prohibits capture into resonances).
Resonances between two frequencies of the problem are possible, and, as these frequencies change along solutions of the perturbed system, it is possible that the ratio between the frequencies remains near a resonant value for times $\sim \varepsilon^{-1}$. This phenomenon is called \emph{capture into resonance}, and there are examples (e.g., in~\cite{neishtadt2014}) when this happens for a set of initial data with measure $~\sim \sqrt{\varepsilon}$. Capture into resonances is the reason for the exceptional set with measure $O(\sqrt{\varepsilon})$.
Most trajectories are not captured into resonances, but still passing through resonances leads to a jump of order $\sqrt{\varepsilon}$ (for most trajectories), this is called \emph{scattering on resonance}. This is why accuracy of averaging method is only $O(\sqrt{\varepsilon} \ln \varepsilon)$ (extra logarithm appears because for some trajectories outside the exceptional set scattering is larger than $O(\sqrt{\varepsilon})$).

Results that hold for multi-frequency systems are weaker.
Very general results~\cite{anosov1960averaging, kasuga1961adiabatic} about averaging in slow-fast systems imply that for any $\rho > 0$ the measure of the set of initial data such that accuracy of averaging method is worse than $\rho$ is $o(1)_{\varepsilon \to 0}$.
Restriction of the generality allows to estimate how this measure depends on $\rho$ and $\varepsilon$ (\cite[\S 6.1.9]{arnold2007mathematical} and references therein).

The phase space of the unperturbed system is often divided into several domains by separatrices. Solutions of the perturbed system can cross separatrices of the unperturbed system and move between these domains.
Effects of separatrix crossing are well studied for one-frequency systems (\cite{neishtadt17} and references therein).
Importantly, separatrix crossing leads to probabilistic phenomena.
For example, trajectory of a point moving in a double-well potential (unperturbed system) with small friction (perturbation) eventually exhibits separatrix crossing and remains bounded in one of the two wells. Trajectories caught in each well are finely mixed in the phase space when friction is small, thus capture in each well can be considered a "random event" with definite probability (see the discussion in~\cite{arnold63},~\cite{neishtadt17}).
One can modify averaged system to cover trajectories crossing separatrices: when the solution of averaged system in some domain reaches the separatrices (and thus averaged system in this domain is no longer defined), one can write averaged system in another domain bounded by the same separatrices and continue this solution using averaged system in the other domain.
When capture in several domains is possible after separatrix crossing, there are modified averaged systems describing capture into each such domain.
The evolution of most initial data (with exceptional set having measure $O(\varepsilon^r)$ for any $r > 0$, this set corresponds to solutions that come too close to the saddle of perturbed system) is described by a solution of modified averaged system with accuracy $O(\varepsilon |\ln \varepsilon|)$ (\cite{neishtadt17}).

Separatrix crossing is much less studied for perturbations of Hamiltonian systems with two and more frequencies.
Separatrix crossing for time-periodic perturbations of one-frequency systems of special form (periodically forced by a single harmonic and weakly damped motion in a double-well potential) is considered in~\cite{brothers1998slow} under an assumption that periodic forcing is sufficiently small so that captures into resonances close to separatrices do not occur. Authors use multiphase averaging between resonances and single-phase averaging near resonances to obtain formulas for the boundaries of the sets of initial data captured in each well (without rigorous justification).
Stochastic perturbations of time-periodic perturbations of one-frequency systems were studied in~\cite{wolansky1990limit}.

The goal of this paper is to obtain realistic estimates for the accuracy of averaging method for time-periodic perturbations of one-frequency systems with separatrix crossing. The unperturbed Hamiltonian can depend on a parameter that slowly changes for perturbed system.
This is a generic situation when there are both resonances and separatrix crossing. Perturbations of generic two-frequency systems can be reduced to this case, see Section~\ref{s:two-freqs} below.
Phase angle on the closed phase trajectories of
the unperturbed system and time are two angle variables, and resonances
between their frequencies are possible.
Far from separatrices results on two-frequency systems mentioned earlier are applicable.
Resonances accumulate on separatrices and resonances near separatrices have to be studied separately.
We prove that accuracy of averaging method $O(\sqrt{\varepsilon} |\ln \varepsilon|)$ holds for most initial conditions, the exceptional set has measure $O(\sqrt{\varepsilon} |\ln^5 \varepsilon|)$.
We also prove formulas for "probabilities" of proceeding into different domains after separatrix crossing similar to the formulas that hold in one-frequency case~\cite{neishtadt17}.

A natural and frequently encountered in applications subclass of systems we consider is the motion of a particle in a double-well or periodic potential with small friction and time-periodic forcing (there might also be a slow change of parameters).
One example of such system is the system describing planar librational movement of an arbitrary shaped satellite in an elliptic orbit (cf., e.g., \cite[Problem~1.2.19]{wiggins1990introduction}; dissipation caused by tidal friction can be added to this problem).
More examples can be found in~\cite{brothers1998slow}.

Perturbations of two-frequency Hamiltonian systems can be reduced to our case, cf. Section~\ref{s:two-freqs}.
Examples of two-frequency integrable systems with separatices include Euler top, geodesic flows on ellipsoid or a surface of revolution~\cite{bolsinov2000integrable}, Neumann problem~\cite{neumann1859problemate} (i.e., movement on a sphere in a quadratic potential).
One can also consider Kovalevskaya top (a $3$-frequency integrable system), the coordinate corresponding to rotation around vertical axis is cyclic, so if the perturbation does not depend on this coordinate, the problem is reduced to a perturbation of a two-frequency system.
Integrable systems with separatrices often appear as model problems arising after an asymptotic approximation in the study of non-integrable systems, e.g., a rigid body with vibrating suspension point (\cite{markeyev2011equations} and references therein), normal forms near equilibria and periodic trajectories (\cite[\S 8.3, \S 8.4]{arnold2007mathematical} and reference therein), Hamiltonian systems near resonances (e.g.,~\cite{sidorenko2014quasi} and references therein).

An important application of separatrix crossing is multiturn extraction in accelerator physics~\cite{cappi2002novel}.
Consider a particle beam moving in a circular accelerator in horizontal direction (vertical direction is ignored in this problem). Accelerator tune measures the number of oscillations the beam makes on each pass around the accelerator. The idea of multiturn extraction is to vary the tune near a resonant value (e.g., $1:4$), this generates several well-separated beams of particles.
The dynamics is represented by iterations of one-turn transfer map, some power of this map ($4$th power for $1:4$ resonance) is close to identity. This power of the transfer map can be written as unit-time flow of a vector field that slowly depends on time (as the map itself slowly depends on time) with small perturbation depending on time fastly and periodically. Initially, phase portrait of the unperturbed vector field has no separatices, but then, as the tune slowly changes, saddles connected by separatrices appear near the origin and domains bounded by these separatrices begin to grow, capturing the initial beam of particles into different domains bounded by separatrices and splitting it into several smaller beams. Thus the splitting of particle beam can be modeled by separatrix crossings for time-periodic perturbations of one-frequency systems. See also~\cite{bazzani2014analysis}, where possibility of extending the results of adiabatic theory from differential equations to quasi-integrable area-preserving maps is discussed.

The structure of this paper is as follows.
In Section~\ref{s:prelim} we discuss averaging method and two main obstacles to its use, resonances and separatices. Then we discuss results of this paper in a less formal manner.
In Section~\ref{s:scheme-proof} we briefly discuss main ideas of the proofs and differences between resonance crossing far from separatices and near separatices.
Then in Section~\ref{s:results} the results are stated. The rest of the paper contains proofs, plan of these parts can be found in Section~\ref{s:proof_plans}.

\section{Averaging, resonances, separatrix crossings} \label{s:prelim}
\subsection{Averaging method}
Consider a Hamiltonian system
\begin{equation} \label{e:i-u-p}
  \dot q = \pdv{H}{p}, \qquad \dot p = - \pdv{H}{q}, \qquad \dot z = 0.
\end{equation}
Here $p, q \in \mathbb R^m$ and the Hamiltonian $H(p, q, z)$ depends on a scalar or vector parameter $z$. We will call~\eqref{e:i-u-p} the \emph{unperturbed system}.
Suppose that this system is completely integrable (this always holds if $m=1$) and in some domain of the phase space one can introduce action-angle variables $I \in \mathbb R^m$, $\varphi \in \mathbb T^m = \mathbb R^m / 2\pi \mathbb Z^m$. Then~\eqref{e:i-u-p} rewrites as
\begin{equation}
  \dot{I} = 0, \; \dot{\varphi} = \omega(I, z), \; \dot z = 0,
\end{equation}
where $\omega = \pdv{H}{I}$ is the vector of frequencies.
We will call the system~\eqref{e:i-u-p} $m$-frequency system.
Let us add a small perturbation $\varepsilon f$:
\begin{equation}
  \dot{q} = \pdv{H}{p} + \varepsilon f_q(p, q, z, \varepsilon), \qquad
  \dot{p} = - \pdv{H}{q} + \varepsilon f_p(p, q, z, \varepsilon), \qquad
  \dot{z} = \varepsilon f_z(p, q, z, \varepsilon).
\end{equation}
This rewrites in the action-angle variables as
\begin{equation} \label{e:avg-perturbed-energy-angle}
  \dot{I} = \varepsilon f_I(I, \varphi, z, \varepsilon), \qquad
  \dot{\varphi} = \omega(I, z) + \varepsilon f_\varphi(I, \varphi, z, \varepsilon), \qquad
  \dot{z} = \varepsilon f_z(I, \varphi, z, \varepsilon),
\end{equation}
where $f_I, f_\varphi, f_z$ are the components of $f$ in the action-angle variables:
\begin{equation}
  f_y = f_q \pdv{y}{q} + f_p \pdv{y}{p} + f_z \pdv{y}{z}, \qquad y=I, \varphi.
\end{equation}
We see that $I$ and $z$ are \emph{slow} variables of the perturbed system ($\dot I, \dot z = O(\varepsilon)$ far from separatrices), $\varphi$ is \emph{fast} variable.
Evolution of slow variables can be approximately described using averaged system
\begin{equation}
  \dot I = \varepsilon \langle f_I(I, \varphi, z, 0) \rangle_\varphi, \qquad
  \dot z = \varepsilon \langle f_z(I, \varphi, z, 0) \rangle_\varphi.
\end{equation}
Here $\langle \cdot \rangle_\varphi$ denotes averaging over the angle variables $\varphi$.
For perturbations of one-frequency systems far from separatices this works for all initial data with accuracy $O(\varepsilon)$ for times $\sim \varepsilon^{-1}$ (\cite{fatou1928mouvement},~\cite{bogolyubov61}).

There are two major obstacles to the use of averaging. First, when the number of frequencies is at least two, \emph{resonances} between frequencies of unperturbed system are possible, then the values of $\varphi$ for solutions of unperturbed system do not span the whole $\mathbb T^m$. Second, solutions of perturbed system can \emph{cross separatrices} of unperturbed system and move from one domain foliated by Liouville tori $I = const$ to another such domain. One action-angle chart cannot cover such trajectories; moreover, action-angle variables are singular on separatrices. In the following two subsections we discuss each of these obstacles in more detail.

\subsection{Separatrix crossing in one-frequency systems} \label{ss:sep-1-freq}
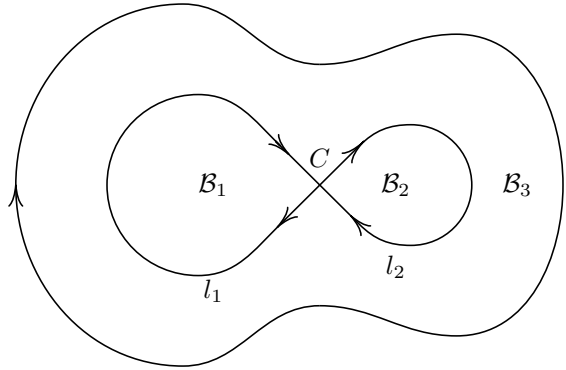
\begin{figure}[H]
  \centering
  \begin{tikzpicture}[scale=0.4]
    \draw[-{To[length=3mm,width=2mm]},semithick](1.1, -1.1)--(1,-1);
    \draw[-{To[length=3mm,width=2mm]},semithick](-1.4, -1.4)--(-1.41,-1.41);
    \draw[-{To[length=3mm,width=2mm]},semithick](-1.1, 1.1)--(-1, 1);
    \draw[-{To[length=3mm,width=2mm]},semithick](1.4, 1.4)--(1.41, 1.41);
    \draw[-{To[length=3mm,width=2mm]},semithick][semithick] (0, 0)
      to [out=45, in=45] (1.3, 1.3) to[out=45, in=180] (3, 2) to[out=0, in=90] (5, 0)
      to[out=-90, in=0] (3, -2) to[out=180, in=-45] (1.3, -1.3) to[out=135, in=135] cycle;
    \draw[semithick] (0, 0) to[out=135, in=135] (-1.3, 1.3)
      to[out=135, in=0] (-4, 3) to[out=180, in=90] (-7, 0)
      to[out=-90, in=180] (-4, -3) to[out=0, in=-135] (-1.3, -1.3) to[out=45, in=45] cycle;

    \draw[semithick] (0, 4) to[out=0, in=180] (4.5, 5) to[out=0, in=90] (8, 0)
      to[out=-90, in=0] (4.5, -5) to[out=180, in=0] (0, -4);
    \draw[semithick] (0, 4) to[out=180, in=0] (-4.5, 6) to[out=180, in=90] (-10, 0);
    \draw[-{To[length=3mm,width=2mm]},semithick] (-4.5, -6) to[out=180, in=-90] (-10, 0.03);
    \draw[semithick] (-4.5, -6) to[out=0, in=180] (0, -4);

    \draw (0, 0.3) node[above]{$C$};
    \draw (-3.5, 0) node{$\mathcal B_1$};
    \draw (-3.5, -3.5) node{$l_1$};
    \draw (2.5, 0) node{$\mathcal B_2$};
    \draw (2.5, -2.7) node{$l_2$};
    \draw (6.5, 0) node{$\mathcal B_3$};
  \end{tikzpicture}
  \caption{The unperturbed system.}
  \label{f:up}
\end{figure}

Consider one-frequency systems, i.e., in~\eqref{e:i-u-p} $p, q \in \mathbb R$. Suppose that for all $z$ the unperturbed system has a saddle $C(z)$ with two separatrix loops $l_1$ and $l_2$ forming a figure eight (Figure~\ref{f:up}).
Solutions of perturbed system can cross separatrices of the unperturbed system.
Set $h(p, q, z) = H(p, q, z) - H(p_C(z), q_C(z), z)$, where $(p_C, q_C) = C(z)$ is the saddle.
Then $h = 0$ on separatrices, assume $h > 0$ in the domain $\mathcal B_3$ (outside separatices) and $h<0$ in the domain $\mathcal B_1 \cup \mathcal B_2$ (inside).
We can use energy $h$ instead of action $I$, then $h, z$ are new slow variables (let us call them \emph{energy-angle variables}). We can write averaged system in this varables.
\begin{figure}[h]
\centering
\includegraphics[width=0.5\textwidth]{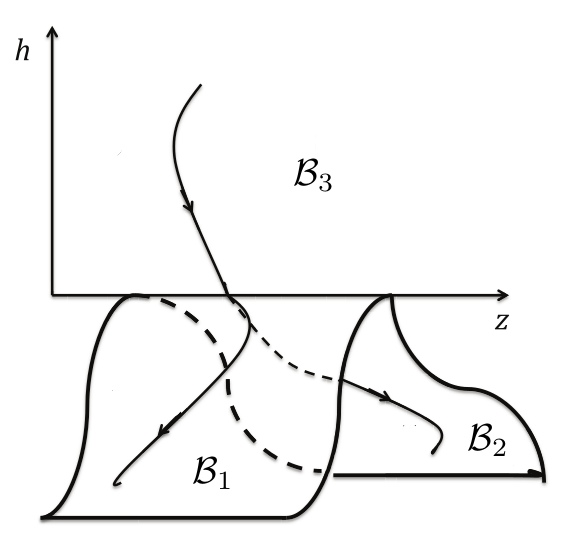}
\caption{Phase space of averaged system. \\
Adapted from~\cite{neishtadt17}.} %
\label{f:avg-phase-space}
\end{figure}
Gluing together averaged systems in $\mathcal B_3$ and $\mathcal B_1$ (or $\mathcal B_2$) by $h=0$ (cf. Figure~\ref{f:avg-phase-space}), we obtain averaged system describing transition from $\mathcal B_3$ to $\mathcal B_1$ (or $\mathcal B_2$).

Averaging method works (\cite{neishtadt17}) for most initial data with measure of exceptional set $O(\varepsilon^r)$, where $r$ can be as large as needed.
Accuracy of averaging method $O(\varepsilon)$ before separatrix crossing and $O(\varepsilon |\ln \varepsilon|)$ after separatrix crossing holds for initial data outside the exceptional set (again for times $\sim \varepsilon^{-1})$. This means that after separatrix crossing evolution of slow variables is approximately described by averaging system describing either transition to $\mathcal B_1$ or to $\mathcal B_2$.
Let us note that the exceptional set is formed by points passing very close to the saddle of perturbed system and the number $r$ above can be set as large as needed, but larger $r$ gives worse constant in the estimate $O(\varepsilon |\ln \varepsilon|)$ for accuracy of averaging method.

Let us say that the outcome of separatrix crossing for some initial data in $\mathcal B_3$ is $1$ if the corresponding solution of perturbed system moves to $\mathcal B_1$ and $2$ if it moves to $\mathcal B_2$.
Initial data in $\mathcal B_3$ with different outcomes are finely mixed, $O(\varepsilon)$ change in initial data is enough to change the outcome. Thus outcome of separatrix croissing is often treated in literature as "random event" with some "probability". We state one precise definition of such probability in Section~\ref{s:main-theorem} below, see also~\cite{neishtadt17} for another definion of such probability and more discussion of this topic.

Particular case when perturbed system is also Hamiltonian (e.g., slow-fast Hamiltonian systems can be written in form~\eqref{e:i-u-p} if slow variables are treated as parameter $z$) is very important and frequently encounered in applications. For this case the action $I$ (in one-frequency case that we consider $2\pi I$ is the area bounded by closed trajectory of unperturbed system) remains constant along solutions of averaged system, so it is called \emph{adiabatic invariant}. Separatrix crossings are still possible, as the area bounded by separatrices changes for perturbed system due to the change of $z$. Separatrix crossing leads to a jump of adiabatic invariant. This jump has magnitude $O(\varepsilon |\ln \varepsilon|)$, there are formulas~\cite{timofeev1978constancy, cary1986adiabatic, neishtadt86, neishtadt1987change} for the value of this jump.

\subsection{Two-frequency systems far from separatrices}
Starting with two-frequency systems, resonances between the frequencies are possible.
Recall that the evolution of fast variables (for two-frequency systems) is given by
\[
  \dot \varphi_1 = \omega_1(I, z), \qquad \dot \varphi_2 = \omega_2(I, z),
\]
where $I = (I_1, I_2)$.
Resonances are given by
\[
 \omega_2/\omega_1 = s_2/s_1, \qquad s_1 \in \mathbb Z_{>0}, \; s_2 \in \mathbb Z.
\]
For each $(s_1, s_2)$ such relation holds on a subset of the space of slow variables $(I_1, I_2, z)$ that is called \emph{resonant surface}. When slow variables are on a resonant surface, the evolution of fast variables for unperturbed system does not span the whole $\mathbb T^2$, but a one-dimensional curve on $\mathbb T^2$. This leads to deviations of the evolution of slow variables from trajectories of averaged system.
\begin{figure}[h]
    \centering
    \begin{minipage}{0.5\textwidth}
        \centering
        \includegraphics[width=0.8\textwidth]{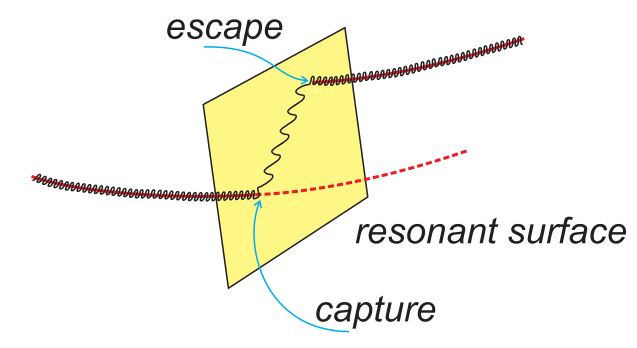} %
        \caption{Capture into a resonance.}
    \end{minipage}\hfill
    \begin{minipage}{0.5\textwidth}
        \centering
        \includegraphics[width=0.8\textwidth]{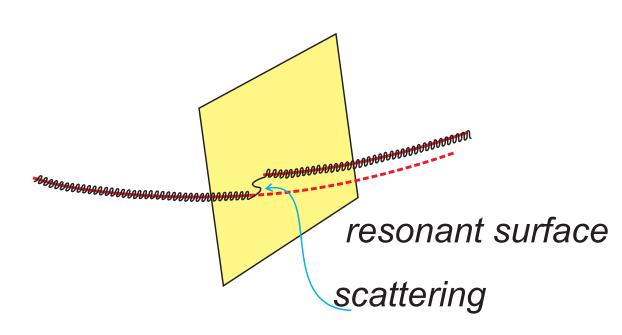} %
        \caption{Scattering on a resonance.}
    \end{minipage}
    \centering{Reproduced from~\cite{neishtadt2019mechanisms}.} %
\end{figure}
Most solutions of perturbed system passing through a resonant surface exhibit \emph{scattering on a resonance}: "random" jump of slow variables of magnitude $O(\sqrt{\varepsilon})$.
Some solutions may be \emph{captured into a resonance}: the solution stays near the resonant surface for time $\sim \varepsilon^{-1}$. Such solution can deviate by $~1$ from trajectories of averaged system. However, measure of initial data that can be captured into resonance is small and evolution of most initial data is still approximately described by averaging method.
Under some genericity condition
\begin{itemize}
  \item Accuracy of averaging method $O(\sqrt{\varepsilon} |\ln \varepsilon|)$ holds for most initial data for times $\sim \varepsilon^{-1}$
  \item Exceptional set has measure $O(\sqrt{\varepsilon})$.
\end{itemize}
This is proved in~\cite{neishtadt1975passage}, see also review~\cite{neishtadt2014} and references therein and an earlier work~\cite{arnold65}, where two-frequency systems were studied under a condition that prohibits capture into resonances.

\subsection{Two-frequency systems near separatrices (our case)}
We consider small time-periodic (with period $2\pi$) perturbations of one-frequency systems with separatrices. Unperturbed system is as in Subsection~\ref{ss:sep-1-freq} (cf. Figure~\ref{f:up}), while the perturbation now depends on time:
\begin{equation}
  \dot{q} = \pdv{H}{p} + \varepsilon f_q(p, q, z, t, \varepsilon), \;
  \dot{p} = - \pdv{H}{q} + \varepsilon f_p(p, q, z, t, \varepsilon), \;
  \dot{z} = \varepsilon f_z(p, q, z, t, \varepsilon),
\end{equation}
where $f_p, f_q, f_z$ are $2\pi$-periodic in $t$.
The time $t$ together with the angle $\varphi$ of the unperturbed one-frequency system form two angle variables of the perturbed system.
This is the simplest case where both resonances and separatrix crossing are encountered.
Perturbations of generic two-frequency systems can be reduced to this case, see Section~\ref{s:two-freqs} below.

As for one-frequency systems near separatrices, let us use $h$ instead of $I$, then $h$ and $z$ are slow variables of perturbed system.
Averaged system is given by
\begin{equation}
  \dot h = \varepsilon \langle f_h(h, \varphi, z, 0) \rangle_{\varphi, t}, \qquad
  \dot z = \varepsilon \langle f_z(h, \varphi, z, 0) \rangle_{\varphi, t},
\end{equation}
where $(f_h, f_\varphi, f_z)$ is the perturbation written in variables $(h, \varphi, z)$. As for one-frequency systems, we can glue averaged systems in different domains and obtain averaging systems describing transition from $\mathcal B_3$ to $\mathcal B_1$ and $\mathcal B_2$.

We assume that (full list of assumptions can be found in Section~\ref{s:results})
\begin{itemize}
\item $H$ and $f$ are analytic
\item $h$ decreases along solutions of averaged system (thus we have transitions from $\mathcal B_3$ to $\mathcal B_1$ and $\mathcal B_2$)
\item some genericity condition similar to the one needed for two-frequency systems far from separatrices holds.
\end{itemize}
We show that (see Section~\ref{s:results} for precise statement of results)
\begin{itemize}
  \item Evolution of most initial data in $\mathcal B_3$ is described by averaged system describing transition from $\mathcal B_3$ to $\mathcal B_1$ or to $\mathcal B_2$ with accuracy $O(\sqrt{\varepsilon} |\ln \varepsilon|)$ over times $\sim \varepsilon^{-1}$.
  \item Exceptional set has measure $O(\sqrt{\varepsilon} |\ln^5 \varepsilon|)$.
  \item Formulas for probabilities of capture in $\mathcal B_1$ and $\mathcal B_2$ similar to one-frequency case hold.
  \item Resonances near separatrices have smaller effect on the dynamics. Consider a part of a resonance surface with $h \sim h_0$, then it can capture measure $O(\sqrt{\varepsilon h_0})$ (up to some power of $\ln h_0$) and the size of scattering on such resonance for trajectories not captured is $O(\sqrt{\varepsilon h_0}|\ln \varepsilon|)$ (up to some power of $\ln h_0$).
\end{itemize}

\section{Statement of results} \label{s:results}

\subsection{Our setting} \label{ss:setting}
Consider a Hamiltonian system with one degree of freedom
\begin{equation} \label{e:unperturbed}
  \dot{q} = \pdv{H}{p}, \qquad \dot{p} = - \pdv{H}{q},
\end{equation}
where $p, q \in \mathbb R$ and the Hamiltonian $H(p, q, z)$ depends on a vector parameter $z = (z_1, \dots, z_n) \in \mathbb R^n$.
We will call this system the \emph{unperturbed system}.
Denote by $\mathcal Z_0 \subset \mathbb R^n$ an open set of parameters, we will only consider $z \in \mathcal Z_0$.
We assume that for all $z \in \mathcal Z_0$ the Hamiltonian $H$ has a saddle $C(z)$ with two separatrix loops $l_1(z)$ and $l_2(z)$ forming a figure eight (Figure~\ref{f:up}) and $C$ is a non-degenerate critical point of $H$.
Denote by $\mathcal B \subset \mathbb R^2_{p, q} \times \mathcal Z_0$ some open neighborhood of $\bigcup_{z \in \mathcal Z_0} l_1(z) \cup l_2(z)$.

The separatrices $\bigcup_{z \in \mathcal Z_0} l_1(z) \cup l_2(z)$ cut $\mathcal B$ into
three open domains: $\mathcal B_1$ and $\mathcal B_2$ inside each separatrix loops and $\mathcal B_3$ outside the union of the separatrices.
Denote
\[
  h(p, q, z) = H(p, q, z) - H(C(z), z).
\]
Then $h = 0$ on the separatrices for all $z$. We will assume
\begin{itemize}
  \item $h > 0$ in $\mathcal B_3$ and $h < 0$ in $\mathcal B_1 \cup \mathcal B_2$ (if the sign is opposite, one can change the sign by exchanging $p$ and $q$)
  \item $\mathcal B_i$ are foliated by the level sets of~$H$. This means that, with a slight abuse of notation, we may write $(h, z) \in \mathcal B_i$
  \item $H$ is analytic on some \emph{compact} set $\tilde{\mathcal B}$ with $\overline{\mathcal B} \subset \Int \tilde{\mathcal B}$ and $\tilde{\mathcal B}$ is also foliated by the level sets of~$H$.
\end{itemize}

Take a compact $\mathcal Z \subset \mathcal Z_0$.
Then for small enough $c_h, c_z > 0$ for any $z_* \in \mathcal Z$ for any $h, z$ with $|h| < c_h$, $\norm{z - z_*} < c_z$ we have $(h, z) \in \mathcal B$.

\noindent Consider the perturbed system
\begin{align}
\begin{split} \label{e:perturbed-pq}
  \dot q &= \pdv{H}{p} + \varepsilon f_q(p, q, z, \lambda, \varepsilon), \\
  \dot p &= - \pdv{H}{q} + \varepsilon f_p(p, q, z, \lambda, \varepsilon), \\
  \dot z &= \varepsilon f_z(p, q, z, \lambda, \varepsilon), \\
  \dot \lambda &= 1.
\end{split}
\end{align}
Here
\begin{itemize}
  \item The perturbation depends on the time $\lambda$. We used the notation $t$ before, $\lambda$ is introduced to distinguish the time $\lambda$ the perturbation depends on from the time $t$ of the unperturbed system that will be used below to parametrize trajectories of this system
  \item $f_p, f_q, f_z$ are $2\pi$-periodic in $\lambda$
  \item $f_p, f_q, f_z$
  are $C^2$ in $\tilde{\mathcal B} \times \mathbb R_\lambda \times [-\varepsilon_0, \varepsilon_0]$ for some $\varepsilon_0 > 0$
  \item for $\varepsilon=0$ the functions $f_p, f_q, f_z$ are real-analytic in $\tilde{\mathcal B} \times \mathbb R_\lambda$.
\end{itemize}

\subsection{Averaged system} \label{s:avg}
The rate of change of $h$ along the solutions of the perturbed system is $\varepsilon f_h$, where
\begin{equation}
  f_h = \pdv{h}{p} f_p + \pdv{h}{q} f_q + \pdv{h}{z} f_z.
\end{equation}
The variables $h$ and $z$ are slow variables of the perturbed system, their evolution can be approximately tracked using averaged system.
In each domain $\mathcal B_i$ denote by $T(h, z)$ and $\omega(h, z)$ the period and the frequency of the solution of the unperturbed system~\eqref{e:unperturbed} with given $h$, $z$.
The averaged system is given by the equations
\begin{equation} \label{e:avg-h-z}
  \dot{\overline h} = \varepsilon f_{h, 0}(\overline h, \overline z), \qquad
  \dot{\overline z} = \varepsilon f_{z, 0}(\overline h, \overline z),
\end{equation}
where
\begin{equation}
  f_{h, 0} = \frac{1}{2\pi T} \int_0^{2\pi} \oint f_h|_{\varepsilon=0} dt d\lambda, \qquad
  f_{z, 0} = \frac{1}{2\pi T} \int_0^{2\pi} \oint f_z|_{\varepsilon=0} dt d\lambda
\end{equation}
denote averages of $f_h$ and $f_z$ over the angle variables $\varphi, \lambda$.
The inner integrals above are taken along the closed trajectory of the unperturbed system given by $h = \overline h$, $z = \overline z$ (inside $\mathcal B_i$) and this trajectory is parametrized by the time $t$ of the unperturbed system.
Recall that $l_1$ and $l_2$ denote the separatrices of the unperturbed system.
Denote
\begin{equation} \label{e:Theta-3}
  \Theta_i(z) = -\frac{1}{2\pi} \int_{0}^{2\pi} \oint_{l_i} f_h|_{\varepsilon=0}(p(t), q(t), z, \lambda) dt d\lambda
  \qquad \text{for} \; i=1, 2; \qquad \Theta_3 = \Theta_1 + \Theta_2.
\end{equation}
(the separatrices are parametrized by the time $t$ of unperturbed system). These integrals converge, see~\cite[Section~2.2]{neishtadt17}.
We will assume $\Theta_i > 0$, $i=1,2,3$. Note that near separatrices we have $T f_{h, 0} \approx -\Theta_i$ in $\mathcal B_i$.
Thus near separatrices \emph{$h$ decreases} in all $\mathcal B_i$ along the solutions of averaged system.
Moreover, for some $K_h > 0$ for small enough $|h| > 0$ we have
\begin{equation} \label{e:K-h}
  f_{h, 0} < -K_h |\ln^{-1} h| < 0.
\end{equation}
Once solution of the averaged system in $\mathcal B_3$ reaches $h=0$, one can continue this solution using the averaged system in $\mathcal B_1$ or $\mathcal B_2$ (cf. Figure~\ref{f:avg-phase-space}). We will say that such "glued" solutions correspond to capture into $\mathcal B_1$ or $\mathcal B_2$, respectively. This is discussed in more detail in~\cite[Section 2.3]{neishtadt17}.

\subsection{Condition $B'$} \label{s:cond-B}
Let $\varphi$ be the angle variable (from the pair of action-angle variables) of the unperturbed system defined in $\mathcal B_3$.  Pick the transversal (to the solutions of unperturbed system) $\varphi = 0$ so that for all $z$ it is a smooth curve that crosses $l_1$ at some point $a_1(z) \ne C(z)$. It is easy to check that the transversal $\varphi = \pi$ crosses $l_2$ at some point that we denote by $a_2(z)$.
Let us define the coordinates $t_1$, $t_2$ on $l_1$ and $l_2$ as the time (for the unperturbed system) passed after the point $a_1$ and $a_2$, respectively.

Denote
\begin{equation}
  M_i(Q) = \int_{l_i} f_h(h{=}0, z, t_i{=}t, \lambda{=}t {-} Q, \varepsilon{=}0) dt \qquad \text{for } i=1,2.
\end{equation}
This is the famous Melnikov function~\cite{mel1963stability} used to describe separatrix splitting, it is $2\pi$-periodic.
Given $s = (s_1, s_2) \in \mathbb Z_{>0}^2$, set
\begin{align}
\begin{split}
  F^*_{s, i}(Q, z) = &-(2\pi)^{-1} \Big\langle
    M_i \Big(Q - 2\pi \frac{j}{s_2}\Big)
  \Big\rangle_{j=0, \dots, s_2-1} \qquad \text{for } i=1,2; \\
  F^*_{s, 3}(Q, z) = & F^*_{s, 1}(Q, z) + F^*_{s, 2}\Big(Q + 2\pi \frac{\{s_1/2\}}{s_2}, z\Big).
\end{split}
\end{align}
Here $\{\cdot\}$ denotes the fractional part and $\langle \psi(j) \rangle_{j=0, \dots, k-1}$ denotes the average $\frac{\psi(0) + \dots + \psi(k-1)}{k}$. The functions $F^*_{s, i}(Q)$ are periodic with period $2\pi/s_2$.

For a function $F(Q, z)$ set $V_F(Q, z) = \int_0^Q F(\tilde Q, z) d \tilde Q$.
We will need the following

\noindent \textbf{Condition $\mathbf{B'(V)}$}. All extrema of the function $V(Q)$ are non-degenerate (i.e. for all $q$ such that $\pdv{V}{Q} = 0$ we have $\pdv[2]{V}{Q} \ne 0$).
Moreover, at different local maxima of $V$ the values of $V$ are different.

\noindent \textbf{Condition $\mathbf{B'(s, z, i)}$}. The function $V(Q) = V_{F^*_{s, i}}(Q, z)$ satisfies the condition $B'(V)$ above.

\noindent For fixed $s$ this is a codimension one genericity condition on $(f_p, f_q, f_z)$  and $z$.

\noindent
\textbf{Condition $\mathbf{B'(z, i)}$}.
For any $s = (s_1, s_2)$ condition~$B'(s, z, i)$ holds.

\noindent The lemma below means that $B'(z, i)$ is also a codimension one genericity condition on $(f_p, f_q, f_z)$  and $z$, as $V_F$ has no extrema if $F > 0$ and thus satisfies condition $B'(V)$.
\begin{lemma}
  Given a uniform bound
  \[
    \norm{f_p}_{C^2}, \norm{f_q}_{C^2}, \norm{f_z}_{C^2} \le K,
    \qquad
    \frac{\Theta_i}{2\pi} > K^{-1} \qquad \text{for } i=1,2,3
  \]
  ($K > 0$), there exists $S(K)$ such that $F^*_{s, i} > K^{-1}/2$ if $s_2 > S(K)$.
\end{lemma}
\begin{proof}
  We have $f_h(C(z), z) = 0$ by~\cite[Lemma 2.1]{neishtadt17}.
  For definiteness, consider the separatrix $l_1$;
  let $(p(t_1), q(t_1))$ denote the point on $l_1$ with given value of $t_1$.
  As $(p(t_1), q(t_1))$ exponentially converges to $C$ for $t_1 \to +\infty$ and $t_1 \to -\infty$, there exists $T_1$ such that
  \begin{equation} \label{e:tmp-4891}
    \int_{t_1 = -\infty}^{-T_1} \max_{\lambda}|f_h(h{=}0, t_1, \lambda)|
    +
    \int_{t_1 = T_1}^\infty \max_{\lambda}|f_h(h{=}0, t_1, \lambda)|
     < K^{-1}/100.
  \end{equation}
  For $|t_1| < T_1$ (i.e. far from $C$) the transition between coordinates $h, t_1$ and $p, q$ is smooth, so $f_h(h{=}0, t_1, \lambda)$ is a smooth function of $t_1, \lambda$ with bounded $C^2$-norm.
  We have
  \[
    \int_{-T_1}^{T_1} \Big\langle f_h\Big(
    h{=}0, z, t_1, \lambda{=}t_1 {-} Q {+} 2\pi \frac{j}{s_2} \Big) \Big\rangle_{j=0, \dots, s_2-1} dt_1
    = \int_{-T_1}^{T_1} \langle f_h(h{=}0, z, t_1, \lambda) \rangle_\lambda + O(s_2^{-2}),
  \]
  as averaging over $j$ gives trapezoidal rule approximation for averaging over $\lambda$.
  The integral on the right-hand side is approximately $-\Theta_1$ with error at most $K^{-1}/100$ by~\eqref{e:tmp-4891}. Take $S$ such that the $O(s_2^{-2})$ term is less then $K^{-1}/100$, then for $s$ with $s_2 > S$ we have $|2\pi F^*_{s, 1} - \Theta_1| < 3K^{-1}/100$. We can obtain similar estimate for $|2\pi F^*_{s, 2} - \Theta_2|$, together they yield the estimate on $|2\pi F^*_{s, 3} - \Theta_3|$.
\end{proof}

\subsection{Main results} \label{s:main-theorem}
Denote $\overline X = (h, z)$,
$\mathcal A_i = \mathcal B_i \times [0, 2\pi]_\lambda$,
$\mathcal A = \mathcal B \times [0, 2\pi]_\lambda$.
Take a point
\[
  \hat{\overline X}_0 = (\hat{\overline h}_0, \hat{\overline z}_0) \in \mathcal B_3
\]
and $\Lambda > 0$.

For $i=1,2$ denote by $\hat{\overline X}_i(\lambda)$ the solution of averaged system~\eqref{e:avg-h-z} describing capture from $\mathcal B_3$ to $\mathcal B_i$ with $\hat{\overline X}_i(0) = \hat{\overline X}_0$.
Denote by $\mathcal Z_B \subset \mathcal Z$ the set of $z$ that satisfy condition~$B'(z, i)$ for all $i=1, 2, 3$.

Suppose that, in addition to assumptions from Subsections~\ref{ss:setting} and~\ref{s:avg},
\begin{enumerate}
  \item $\omega$ decreases along solutions of the averaged system: $\pdv{\omega}{z} f_{z, 0} + \pdv{\omega}{h} f_{h, 0} < 0$ in $\mathcal B_1 \cup \mathcal B_2 \cup \mathcal B_3$;
  \item solutions of averaged system stay in $\mathcal B$, i.e., there exists $c_{\mathcal B} > 0$ such that for any $\lambda \in [0, \varepsilon^{-1} \Lambda]$ the points $\hat{\overline X}_i(\lambda)$ are in $\mathcal B$ and at least $c_{\mathcal B}$-far from the border of $\mathcal B$;
  \item solutions of averaged system cross separatices, i.e., $\hat{\overline X}_1(\varepsilon^{-1}\Lambda) \in \mathcal B_1$,
  $\hat{\overline X}_2(\varepsilon^{-1}\Lambda) \in \mathcal B_2$.
  \item separatrix crossing happens when $z \in \mathcal Z_B$. To state this more precisely, denote by $\lambda_*$ the time when $\hat{\overline X}_i$ cross separatrices (i.e. $\hat{\overline h}_i(\lambda_*) = 0$) and set $z_* = \hat{\overline z}_i(\lambda_*)$ (note that $\lambda_*$ and $z_*$ are the same for $i=1,2$). Then $z_* \in \mathcal Z_B$;
  \item For some small
  enough\footnote{It should be so small that if $\lambda$ is $\varepsilon^{-1} \Lambda_1$-close to $\lambda_*$, we have $|\overline h_i(\lambda)| < c_h/2$ and $\norm{\overline z_i(\lambda) - z_*} < c_z/2$, where $c_h, c_z$ should be so small that we can apply Theorem~\ref{t:sep-pass} below.}
  $\Lambda_1$ for $i=1,2$ the solutions $\hat{\overline X}_i$ satisfy condition $B$ from~\cite[Section~2]{neishtadt2014}
  for\footnote{In the notation of~\cite{neishtadt2014} these solutions should not cross the set $\theta$.}
  $\lambda \in [0, \lambda_* - \varepsilon^{-1} \Lambda_1]$
  and
  $\lambda \in [\lambda_* + \varepsilon^{-1} \Lambda_1, \varepsilon^{-1} \Lambda]$. This is a genericity condition, as explained in~\cite{neishtadt2014}.
\end{enumerate}

\noindent Let $B_r(\overline X) \subset \mathbb R^{n+1}_{h, z}$ denote the open ball with center $\overline X$ and radius $r$.
Denote
\[
  A_r(\overline X) = B_r(\overline X) \times [0, 2\pi]^2_{\lambda, \varphi} \subset \mathcal A.
\]
Denote by $m$ the Lebesgue measure on $\mathbb R^{n+2}_{p, q, z} \times [0, 2\pi]_\lambda$.
\begin{theorem} \label{t:main}
  There exist $C, r > 0$ such that for any small enough $\varepsilon$ there exists
  \[
    \mathcal E \subset A_r(\hat{\overline X}_0) \subset \mathcal A_3
    \qquad \text{with} \qquad
    m(\mathcal E) \le C \sqrt{\varepsilon} |\ln^5 \varepsilon|
  \]
  such that the following holds for any $X_0 \in A_r(\hat{\overline X}_0) \setminus \mathcal E$.

  Set $\lambda_0 = \lambda(X_0)$.
  Denote by $X(\lambda)$ the solution of the perturbed system~\eqref{e:perturbed-pq} with $X(\lambda_0) = X_0$.
  Then $X(\lambda_0 + \varepsilon^{-1} \Lambda) \in \mathcal A_i$ for some $i=1,2$.
  Denote by $\overline X(\lambda)$ the solution of the averaged system~\eqref{e:avg-h-z} describing transition from $\mathcal B_3$ to $\mathcal B_i$ with
  $\overline X(\lambda_0) = (h(X_0), z(X_0))$.
  Then for any $\lambda \in [\lambda_0, \lambda_0 + \varepsilon^{-1} \Lambda]$ we have
  \begin{equation} \label{e:mt-close}
    |h(X(\lambda)) - h(\overline X(\lambda))| < C \sqrt{\varepsilon} |\ln \varepsilon|, \qquad
    \norm{z(X(\lambda)) - z(\overline X(\lambda))} < C \sqrt{\varepsilon} |\ln \varepsilon|.
  \end{equation}
\end{theorem}
\noindent This theorem is proved in Section~\ref{s:decompose}, it is reduced to technical Theorem~\ref{t:sep-pass} below on crossing a small neighborhood of separatrices.

Let us now discuss "probabilities" of capture in $\mathcal A_1$ and $\mathcal A_2$. Given $X_0 = (p_0, q_0, z_0, \lambda_0) \in \mathcal A_3$ and small $\delta > 0$, denote $I_0 = I(p_0, q_0, z_0)$, $\varphi_0 = \varphi(p_0, q_0, z_0)$, $h_0 = h(p_0, q_0, z_0)$ (here $I, \varphi$ are action-angle variables of unperturbed system in $\mathcal A_3$).
Let us define the set $U^\delta \subset \mathcal A_3$ by
\begin{equation}
  U^\delta = \{ (I, z, \varphi, \lambda):
  |I-I_0|, \norm{z-z_0}, |\varphi-\varphi_0|, |\lambda-\lambda_0| < \delta \}.
\end{equation}
Solutions of the perturbed system with initial data in $U^\delta \setminus \mathcal E$ are described by solutions of averaged system describing transition to $\mathcal B_1$ (we say that such initial data is captured in $\mathcal A_1$) or transition to $\mathcal B_2$ (we say that such initial data is captured in $\mathcal A_2$).
Denote by
$U^\delta_1, U^\delta_2 \subset U^\delta \setminus \mathcal E$
the sets of initial data captured in $\mathcal A_1$ and $\mathcal A_2$, respectively.
\begin{definition}[V.I. Arnold, \cite{arnold63}]
  The probability of capture in $\mathcal B_i$, $i=1,2$ is
  \begin{equation}
    P_i(X_0) = \lim_{\delta \to 0} \lim_{\varepsilon \to 0} \frac{m(U^\delta_i)}{m(U^\delta)}.
  \end{equation}
\end{definition}

\begin{proposition} \label{p:prob}
  \begin{equation}
    P_i(X_0) = \Theta_i(z_*)/\Theta_3(z_*), \qquad i=1,2.
  \end{equation}
  Here $z_*$ is the value of $z$ when the solution of averaged system with initial data $(h_0, z_0)$ crosses separatrices.
\end{proposition}
\noindent This proposition is proved in Section~\ref{s:prob}.

\begin{remark} \label{r:res-effects-decrease}
  Resonances near separatrices have smaller effect on the dynamics. Consider a part of a resonance surface with $h \sim h_0$, then it can capture measure $O(\sqrt{\varepsilon h_0})$ (up to some power of $\ln h_0$) and the size of scattering on such resonance for trajectories not captured is $O(\sqrt{\varepsilon h_0}|\ln \varepsilon|)$ (up to some power of $\ln h_0$).
\end{remark}
\noindent This informal statement is stated more precisely in Section~\ref{s:est-res-cross} as Remark~\ref{r:res-effects-decrease-formal} and proved there.

\begin{remark}
  Our results are applicable for Hamiltonian systems with two and a half degrees of freedom, i.e., with the Hamiltonian
  \[
    H(p, q, y, x, t) = H_0(p, q, y, x) + \varepsilon H_1(p, q, y, x, t),
  \]
  where $H_1$ is $2\pi$-periodic in time $t$.
  Here pairs of conjugate variables are $(q, p)$ and $(\varepsilon^{-1} x, y)$; $(q, p)$ are fast variables and $(x, y)$ are slow variables.
  Indeed, we can take $z = (x, y)$, then Hamilton's equations will be of form~\eqref{e:perturbed-pq}.
  Then the action $I(h, z)$ is the \emph{adiabatic invariant}, it stays constant along the solutions of averaged system. Theorem~\ref{t:main} shows that for most initial data $I$ is preserved after crossing separatices with accuracy\footnote{Far from separatrices this follows from the statement of Theorem~\ref{t:main}, as $\pdv{I}{h} = \omega^{-1}$ is bounded. But this also holds near separatrices, as in the proof of this theorem the difference in $I$ is estimated, cf. Section~\ref{s:approach-proof}.} $\sim \sqrt{\varepsilon} |\ln \varepsilon|$.

  For this case the values $\Theta_i$, $i=1,2$ can be computed as
  $\Theta_i(y, x) = \{ S_i, H(C) \}$ (\cite{arnold2007mathematical}),
  where $S_i$ is the area of the domain $\mathcal B_i$ (cf. Figure~\ref{f:up}), $H(C) = H(p_C, q_C, y, x)$ denotes the value of the Hamiltonian $H$ at the saddle $C(y, x) = (p_C(y, x), q_C(y, x))$, and $\{\;,\;\}$ is the Poisson bracket in the variables $y$, $x$.
\end{remark}

In~\cite{arnold65} two-frequency systems (far from separatrices) were considered under the following condition prohibiting capture into resonances: $\dot \omega < C < 0$ along the solutions of \emph{perturbed system}. This condition cannot hold near the separatrices, as $\omega$ is undefined on the separatrices, but a similar condition can be that $f_h < 0$ except in the saddle, where $f_h = 0$.
\begin{remark} \label{r:no-capture}
  Suppose $f_h(p, q, z, \lambda)|_{\varepsilon=0} < 0$ for all $\lambda$ if $(p, q) \ne C(z)$. Then Theorem~\ref{t:main} holds with stronger estimates: for any $r>0$ we can take exceptional set with measure $O(\varepsilon^r)$, then accuracy of averaging method (i.e., rhs in~\eqref{e:mt-close}) $O(\sqrt{\varepsilon})$ holds for initial data outside exceptional set.
\end{remark}
\noindent A sketch of proof of this remark is given in Section~\ref{s:no-capture-proof}.
Note that the exceptional set in this remark is formed by initial data near separatrices of the saddle of \emph{perturbed} system (these separatrices wind around the figure eight and are present even far from separatrices of unperturbed system).
Estimate for accuracy of averaging method in Remark~\ref{r:no-capture} cannot be improved, as scattering on resonance of amplitude $\sim \sqrt{\varepsilon}$ is possible.

\subsection{Two-frequency systems} \label{s:two-freqs}
Consider an integrable two-frequency system
\begin{equation}
  \dot q_1 = \pdv{H}{p_1}, \qquad
  \dot q_2 = \pdv{H}{p_2}, \qquad
  \dot p_1 = -\pdv{H}{q_1}, \qquad
  \dot p_2 = -\pdv{H}{q_2}, \qquad
  \dot z = 0
\end{equation}
with Hamiltonian $H$ (depending on a parameter $z$) and another first integral $F$.
Denote by $v$ the Hamiltonian vector field given by $H$.
Separatrices are singularities of the Liouville foliation.
An isoenergy level $\{ H=h \}$ is called \emph{topologically stable}~\cite[\S3.3]{bolsinov2004integrable} if for sufficiently small variations of energy level $\{ H=h+\delta \}$ the Liouville foliations on these isoenergy levels are equivalent (i.e., there exists a diffeomorphism that maps one foliation into another).
Consider perturbed and averaged systems.
By separatrix crossing we mean that solution of averaged system crosses a singular leaf $H=h_0$, $F=f_0$, $z=z_0$ of the Liouville foliation.
Suppose that separatrix crossing for averaged system happens on a topologically stable energy level $H=h_0$, $z=z_0$.
We assume that the restrictions of $F$ on isoenergy surfaces are Bott functions, such singularities are typical in real problems in physics and mechanics  (\cite[\S1.8.1]{bolsinov2004integrable}).

Under these assumptions perturbation of two-frequency integrable system near separatrices can be reduced to a periodic perturbation of one-frequency integrable system depending on an additional parameter (it is denoted by $h$ in~\eqref{e:one-freq-reduced} below). Thus Theorem~\ref{t:main} can be applied to perturbations of two-frequency systems.
\begin{lemma} \label{l:two-freq-one-freq}
  There exist ($z$-dependent) new canonical coordinates $(p, q, h, s)$, $s \in [0, 2\pi]$ and a cover (a bijection or a double cover)
  $C: p, q, h, s, z \mapsto p_1, p_2, q_1, q_2, z$ defined in a neighborhood of the singular leaf $\{ H=h_0, \; F=f_0, \; z=z_0 \}$ such that the dynamics of $(p, q, h, s, z)$ given by the unperturbed system with $s$ as new time is
  \begin{equation} \label{e:one-freq-reduced}
    s' = 1, \qquad
    h' = 0, \qquad
    p' = - \pdv{S}{q}, \qquad
    q' = \pdv{S}{p}, \qquad
    z' = 0
  \end{equation}
  with some Hamiltonian $S(p, q, h, z)$ depending on parameters $h$ and $z$, here $\cdot'$ denotes the derivative with respect to $s$.
\end{lemma}
\noindent This lemma is based on the description of two-frequency integrable systems in~\cite{bolsinov2004integrable}, its proof (by A.V. Bolsinov) is presented in Appendix~\ref{a:two-freq-reduction}. In this appendix we also write perturbed system in the new coordinates, it has the form~\eqref{e:perturbed-pq}.

\section{Scheme of proof} \label{s:scheme-proof}
\subsection{Far from separatrices}
\begin{figure}[h]
  \begin{center}
  \includegraphics[scale=0.4]{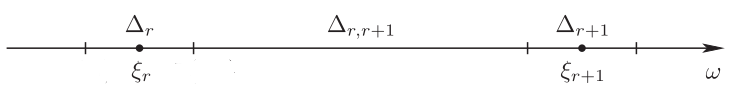}
  \caption{Resonant and non-resonant zones. \\ Adapted from~\cite{neishtadt2014}. } %
  \label{f:res-zones}
  \end{center}
\end{figure}

Let us recall the scheme of proof of the result on averaging in two-frequency systems far from separatrices from~\cite{neishtadt2014}. More details can be found in~\cite[Section~5]{neishtadt2014}.
General two-frequency case is reduced to time-periodic perturbations of one-frequency systems.
Denote by $\omega$ the frequency of unperturbed system.
Effect of a resonance $\omega = s_2/s_1$ is determined by Fourier coefficients $f_{ks}$, $s=(s_1, s_2)$, $k \in \mathbb Z_{>0}$ of the perturbation $f(\varphi, t) = (f_h, f_\varphi, f_z)$.
As Fourier coefiicients of analytic functions decay exponentially, effect of resonances also exponentially decays when $s$ grows.
It is enough to consider $\sim \ln^2 \varepsilon$ resonances with $|s_1|, |s_2| \lesssim |\ln \varepsilon|$, as total effect of other resonances is negligebly small.
Let us represent each of these resonances by a point $s_2/s_1$ on the line of possible values of $\omega$ (cf. Figure~\ref{f:res-zones}) and surround it by a \emph{resonance zone} with width of order
\begin{equation}
  \delta_s = \sqrt{\varepsilon a_s/s_1} + \varepsilon s_1, \qquad a_s = e^{-K |s|}.
\end{equation}
Here $K > 0$ is some constant; $|s| = |s_1| + |s_2|$; the numbers $a_s$ are bounds on Fourier coefficients $f_s$.
Two neighboring resonant zones $\Delta_{r}$ and $\Delta_{r+1}$ given by $\omega \approx \xi_r$ and $\omega \approx \xi_{r+1}$ ($\xi_r, \xi_{r+1} \in \mathbb Q$) are shown in Figure~\ref{f:res-zones} together with a \emph{non-resonant} zone $\Delta_{r, r+1}$ between them.
Total width of all resonant zones is $O(\sqrt{\varepsilon})$.

It is assumed that $\omega$ decreases along solutions of averaged system.
Dynamics between resonant zones is described by standard coordinate change used to justify averaging method. We will discuss methods used to study dynamics inside resonant zones in Subsection~\ref{ss:res-far-sep} below.
It turns out that  there are only finitely many resonances such that capture into resonance is possible (as for capture we should have $a_s \sim 1$) and each such resonant can capture measure $O(\sqrt{\varepsilon})$ with scattering of order $\sqrt{\varepsilon} |\ln \varepsilon|$ for the remaining trajectories.
The size of scattering on other resonances is of the same order as the width of corresponding resonance zone.
Thus total scattering on all resonances is $O(\sqrt{\varepsilon} |\ln \varepsilon|)$.

\subsection{Resonant zones near separatrices}
Near the separatrices the angle variable $\varphi$ behaves badly, as solutions of unperturbed system spend most time near the saddle $C$. Many functions used in the averaging method are unbounded (e.g., $f_\varphi$ in \eqref{e:avg-perturbed-energy-angle}).
Also, for a smooth function $\psi(p, q, z)$ after transition to the energy-angle variables the partial derivative $\pdv{\psi}{h}$ can be unbounded, as it is taken for fixed $\varphi$.
Estimates on these functions are required to use the averaging
method.

Estimates on Fourier coefficients of the perturbation $f$ are needed to determine resonant zones.
Denote $Y = (p, q, z)$.
For fixed $h, z$ we prove that the function $Y(\varphi)$ can be continued to the complex domain
\begin{equation} \label{e:Y-continuation}
  |\Im \varphi| \lesssim \omega \sim |\ln^{-1} h|.
\end{equation}
This gives estimates on Fourier coefficients of $f$ (cf. Lemma~\ref{l:fourier} below). Then we define resonant zones and non-resonant zones between them. The definition is same as far from separatices, but the width of resonant zones is different:
\begin{equation}
  \delta_s \approx
  \sqrt{\varepsilon b_s h^{-1}}
  + \varepsilon h^{-1} \ln^2 \varepsilon,
\end{equation}
where $\approx$ means that this formula holds up to multiplying both summands by some powers of $\ln h$, and $b_s = e^{-C_F s_2}$ with some constant $C_F$ is a bound on norm $\norm{f_s}$ of the Fourier coefficient. Full formula for $\delta_s$ is given below, cf.~\eqref{e:delta}.
The formula above is for the width of resonant zones in $\omega$, the width in $h$ can be computed using $\pdv{\omega}{h} \sim h^{-1} \ln^{-2} h$, it is (up to some powers of $\ln h$)
\begin{equation}
  \delta_{s, h} \approx
  \sqrt{\varepsilon b_s h}
  + \varepsilon \ln^2 \varepsilon.
\end{equation}
We see that the width of resonant zones in the phase space decreases near separatices.
Further differences in the structure of resonant zones compared with what happens far from separatices are as follows.
\begin{itemize}
  \item We consider resonances with $|s_1|, |s_2| \lesssim \ln^2 \varepsilon$, their number is $\sim \ln^4 \varepsilon$.
  \item There are infinitely many resonant zones $\omega = s_2/s_1$ such that capture into resonance is possible (i.e, number of such zones grows when $\varepsilon \to 0$). Capture is only possible when $|s_2|$ is small, but $|s_1|$ can be large.
  \item Total width (in $h$) of all resonant zones is still $O(\sqrt{\varepsilon})$, as far from separatices.
  \item Methods describing passage through resonant zones work when $|h| \gtrsim \varepsilon |\ln^5 \varepsilon|$, passage through the zone $|h| \lesssim \varepsilon |\ln^5 \varepsilon|$ is thus considered separately (we call it \emph{immediate neighborhood of separatrices}). We only consider resonant and non-resonant zones with $|h| \gtrsim \varepsilon |\ln^5 \varepsilon|$.
\end{itemize}

Dynamics in non-resonant zones is studied with help of the standard coordinate change used to justify averaging method, as in~\cite{neishtadt2014}. This coordinate change is written using decomposition of the perturbation $f$ in Fourier series, estimates on Fourier coefficients follow from analytic continuation in the domain~\eqref{e:Y-continuation}. We use estimates on components of the perturbation in energy-angle and their partial derivatives obtained in~\cite{neishtadt2020phase}.

Dynamics in the immediate neighborhood of separatrices might be complicated, as different resonant zones begin to overlap when $h \lesssim \varepsilon \ln^2 \varepsilon$. Resonance overlap is the celebrated Chirikov criterion for chaotic dynamics~\cite{chirikov1960resonance}.
However, this zone is small, and a volume argument based on the fact that the flow of perturbed system changes volume slowly (as perturbation has divergence $O(\varepsilon)$) can be used to show that most solutions leave the zone $|h| < \varepsilon |\ln^5 \varepsilon|$ after time $O(\varepsilon^{-1/2}|\ln \varepsilon|)$ passes, thus leading to $O(\sqrt{\varepsilon}|\ln \varepsilon|)$ deviation from the solution of averaged system. This argument is where exceptional set of measure $O(\sqrt{\varepsilon} |\ln^5 \varepsilon|)$ appears, the power of logarithm comes from technical details of the proof and might potentially be improved (but still exceptional set should have measure at least $O(\sqrt{\varepsilon} |\ln \varepsilon|)$, as estimates far from separatices cannot be improved~\cite{neishtadt2014}).

The hardest part of this paper is the study of dynamics in resonant zones near separatrices.
Let us first recall how resonant zones far from separatices are studied and then comment on the differences arising near separatrices.

\subsection{Passage through resonant zones far from separatrices} \label{ss:res-far-sep}
Dynamics near resonances can be reduced to an auxiliary system.
This approach is widely used (cf. references in introduction of~\cite{neishtadt2014}), our exposition loosely follows~\cite{neishtadt2014} with some parts modified to be closer to the way we treat resonance zones near separatrices in the current paper.
Consider the following simplified case, where parameter $z$ and dependence of the perturbation on $\varepsilon$ are removed.
\begin{equation}
  \dot I = \varepsilon f_I(I, \varphi, t), \qquad
  \dot \varphi = \omega(I) + \varepsilon f_\varphi(I, \varphi, t).
\end{equation}
Fix rational $\hat \omega$, consider resonance $\omega \approx \hat \omega$.
Define $\hat I$ by $\omega(\hat I) = \hat \omega$.
Let us introduce new variables
\begin{equation}
  \gamma = \varphi - \hat \omega t, \qquad
  J = I - \hat I.
\end{equation}
Near the resonance $J, \dot \gamma \approx 0$. Denote $\hat {\pdv{\omega}{I}} = \pdv{\omega}{I}\big|_{\hat I}$. We have
\begin{equation}
  \dot \gamma = \omega - \hat \omega + O(\varepsilon) = \hat {\pdv{\omega}{I}} J + O(\varepsilon) + O(J^2),
  \qquad
  \dot J = \varepsilon f_I(I, \gamma + \hat \omega t, t) =  \varepsilon f_I(\hat I, \gamma + \hat \omega t, t) + O(\varepsilon J).
\end{equation}
Set $\alpha = \sqrt{\varepsilon / \hat {\pdv{\omega}{I}}}$, $\beta = \sqrt{\varepsilon \hat {\pdv{\omega}{I}}} \sim \sqrt{\varepsilon}$, $P = J / \alpha$, $Q = \gamma$. Assume $P \lesssim 1$ (as we consider only what happens near resonance, this shows why width of resonant zones is $O(\sqrt{\varepsilon})$).
We get
\begin{equation}
  \dot Q = \beta P + O(\varepsilon), \qquad
  \dot P = \beta f_I(\hat I, Q + \hat \omega t, t) + O(\varepsilon).
\end{equation}
We see that $t$ is fast variable compared with $P, Q$. Let us apply averaging over $t$, we omit justification of the use of averaging method here (it goes close to the standard justification of averaging using coordinate change) and simply replace the true system by averaged system. Denote $F(Q) = \langle  f_I(\hat I, Q + \hat \omega t, t) \rangle_{t \in [0, 2\pi]}$. We get the system
\begin{equation}
  \dot Q = \beta P + O(\varepsilon), \qquad
  \dot P = \beta F(Q) + O(\varepsilon).
\end{equation}
Taking new time $\tau = \beta t$ and denoting $a' = \dv{a}{\tau}$, we get an auxiliary system describing dynamics near the resonance:
\begin{equation} \label{e:aux-8432}
  Q' = P + O(\sqrt{\varepsilon}), \qquad
  P' = F(Q) + O(\sqrt{\varepsilon}).
\end{equation}
Here
\begin{itemize}
  \item $F(Q)$ is $2\pi$-periodic
  \item $P = (I - \hat I) / \sqrt{\varepsilon / \pdv{\omega}{I}(\hat I)}$ measures how far we are from the resonance.
\end{itemize}
We will call the system~\eqref{e:aux-8432} \emph{auxiliary system} describing passage through resonances. This system can be considered as a simple Hamiltonian system
\begin{equation} \label{e:aux-unperturbed-8432}
  Q' = P, \qquad
  P' = F(Q)
\end{equation}
describing a movement of a unit mass particle with coordinate $Q$ and velocity $P$ under the force $F(Q)$ with extra perturbation of order $\sqrt{\varepsilon}$.
\begin{figure}[h]
    \centering
    \includegraphics[width=0.8\textwidth]{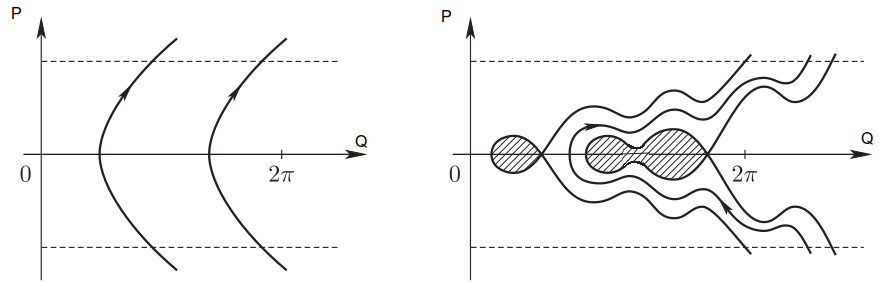} %
    \caption{Possible phase portraits of unperturbed auxiliary system. \\
    Adapted from~\cite{neishtadt2014}. } %
    \label{f:pendulum}
\end{figure}
There are two possibilities. First, $F(Q)$ can have the same sign for all $Q$, the case $F>0$ is depicted in Figure~\ref{f:pendulum}, left. Then capture into resonance is impossible and trajectory leaves resonant zone $|P| \lesssim 1$ after time $O(1)$ for auxiliary system, corresponding to time $O(1/\sqrt{\varepsilon})$ for initial system and evolution of slow variables of order $O(\sqrt{\varepsilon})$. We will say that such resonances are \emph{weak}\footnote{We adapt the terminology from~\cite{neishtadt2014} so that it can be used near separatices. In~\cite{neishtadt2014} resonances were diveded into high-order, weak, and strong, and capture was possible only into strong resonances. We drop high-order resonances and divide resonances simply into strong (capture is possible) and weak (capture is impossible).}. The second possibility is that $F$ changes sign. Points with $F=0$ correspond to equilibria of unperturbed auxiliary system. An example is depicted in Figure~\ref{f:pendulum}, right. Some of these equilibria are saddles, and crossing separatices of these saddles due to perturbation can lead to capture into resonance. For example, in Figure~\ref{f:pendulum} (right) a trajectory going near separatrix of the right saddle can cross this separatrix and enter the dashed domain, staying in the dashed domain. This trajectory of auxiliary system corresponds to a trajectory of the initial system staying $O(\sqrt{\varepsilon})$-close to resonance.

\subsection{Passage through resonant zones near separatrices}
As in the previous subsection, we discuss simplified case without the parameter $z$ to make main ideas more transparent. The perturbation $f$ grows near separatrices in variables $I, \varphi$, but the divergence of $f$ remains $O(1)$, as the coordinate change $p, q \mapsto I, \varphi$ is volume-preserving. This allows us to separate the perturbation (in action-angle variables) into a Hamiltonian part that grows near separatices and a non-Hamiltonian part that is bounded even close to separatrices.
The Hamiltonian part of perturbation has amplitude $O(\sqrt{\varepsilon / h})$ (up to multiplying by some power of $\ln \varepsilon$) and the non-Hamiltonian part has amplitude $O_*(\sqrt{\varepsilon h})$ (here and thereafter the notation $O_*$ denotes that the estimate holds up to multiplying by some power of $\ln h$).

Then rescaling near resonance together with averaging over the remaining fast variable is applied.
As the Hamiltonian part of the perturbation is fairly large, estimates on accuracy of
single step of averaging method are not enough, for example, to get good estimates of measure captured into resonances: the small parameter for averaging method is $\approx \sqrt{\varepsilon / h}$, while we want estimates with accuracy $\approx \sqrt{\varepsilon h}$. However, it is possible to use many steps of averaging method instead of just one to get better accuracy. We use the result~\cite{neishtadt1984separation} on many-step averaging.
This result needs estimates for complex continuation of the pertubation in action-angle variables, we obtain such estimates. This allows to get auxiliary system describing dynamics near resonances.

The resulting auxiliary system is similar to auxiliary system far from separatrices. The unperturbed system is the same, but the perturbation is now divided into $O_*(\sqrt{\varepsilon / h})$ Hamiltonian part and $O_*(\sqrt{\varepsilon h})$ non-Hamiltonian part. Hamiltonian perturbation alone does not lead to capture into resonances, as separatrix loops (cf. Figure~\ref{f:pendulum}) survive.
So amplitude of non-Hamiltonian perturbation determines the measure of captured trajectories, this measure is $O_*(\sqrt{\varepsilon h})$.
The magnitude of scattering on resonances is $O_*(\sqrt{\varepsilon h} |\ln \varepsilon|)$, as width of resonance zones near separatices is $O_*(\sqrt{\varepsilon h})$.

Finally, let us mention one of the technical details. Certain condition on non-degeneracy of unperturbed auxiliary system should hold for all strong resonances, it is needed for estimates on passage through resonance zones.
Far from separatrices there are only finitely many strong resonances, so for generic systems this condition is satisfied for all strong resonances. But strong resonances accumulate on separatrices and their number grows when $\varepsilon \to 0$. To deal with this problem, we show that there are only finitely many limit auxiliary systems near separatices.

\section{Plan of the rest of the paper} \label{s:proof_plans}
In the rest of paper we prove Theorem~\ref{t:main}.
In Section~\ref{s:decompose} a technical theorem on crossing a small neighborhood of separatrices is stated (Theorem~\ref{t:main} follows from it, as dynamics far from separatrices is covered by~\cite{neishtadt2014}) and this technical theorem is splitted into three lemmas on approaching separatrices, crossing separatrices, and moving away from separatrices.
In Section~\ref{s:an-perturbed} estimates on functions describing the perturbed system in action-angle variables and their complex continuation are gathered.
In Section~\ref{s:res-nonres-zones} estimates on Fourier coefficients of the perturbation are obtained, resonant and non-resonant zones are defined, and lemmas on passage through resonant and non-resonant zones are stated. Then these lemmas are used to prove Lemma on approaching separatrices in Section~\ref{s:approach-proof}. Lemma on moving away from separatrices can be proved in the same way.
In Section~\ref{ss:nonres_proof} Lemma on crossing non-resonant zones is proved. In Section~\ref{s:aux} auxiliary system describing movement in resonant zones is obtained and in Section~\ref{s:res-zones-proofs} lemmas on crossing resonant zones are proved.
Lemma on passing separatrices is proved in Section~\ref{s:pass-sep}, thus completing the proof of the technical theorem.
In Section~\ref{s:prob} formula for probabilities of capture into different regions (Proposition~\ref{p:prob}) is proved.
Finally, in Section~\ref{s:no-capture-proof} we sketch a proof of Remark~\ref{r:no-capture}, where our main results is strengthened for a special class of perturbations such that capture into resonance is impossible.

\section{Approaching separatrices and passing through separatrices} \label{s:decompose}
In this section we state a technical theorem on crossing a small neighborhood of separatrices and reduce Theorem~\ref{t:main} to this technical theorem. Then we split the technical theorem into three lemmas on approaching separatrices, crossing separatrices, and moving away from separatrices.
\begin{itemize}
  \item For given $X_{init} \in \mathcal A_3$, denote $\lambda_{init} = \lambda(X_{init})$ and let $X(\lambda)$ be the solution of the perturbed system~\eqref{e:perturbed-pq} with initial data $X(\lambda_{init}) = X_{init}$.
  \item Given $\overline X_0 \in \mathbb R^{n+1}_{h, z}$ and $\lambda_0$, denote by $\overline X(\lambda)$ the solution of the averaged system (one needs to specify in which domain $\mathcal B_i$ or in which union of these domains) with initial data $\overline X(\lambda_0) = \overline X_0$.
\end{itemize}

\begin{theorem} \label{t:sep-pass}
  Given any $z_* \in \mathcal Z_B$, for any small enough $c_z, c_h > 0$ for any $C_0, \Lambda > 0$ there exists $C > 0$ such that for any small enough $\varepsilon$ there exists $\mathcal E \subset \mathcal A_3$ with
  \[
    m(\mathcal E) \le C \sqrt{\varepsilon} |\ln^5 \varepsilon|
  \]
  such that the following holds for any $X_{init} \in \mathcal A_3 \setminus \mathcal E$.

  Suppose at some time
  \[
    \lambda_0 \in [\lambda_{init}, \lambda_{init} + \varepsilon^{-1} \Lambda]
  \]
  the point $X_0 = X(\lambda_0)$ satisfies
  \[
    X_0 \in \mathcal A_3, \qquad
    \norm{z(X_0) - z_*} \le c_z, \qquad
    h(X_0) = c_h.
  \]
  Then there exists $i=1,2$ and $\lambda_1 > \lambda_0$ such that
  \[
    X(\lambda_1) \in \mathcal A_i, \qquad
    h(X(\lambda_1)) = -c_h.
  \]

  Take any $\overline X_0 \in \mathbb R^{n+1}_{h, z}$ with
  \[
    \norm{\overline X_0 - (h(X_0), z(X_0))} < C_0  \sqrt{\varepsilon} |\ln \varepsilon|
  \]
  and consider the solution $\overline X(\lambda)$ of averaged system corresponding to capture from $\mathcal B_3$ to $\mathcal B_i$ with initial data $\overline X(\lambda_0) = \overline X_0$.
  Then for any $\lambda \in [\lambda_0, \lambda_1]$ we have
  \begin{equation} \label{e:thm-est}
    |h(X(\lambda)) - h(\overline X(\lambda))| < C \sqrt{\varepsilon} |\ln \varepsilon|, \qquad
    \norm{z(X(\lambda)) - z(\overline X(\lambda))} < C \sqrt{\varepsilon} |\ln \varepsilon|.
  \end{equation}
\end{theorem}
Let us now prove the main theorem using the technical theorem to cover passage near separatrices and~\cite[Theorem 1 and Corollary 3.1]{neishtadt2014} far from separatrices. Let us now state this result from~\cite{neishtadt2014} using our notation.
\begin{theorem}[{\cite{neishtadt2014}}]
  Pick $\hat{\overline X}_0$ and $\Lambda > 0$. Suppose that solution of the averaged system $\hat{\overline X}$ with initial data $\hat{\overline X}(0) = \hat{\overline X}_0$ stays far from the separatrices for $\lambda \in [0, \varepsilon^{-1} \Lambda]$ and satisfies certain conditions (discussed right after the statement of theorem).

  Then for small enough $r > 0$ for any small enough $\varepsilon > 0$ there exists $\mathcal E \subset A_r(\hat{\overline X}_0)$ with $m(\mathcal E) = O(\sqrt{\varepsilon})$ such that for any
  $X_0 \in A_r(\hat{\overline X}_0) \setminus \mathcal E$ we have
  \[
  \abs{ h(X(\lambda)) - h(\overline X(\lambda)) } = O(\sqrt \varepsilon |\ln \varepsilon|),
  \qquad
  \norm{ z(X(\lambda)) - z(\overline X(\lambda)) } = O(\sqrt \varepsilon |\ln \varepsilon|)
  \]
  for $\lambda \in [\lambda_0, \lambda_0 + \varepsilon^{-1} \Lambda]$,
  where $\lambda_0 = \lambda(X_0)$, $X(\lambda)$ is the solution of perturbed system with initial data $X(\lambda_0) = X_0$ and
  $\overline X(\lambda)$ is the solution of averaged system with initial data $\overline X(\lambda_0) = (h(X_0), z(X_0))$.
\end{theorem}
\noindent For a full statement of the conditions, we refer the reader to~\cite[Section 2]{neishtadt2014}. When we apply this theorem below, these conditions are satisfied, the conditions in Section~\ref{s:main-theorem} are written for this purpose.

We will also need the lemma below, it is proved in Appendix~\ref{a:proof-aux}.
\begin{lemma} \label{l:volume}
  For any $\Lambda > 0$ there exists $C > 1$ such that the flow $g^\lambda$ of~\eqref{e:perturbed-pq} satisfies the following.
  Take open
  $A \subset \mathbb R^{2+n}_{p, q, z} \times [0, 2\pi]_\lambda$
  with $m(A) < \infty$. Then for any $\lambda \in [-\varepsilon^{-1} \Lambda, \varepsilon^{-1} \Lambda]$ we have
  \[
    m(g^\lambda(A)) \le C m(A).
  \]
\end{lemma}
\begin{proof}[Proof of Theorem~\ref{t:main}]
  Take $c_h, c_z > 0$ such that we can apply Theorem~\ref{t:sep-pass} with these constants.
  Recall that
  $\hat{\overline X}_i(\lambda) = (\hat{\overline h}_i(\lambda), \hat{\overline z}_i(\lambda))$
  denotes the solution of averaged system describing capture in $\mathcal B_i$ with $\hat{\overline X}_i(0) = \hat{\overline X}_0$.
  Define $\lambda_+$ by $\hat{\overline h}_1(\lambda_+) = \hat{\overline h}_2(\lambda_+) = 2c_h/3$ and $\lambda_{-, i}$ by $\hat{\overline h}_i(\lambda_{-,i}) = -2c_h/3$, $i=1,2$.
  The number $\Lambda_1$ from the conditions for the main theorem is such that if $\lambda$ is $\varepsilon^{-1} \Lambda_1$-close to $\lambda_*$, we have $|\hat{\overline h}_i(\lambda)| < c_h/2$ and $\norm{\hat{\overline z}_i(\lambda) - z_*} < c_z/2$ for $i=1,2$.
  Thus
  \[
    \lambda_+ < \lambda_* - \varepsilon^{-1} \Lambda_1 < \lambda_* + \varepsilon^{-1} \Lambda_1 < \lambda_{-, i}.
  \]
  This means condition $B$ from~\cite[Section 2]{neishtadt2014} is satisfied for
  \[
    \hat{\overline X}_i(\lambda), \lambda \in [0, \lambda_+]
    \qquad \text{and} \qquad
    \hat{\overline X}_i(\lambda), \lambda \in [\lambda_{-, i}, \varepsilon^{-1} \Lambda].
  \]
  For small enough $r$ solutions $\overline X'(\lambda) = (\overline h'(\lambda), \overline z'(\lambda))$ of averaged system with any initial condition
  $\overline X'(0) \in B_r(\hat{\overline X}_0)$ satisfy
  \begin{equation} \label{e:loc-5925}
      \overline h'(\lambda_+) \in  [c_h/2, 3c_h/4],
      \qquad
      \overline h_i'(\lambda_{-, i}) \in [-3c_h/4, -c_h/2].
  \end{equation}

  By~\cite[Corollary 3.1]{neishtadt2014} we have~\eqref{e:mt-close} for $\lambda \in [0, \lambda_+]$, given that $X_0$ is not in some set $\mathcal E_1$ of measure $\lesssim \sqrt{\varepsilon}$.
  Together with~\eqref{e:loc-5925} for small $\varepsilon$ this implies $h(X(\lambda_0 + \lambda_+)) < 4c_h/5$.
  By continuity we have $h(X(\lambda_0 + \lambda'_+)) = c_h$ for some $\lambda'_+ \in [0, \lambda_+]$.
  Thus we can apply Theorem~\ref{t:sep-pass} (with $\lambda_0$ in this theorem equal to $\lambda_0 + \lambda'_+$). This theorem gives (possibly after reducing $r$) a set $\mathcal E_2$ of measure $\lesssim \sqrt{\varepsilon} |\ln^5 \varepsilon|$ such that if
  $X_0 \not \in \mathcal E_1 \cup \mathcal E_2$,
  there is $i=1,2$ and $\lambda'_{-, i}$ such that $X(\lambda_0 + \lambda'_{-, i}) \in \mathcal A_i$, with $h(X(\lambda_0 + \lambda'_{-, i})) = -c_h$ and~\eqref{e:mt-close} holds for
  $\lambda \in [\lambda_0 + \lambda_+, \lambda_0 + \lambda'_{-, i}]$.
  We have $\overline h_i(\lambda_0 + \lambda'_{-, i}) < -4c_h/5$,
  by~\eqref{e:loc-5925} this implies $\lambda'_{-, i} > \lambda_{-, i}$ and so~\eqref{e:mt-close} holds for $\lambda = \lambda_0 + \lambda_{-, i}$.

  Denote
  \[
    \hat{\overline X}_{-,1} = \hat{\overline X}_1(\lambda_{-, 1}), \qquad
    \hat{\overline X}_{-,2} = \hat{\overline X}_2(\lambda_{-, 2}).
  \]
  By~\cite[Corollary 3.1]{neishtadt2014} there exist $r_-$ and $\mathcal E'_{-, i}, i=1,2$ with $m(\mathcal E'_{-, i}) = O(\sqrt{\varepsilon})$ such that for $i=1,2$ solutions starting in $A_{r_-}(\hat{\overline X}_{-,i}) \setminus \mathcal E'_{-, i}$  are approximated by solutions of the averaged system (with the same initial $h, z$) with error
  $O(\sqrt{\varepsilon} \ln \varepsilon)$.
  Let $\mathcal E_{-, i}$ be the preimage of $\mathcal E'_{-, i}$ under the flow of perturbed system over time $\lambda_{-, i}$, we have $m(\mathcal E_{-, i}) = O(\sqrt{\varepsilon})$ by Lemma~\ref{l:volume}.
  We can now write the exceptional set in the current theorem:
  $\mathcal E = \mathcal E_1 \cup \mathcal E_2 \cup \mathcal E_{-, 1} \cup \mathcal E_{-, 2}$.
  Reducing $r$ if needed, we may assume that solutions of averaged system describing capture in $\mathcal B_i$ starting in $A_r(\hat{\overline X}_0)$ at $\lambda=\lambda_0$ are in $A_{r_-/2}(\hat{\overline X}_{-, i})$ at $\lambda=\lambda_0 + \lambda_{-,i}$ for $i=1,2$. Thus $X(\lambda)$ is $O(\sqrt{\varepsilon}|\ln \varepsilon|)$-close to the solution of averaged system with initial data
  \[
    \overline X(\lambda_0 + \lambda_{-, i}) = (h(X(\lambda_0 + \lambda_{-, i})), z(X(\lambda_0 + \lambda_{-, i})))
  \]
  for
  $\lambda \in [\lambda_0 + \lambda_{-, i}, \lambda_0 + \varepsilon^{-1} \Lambda]$.
  The difference between this solution of averaged system and $\overline X_i(\lambda)$ at the moment $\lambda=\lambda_0 + \lambda_{-, i}$ is $O(\sqrt \varepsilon |\ln \varepsilon|)$, it stays of the same order, as the dynamics in slow time takes time $O(1)$ and the averaged system is smooth far from separatrices. Thus we have~\eqref{e:mt-close} for $\lambda \in [\lambda_0 + \lambda_{-, i}, \lambda_0 + \varepsilon^{-1} \Lambda]$. Now we have proved~\eqref{e:mt-close} for all
  $\lambda \in [\lambda_0, \lambda_0 + \varepsilon^{-1} \Lambda]$,
  as required.
\end{proof}

Let us now split Theorem~\ref{t:sep-pass} into a lemma on approaching the separatrices, lemma on crossing immediate neighborhood of separatrices, and lemma on moving away from the separatrices.
\begin{itemize}
  \item Take $\rho = 5$, suppose we are given $C_\rho > 0$. Denote $h_* = C_\rho \varepsilon |\ln^\rho \varepsilon|$.
\end{itemize}
The immediate neighborhood of separatrices is given by $|h| < h_*(\varepsilon)$.
\begin{lemma}[On approaching separatrices] \label{l:approach-sep}
  Take any $z_* \in \mathcal Z_B$, any small enough $c_h, c_z > 0$, any large enough $C_\rho > 0$.
  Then for any $C_0, \Lambda > 0$ there exists $C > 0$ such that for any small enough $\varepsilon$ there exists $\mathcal E \subset \mathcal A_3$ with
  \[
    m(\mathcal E) < C \sqrt{\varepsilon}
  \]
  such that for any $X_{init} \in \mathcal A_3 \setminus \mathcal E$ the following holds.

  Suppose that at some time
  \[
    \lambda_0 \in [\lambda_{init}, \lambda_{init} + \varepsilon^{-1} \Lambda]
  \]
  the point $X_0 = X(\lambda_0)$ satisfies
  \[
    X_0 \in \mathcal A_3, \qquad
    h(X_0) = c_h, \qquad
    \norm{z(X_0) - z_*} \le c_z.
  \]
  Then at some time $\lambda_1 > \lambda_0$ we have
  \[
    X(\lambda_1) \in \mathcal A_3, \qquad h(X(\lambda_1)) = h_*(\varepsilon).
  \]

  Take any $\overline X_0 \in \mathbb R^{n+1}_{h, z}$ with
  \[
    \norm{\overline X_0 - (h(X_0), z(X_0)} < C_0  \sqrt{\varepsilon} |\ln \varepsilon|.
  \]
  and consider the solution $\overline X(\lambda)$ of averaged system in $\mathcal B_3$ with initial data $\overline X(\lambda_0) = \overline X_0$.
  Then for any $\lambda \in [\lambda_0, \lambda_1]$ we have
  \begin{equation}
    |h(X(\lambda)) - h(\overline X(\lambda))| < C \sqrt{\varepsilon} |\ln \varepsilon|, \qquad
    \norm{z(X(\lambda)) - z(\overline X(\lambda))} < C \sqrt{\varepsilon} |\ln \varepsilon|.
  \end{equation}
  Moreover,
  \begin{equation}
    |h(X(\lambda_1)) - h(\overline X(\lambda_1))| < C \sqrt{\varepsilon}.
  \end{equation}
\end{lemma}
\noindent This lemma will be proved in Section~\ref{s:approach-proof}.

\begin{lemma}[On moving away from separatrices]  \label{l:away-sep}
  Take any $z_* \in \mathcal Z_B$, any small enough $c_h, c_z > 0$, any large enough $C_\rho > 0$.
  Then for any $C_0, \Lambda > 0$ there exists $C > 0$ such that for any small enough $\varepsilon$ there exists $\mathcal E \subset \mathcal A_3$ with
  \[
    m(\mathcal E) < C \sqrt{\varepsilon}
  \]
  such that for any $X_{init} \in \mathcal A_3 \setminus \mathcal E$ the following holds.

  Suppose that for $i=1$ or $i=2$ at some time
  \[
    \lambda_0 \in [\lambda_{init}, \lambda_{init} + \varepsilon^{-1} \Lambda]
  \]
  the point $X_0 = X(\lambda_0)$ satisfies
  \[
    X_0 \in \mathcal A_i, \qquad
    h(X_0) = -h_*(\varepsilon), \qquad
    \norm{z(X_0) - z_*} \le c_z.
  \]

  Take any $\overline X_0 = (\overline h_0, \overline z_0) \in \mathbb R^{n+1}_{h, z}$ with
  \[
    |\overline h_0 - h(X_0)| < C_0  \sqrt{\varepsilon}, \qquad
    \norm{\overline z_0 - z(X_0)} < C_0  \sqrt{\varepsilon} |\ln \varepsilon|
  \]
  and consider the solution $\overline X(\lambda)$ of averaged system in $\mathcal B_i$ with initial data $\overline X(\lambda_0) = \overline X_0$.
  Then for any $\lambda \in [\lambda_0, \lambda_{init} + \varepsilon^{-1} \Lambda]$ we have
  \begin{equation}
    |h(X(\lambda)) - h(\overline X(\lambda))| < C \sqrt{\varepsilon} |\ln \varepsilon|, \qquad
    \norm{z(X(\lambda)) - z(\overline X(\lambda))} < C \sqrt{\varepsilon} |\ln \varepsilon|.
  \end{equation}
\end{lemma}
\noindent This lemma is proved similarly to the previous one, thus we omit the proof. In the proofs of these two lemmas we estimate the difference between solutions of perturbed and averaged system in the chart $w = (I, z)$, and this distance is $O(\sqrt{\varepsilon}|\ln \varepsilon|)$.
This and the estimate $\norm{\pdv{h}{w}} = O(\ln^{-1} h)$ in Lemma~\ref{l:d-f-0-I} below explains why near separatrices (when $\ln h \sim \ln \varepsilon$) the distance in $h$ is $O(\sqrt{\varepsilon})$, while far from separatrices this distance is $O(\sqrt{\varepsilon}|\ln \varepsilon|)$.

\begin{lemma}[On passing separatrices] \label{l:near-sep-local}
  Take any $z_* \in \mathcal Z$, any small enough $c_z > 0$, any large enough $C_\rho > 0$. Then for any $\Lambda > 0$ and $\gamma \in \mathbb R$ there exists $C > 0$ such that for any small enough $\varepsilon > 0$
  there exists a set $\mathcal E \subset \mathcal A_3$
  with
  \[
    m(\mathcal E) \le C \sqrt{\varepsilon} |\ln^{\rho - \gamma + 1} \varepsilon|
  \]
  such that for any
  $X_{init} \in \mathcal A_3 \setminus \mathcal E$
  the following holds.

  Suppose that at some time
  \[
    \lambda_0 \in [\lambda_{init}, \lambda_{init} + \varepsilon^{-1} \Lambda]
  \]
  the point $X_0 = X(\lambda_0)$ satisfies
  \begin{equation}
    X_0 \in \mathcal A_3, \qquad
    \norm{z(X_0) - z_*} \le c_z, \qquad
    h(X_0) = h_*(\varepsilon).
  \end{equation}
  Then at some time $\lambda_1 > \lambda_0$ with
  \[
    \varepsilon(\lambda_1 - \lambda_0) \le  \sqrt{\varepsilon} |\ln^\gamma \varepsilon|.
  \]
  we have
  \[
    X(\lambda_1) \in \mathcal A_1 \cup \mathcal A_2, \qquad
    h(X(\lambda_1)) = -h_1.
  \]
\end{lemma}
\noindent This lemma will be proved in Section~\ref{s:pass-sep}.

\begin{proof}[Proof of Theorem~\ref{t:sep-pass}]
  Suppose that $c_z^{(0)}$ is small enough for all three lemmas, take in the theorem $c_z = c_z^{(0)} / 3$. Take $c_h$ in the theorem such that
  \begin{itemize}
    \item it is small enough for all three lemmas;
    \item while any solution of averaged system describing passage from $\mathcal B_3$ to $\mathcal B_1$ or $\mathcal B_2$ passes from $h = c_h$ to $h = -2 c_h$ the total variation of $z$ is at most $c_z/3$.
  This can be done, as $\dv{z}{h} = O(\ln h)$ for solutions of averaged system.
  \end{itemize}
  Take $C_\rho$ large enough for all three lemmas.
  Take $\Lambda^{(0)}$ such that solutions of averaged systems describing capture in $\mathcal B_1$ and $\mathcal B_2$ starting with $h=c_h$ reach $h=-2c_h$ after time less than $\varepsilon^{-1} \Lambda^{(0)}$ passes. We will apply the three lemmas with $\Lambda$ greater than given in the theorem by $\Lambda^{(0)}$.

  Given $C_0$ and $\Lambda$, Lemma~\ref{l:approach-sep} gives us $C, \mathcal E$ that we denote by $C^{(1)}, \mathcal E^{(1)}$.
  For $X_{init} \in \mathcal A_3 \setminus \mathcal E^{(1)}$ the solution $X(\lambda)$ reaches $h = h_*(\varepsilon)$ at some moment that we denote $\lambda^{(1)}$. For $\lambda \in [\lambda_0, \lambda^{(1)}]$ we have~\eqref{e:thm-est} with $C = C^{(1)}$.
  Note that $\norm{z(X(\lambda^{(1)})) - z_*} < c_z^{(0)}$ for small enough $\varepsilon$ due to~\eqref{e:thm-est} and $\norm{z(\overline X(\lambda^{(1)})) - z_*} < c_z/3$  (this holds by our choice of $c_h$).

   Apply Lemma~\ref{l:near-sep-local} with $\gamma=1$, $c_z = c_z^{(0)}$, $C_0 = C^{(1)}$, it gives us $C, \mathcal E$ that we denote by $C^{(2)}, \mathcal E^{(2)}$.
   For $X_{init} \in \mathcal A_3 \setminus (\mathcal E^{(1)} \cup \mathcal E^{(2)})$ the solution $X(\lambda)$ reaches $h = -h_*(\varepsilon)$ at some moment that we denote $\lambda^{(2)}$.
   We have $\varepsilon(\lambda^{(2)} - \lambda^{(1)}) \le \sqrt{\varepsilon} |\ln \varepsilon|$.
   The change of $h, z, \overline h, \overline z$ during this time is bounded by
   $C^{(2)} \sqrt{\varepsilon} |\ln \varepsilon|$ for some $C^{(2)} > 0$.
   Then for $\lambda \in [\lambda^{(1)}, \lambda^{(2)}]$ we have~\eqref{e:thm-est} with $C = C^{(1)} + 2C^{(2)}$. As above, we get $\norm{z(X(\lambda^{(2)})) - z_*} < c_z^{(0)}$.

   Finally, apply Lemma~\ref{l:away-sep} with $c_z = c_z^{(0)}$ and $C_0 = C^{(1)} + 2 C^{(1)}$, it gives us $C, \mathcal E$ that we denote by $C^{(3)}, \mathcal E^{(3)}$. Set $\mathcal E$ in the theorem equal to
   $\mathcal E^{(1)} \cup \mathcal E^{(2)} \cup \mathcal E^{(3)}$.
   For $X_{init} \in \mathcal A_3 \setminus \mathcal E$
   the solution $X(\lambda)$ reaches $h = -c_h$ at some moment that we denote $\lambda^{(3)}$.
   For $\lambda \in [\lambda^{(2)}, \lambda^{(3)}]$ we have~\eqref{e:thm-est} with $C = C^{(3)}$.
   Take in the theorem $C = C^{(3)} + 2C^{(1)} + 2 C^{(2)}$.
\end{proof}

\section{Analysis of the perturbed system} \label{s:an-perturbed}
We now focus on the proof of the lemma on approaching separatrices. Let us consider the perturbed system in $\mathcal A_3$.
\subsection{Action-angle and energy-angle variables}
In the domain $\mathcal B_3$ let us consider the action-angle variables $I(p, q, z), \varphi(p, q, z)$ of the unperturbed system. We choose the angle variable in such a way that $\varphi = 0$ is given by an analytic transversal to one of the separatrices that is far away from the saddle.
In the energy-angle variables $h, \varphi$ the perturbed system rewrites as
\begin{align}
\begin{split} \label{e:init}
  \dot h &= \varepsilon f_h(h, z, \varphi, \lambda, \varepsilon), \\
  \dot z &= \varepsilon f_z(h, z, \varphi, \lambda, \varepsilon), \\
  \dot \varphi &= \omega(h, z) + \varepsilon f_\varphi(h, z, \varphi, \lambda, \varepsilon), \\
  \dot \lambda &= 1.
\end{split}
\end{align}
We will denote by $\pdv{}{z}$ the partial derivative for fixed $h, \varphi$ and by $\pdv{}{z}_{p, q}$ the partial derivative for fixed $p, q$.
We will often use the action $I$ instead of $h$, then the perturbed system rewrites as
\begin{align}
\begin{split}
  \dot I &= \varepsilon f_I(I, z, \varphi, \lambda, \varepsilon), \\
  \dot z &= \varepsilon f_z(I, z, \varphi, \lambda, \varepsilon), \\
  \dot \varphi &= \omega(I, z) + \varepsilon f_\varphi(I, z, \varphi, \lambda, \varepsilon), \\
  \dot \lambda &= 1,
\end{split}
\end{align}
where $f_I = \pdv{I}{h} f_h + \pdv{I}{z} f_z$.
We will denote $w = (I, z)$ and $f = (f_I, f_z)$. Then the system above rewrites as
\begin{equation} \label{e:perturbed-aa}
  \dot{w} = \varepsilon f(w, \varphi, \lambda, \varepsilon), \qquad
  \dot{\varphi} = \omega(w) + \varepsilon f_\varphi(w, \varphi, \lambda, \varepsilon), \qquad
  \dot{\lambda} = 1.
\end{equation}
Denote
\begin{align}
\begin{split}
  f_0 = \langle  f(w, \varphi, \lambda, 0) \rangle_{\varphi, \lambda}, \qquad
  f_{z, 0} = \langle  f_z(w, \varphi, \lambda, 0) \rangle_{\varphi, \lambda}, \\
  f_{h, 0} = \langle  f_h(w, \varphi, \lambda, 0) \rangle_{\varphi, \lambda}, \qquad
  f_{I, 0} = \langle  f_I(w, \varphi, \lambda, 0) \rangle_{\varphi, \lambda}.
\end{split}
\end{align}
Here $\langle \psi(\varphi, \lambda) \rangle$ denotes the average $\frac 1 {4\pi^2} \int_0^{2\pi} \int_0^{2\pi} \psi(\varphi, \lambda) d\varphi d\lambda$.
Note that this gives the same $f_{h, 0}$ and $f_{z, 0}$ as the formulas in Section~\ref{s:avg}.
The averaged system can be rewritten using $I$ instead of $h$ as follows
\begin{equation} \label{e:avg-w}
  \dot{\overline w} = \varepsilon f_0(\overline w).
\end{equation}
As $\pdv{\omega}{h} \sim h^{-1} \ln^{-2} h > 0$ (see~\eqref{e:est-real} below),~\eqref{e:K-h} implies that close enough to separatrices for some $K_\omega > 0$
\begin{equation} \label{e:K-omega}
  \pdv{\omega}{w} f_0 = \pdv{\omega}{h}\/ f_{h, 0} + \pdv{\omega}{z}\/ f_{z, 0} < -K_{\omega} h^{-1} \abs{\ln^{-3} h} < 0.
\end{equation}
This means that near separatrices $\omega$ decreases along the solutions of averaged system.

The following estimates on the connection between $I$ and $h$ will be proved in Appendix~\ref{a:appendix}. Denote $v = (h, z)$.
\begin{lemma} \label{l:d-f-0-I}
  We have
  \begin{equation}
    \pdv{h}{I} = \omega, \qquad
    \norm{\pdv{w}{v}} = O(\ln h), \qquad
    \norm{\pdv{h}{w}} = O(\ln^{-1} h), \qquad
    \norm{\pdv{f_{I, 0}}{w}}, \norm{\pdv{f_{z, 0}}{w}} = O(h^{-1} \ln^{-3} h).
  \end{equation}
\end{lemma}

\subsection{Complex continuation} \label{ss:complex}
Taking a finite subcover, it is easy to prove that there is $K_1 > 0$ such that the function $H$ and the unperturbed system~\eqref{e:unperturbed} can be continued to the set
\[
  U_{K_1}(\mathcal B) =  \{ v + v',  v=(p, q, z) \in \mathcal B, \; v' \in \mathbb C^{n+2}_{p, q, z}, \; \norm{v'} < K_1 \},
\]
while the functions
$f_p|_{\varepsilon=0}, f_q|_{\varepsilon=0}, f_z|_{\varepsilon=0}$
can be continued to
\[
  U_{K_1}(\mathcal A) =  \{ X + X',  X=(p, q, z, \lambda) \in \mathcal A, \; X' \in \mathbb C^{n+3}_{p, q, z, \lambda}, \; \norm{X'} < K_1 \}.
\]
\noindent Let us now discuss analytic continuation of the angle variable near the separatrices. The proofs of the statements below can be found in Appendix~\ref{a:AC-proofs}.

\noindent For given $\hat \omega$ define $\hat h(z)$ by the equality $\omega(\hat h, z) = \hat \omega$.
\begin{lemma} \label{l:period-cont}
  For any $c_1 > 0$ there exists $c_2 > 0$ such that for any $(h_0, z_0) \in \mathcal B_3$ with $0 < h_0 < c_2$ the following holds.
  Set $\hat \omega = \omega(h_0, z_0)$ and $\hat T = T(h_0, z_0)$.
  Then $\hat h$ is uniquely defined for all $z$ with $\norm{z - z_0} < c_2$ and the period $T(h, z)$ can be continued to
  \[
    \{ (h, z) \in \mathbb C^{n+1}, \; \norm{z - z_0} < c_2, \; |h - \hat h(z)| < c_2 |\hat h(z)| \}
  \]
  with
  $\abs{T(h, z) - \hat T} < c_1 \hat T$.
  Moreover, in the neighborhood above we have $T = A(h, z) \ln h + B(h, z)$, where $A$ and $B$ are bounded analytic functions with $A \ne 0$ and $\ln$ is the branch of the complex logarithm obtained by analytic continuation of the real logarithm.
\end{lemma}

\begin{lemma} \label{l:cont}
  Denote $r(h, z, \varphi) = (p, q)$. Then there is $c \in (0, c_2)$ (here $c_2$ is the constant from Lemma~\ref{l:period-cont}) such that
  for any $z_0$ and any $\hat \omega \in (0, c)$ with $(\hat h(z_0), z_0) \in \mathcal B_3$
  the function $r(h, z, \varphi)$ can be continued to
  \begin{equation} \label{e:D}
    \mathcal D = \big\{
      (h, z, \varphi) \in \mathbb C^{n+2}, \;
      \norm{z - z_0} < c, \; |h - \hat h(z)| < c |\hat h(z)|, \;
      \abs{\Im \varphi} < c \hat \omega
    \big\}
  \end{equation}
  with $(r, z) \in U_{K_1}(\mathcal B)$ and $r(h, z, \varphi) = r(h, z, \varphi + 2\pi)$.
\end{lemma}

Let us now consider a resonance given by $\omega = s_2/s_1$ for coprime positive integers $s_1, s_2$.
Set $\hat \omega = s_2/s_1$.
Given a function $\psi(h, z)$, we denote $\hat \psi(z) = \psi(\hat h(z), z)$.
Take $z_0$ with $(\hat h(z_0), z_0) \in \mathcal B_3$.
By Lemma~\ref{l:cont} the system~\eqref{e:init} can be continued to the complex domain
\begin{equation} \label{e:D-0}
  \mathcal D_0 = \Big\{z, h, \varphi, \lambda \in \mathbb C^{n+3} :
  \abs{z - z_0} < c_{cont}, \;
  |h - \hat h(z)| < c_{cont} |\hat h(z)|, \;
  \abs{\Im \varphi} < c_{cont} \hat \omega, \;
  \abs{\Im \lambda} < c_{cont}
  \Big\}
\end{equation}
for some $c_{cont} > 0$. The function $T(h, z)$ also continues in this domain with
$|T(h, z) - \hat T| < d_{cont} \hat T$, where $\hat T = 2\pi/\hat \omega$
and $d_{cont} > 0$ can be made as small as needed by reducing $c_{cont}$, by Lemma~\ref{l:period-cont}. The constants $c_{cont}$ and $d_{cont}$ are uniform, they do not depend on $s_2, s_1, h_0, z_0$.

\noindent Finally, we need the following technical lemma.
\begin{lemma} \label{l:int-psi}
  Consider an analytic function $\psi(p, q, z, \lambda)$ (such that $\psi$ can be continued to $U_{K_1}(\mathcal A)$) such that $\psi=0$ at $C(z)$ for all $z$ and $\lambda$. Rewrite this function in the chart $h, z, \varphi, \lambda$. Then for any $(h, z, \varphi_0, \lambda) \in \mathcal D_0$ with $\Re \varphi_0 \in [0, 2\pi]$ we have
  \[
    \omega^{-1} \int_{\varphi=0}^{\varphi_0} \psi(h, z, \varphi, \lambda) d\varphi = O(1).
  \]
  Here the integral is taken along any path homotopic to the segment $[0, \varphi_0]$.
\end{lemma}

\subsection{Estimates}

\begin{lemma}
  We have the following estimates and equalities
  valid in $\mathcal D_0$, the constants in $O$-estimates below do not depend on $s_1, s_2, h_0, z_0$.
  \begin{align} \label{e:est}
  \begin{split}
    &T \sim \ln h, \;
    \pdv{T}{h} \sim h^{-1}, \;
    \pdv{T}{z} = O(\ln h), \\
    &\pdv{\omega}{h} \sim h^{-1} \ln^{-2} h, \;
    \pdv[2]{\omega}{h}{z} = O(h^{-1} \ln^{-2} h), \;
    \pdv[2]{\omega}{h} = O(h^{-2} \ln^{-2} h), \\
    &\pdv{\omega}{I} \sim h^{-1} \ln^{-3} h, \;
    \pdv[2]{\omega}{I} = O(h^{-2} \ln^{-4} h), \\
    &\pdv{h}{I} = \omega, \;
    \pdv{I}{z} = -\frac{1}{2\pi} \int_{t=0}^{T} \pdv{H}{z}\Big|_{p, q}(h, z, \omega t) dt = O(1), \\
    &\pdv{\hat h}{z} = O(h \ln h).
  \end{split}
  \end{align}
  Here and below in this paper we use the expressions $y = O(x)$ and $y \sim x$ for negative or complex $x$ and $y$ as a shorthand for $|y| = O(|x|)$ and $|y| \sim |x|$.
\end{lemma}
\begin{proof}
  We have $\pdv{h}{I} = \omega$ from the Hamiltonian equations in the coordinates $I, \varphi$.
  The estimates for $T$ and $\omega$ and their derivatives follow from the last part of Lemma~\ref{l:period-cont} and the formula for $\pdv{h}{I}$. The formula for $\pdv{I}{z}$ is well known. It can be found, e.g., in~\cite[Corollary~3.2]{neishtadt17}, where the estimate $\pdv{I}{z} = O(1)$ is proved in the real case. In the complex case this estimate follows from the formula for $\pdv{I}{z}$ by Lemma~\ref{l:int-psi}.

  Let us prove the estimate for $\pdv{\hat h}{z}$. We have $0 = \pdv{\hat \omega}{z} =
  \hatw{\pdv{\omega}{h}}
   \pdv{\hat h}{z} + \hatw{\pdv{\omega}{z}}$, thus
  $\pdv{\hat h}{z} = - \frac{\hatw{\pdv{\omega}{z}}}{\hatw{\pdv{\omega}{h}}} = O(h \ln h)$.
  Here we use the notation $\hat \psi(z) = \psi(z, \hat h(z))$.

\end{proof}

Let us move on from estimates on the complex continuation of the perturbed system to estimates on the real perturbed system.
We will use the notation $O_*$. A precise definition can be found in~\cite[Table 1]{neishtadt2020phase}.
Roughly speaking, $g = O_*(h^a \ln^b h)$ means that $g h^{-a} \ln^{-b} h = O(1)$ and fastly decreases near the saddle $C$. We will need the following fact (\cite[Lemma 11.1]{neishtadt2020phase})
\begin{equation} \label{e:int-O-star}
  \int_0^{2\pi} O_*(h^a \ln^b h) d\varphi = O(h^a\ln^{b-1} h).
\end{equation}

\begin{lemma}
Let $\psi(p, q, z)$ denote any smooth function such as\footnote{
We can ignore that $f_h$ and $f_z$ also depends on $\lambda$, as we can just use the estimates on $\psi$ for each fixed value of $\lambda$.
} $f_h$ or $f_z$.
Then we have in the \emph{real} domain $\mathcal B_3$
\begin{align} \label{e:est-real}
\begin{split}
  &\norm{\pdv{\psi}{h}} = O_*(h^{-1} \ln^{-1} h), \qquad
  \norm{\pdv{\psi}{z}} = O(1), \\
  &f_\varphi = O_*(h^{-1} \ln^{-2} h), \qquad
  f_I = O(\ln h), \qquad
  f_{h, 0} \sim \ln^{-1} h, \qquad
  f_{I, 0} = O(1), \qquad f_{z, 0} = O(1), \qquad  \\
  &\norm{\pdv{f}{w}} = O_*(h^{-1} \ln^{-1} h), \qquad
  \pdv{\omega}{h} \sim h^{-1} \ln^{-2} h, \qquad
  \norm{\pdv{\omega}{z}} = O(\ln^{-1} h).
\end{split}
\end{align}
These estimates hold with $O$-estimates bounded from above by uniform constants that depend only on the system and the domain $\mathcal B_3$ and $\sim$-estimates bounded from above and from below by uniform constants.
\end{lemma}
\begin{proof}
  For the proofs of these estimates (except the estimate for $\pdv{f}{w}$) see~\cite[Table 1]{neishtadt2020phase} and the references within.
  We have $\pdv{f_z}{I} = \omega \pdv{f_z}{h} = O_*(h^{-1} \ln^{-2} h)$ and
  \[
    \norm{\pdv{f_I}{I}} = \norm{\pdv{}{I}\/(\omega^{-1} f_h)} =\norm{ \omega^{-1} \pdv{f_h}{I} + \pdv{}{I}\/(\omega^{-1}) f_h} = O_*(h^{-1} \ln^{-1} h).
  \]
  In the same way we get $\norm{\pdv{f_z}{z}} = O(1)$ and $\norm{\pdv{f_I}{z}} = O(\omega^{-1})$.
  This implies the estimate for $\pdv{f}{w}$.
\end{proof}

\section{Resonant and non-resonant zones} \label{s:res-nonres-zones}
\subsection{Fourier series}
Denote by $f_m$, $m = (m_1, m_2) \in \mathbb Z^2$ the Fourier coefficients of the vector-valued function $f|_{\varepsilon=0}$:
\begin{equation}
  f(w, \varphi, \lambda, 0) = f_0(w) + \sum_{|m| \ne 0} f_m(w) e^{i(m_1 \varphi + m_2 \lambda)}.
\end{equation}
The following lemma will be proved in Appendix~\ref{a:fourier_proof}, using the results on complex continuation stated in Subsection~\ref{ss:complex}.
\begin{lemma} \label{l:fourier}
  There is $C_F$ such that for any $m \in \mathbb Z^2$ we have in $\mathcal B_3$
  \begin{equation} \label{e:fourier}
    \norm{f_m} \lesssim \exp\bigg( - C_F |m_2| - C_F \frac{|m_1|}{T} \bigg).
  \end{equation}
  Moreover,
  \begin{align} \label{e:d-fourier}
  \begin{split}
    \norm{\pdv{f_m}{h}} \lesssim |h^{-1} \ln^{-1} h| \; \exp(- C_F |m_2|), \qquad
    \norm{\pdv{f_m}{z}} \lesssim \exp(- C_F |m_2|),
  \end{split}
\end{align}
\end{lemma}

For given $N$ that will be chosen in the lemma below let us denote
\begin{equation}
  R_N f(w, \varphi, \lambda) = \sum_{|m| > N} f_m(w) e^{i(m_1 \varphi + m_2 \lambda)}.
\end{equation}
\begin{lemma} \label{l:fourier-remainder}
  There is $N \sim \ln^2 \varepsilon$ such that for $h > \varepsilon$ we have $\norm{R_N f} < \varepsilon$.
\end{lemma}
\begin{proof}
  As $h > \varepsilon$, we have $T \lesssim |\ln \varepsilon|$.
  Hence, we can take $N \sim \ln^2 \varepsilon$ such that $C_F N/T \ge 2 |\ln \varepsilon|$, then $\exp (- C_F |m_2| - C_F \frac{|m_1|}{T} ) \le \varepsilon^2$ for any $m$ with $\abs{m} > N$. From this it is easy to obtain $\norm{R_N f} < \varepsilon$.
\end{proof}

\subsection{Non-resonant zones} \label{s:res-zones}
A resonance is given by $\omega(h, z) = \xi > 0$, where $\xi = s_2 / s_1$ with $s = (s_1, s_2) \in \mathbb Z_{> 0}^2$; $\abs{s} = \abs{s_1} + \abs{s_2} \le N$ and the numbers $s_1$ and $s_2$ are coprime. Assume that the resonances are enumerated in such a way that $\xi_1 > \xi_2 > \xi_3 > \dots$. Let $s_r$ denote the vector $(s_1, s_2)$ corresponding to $\xi_r$.
Take small $k > 0$ and denote by $\mathcal B_{3, *} \subset \mathcal B_3$ the set of points $(h, z)$ that satisfy $\norm{w(h, z) - w(h{=}0, z_*)} < k$.
The value of $k$ is picked so that $\mathcal B_{3, *}$ is far from $\partial \mathcal B$.
The constants $c_h, c_z$ in Theorem~\ref{t:sep-pass} are small enough; we will assume that $\mathcal B_{3, *}$ contains the set given by $h \in (0, 2c_h)$, $\norm{z - z_*} < 2c_z$.
Let us fix $\gamma = 5$ and consider the sets
\begin{equation} \label{e:Pi}
  \Pi = \{ (h, z) \in \mathcal B_{3, *} : h \ge 2 \varepsilon |\ln^{\gamma} \varepsilon| \},
  \qquad
  \partial \Pi = \{ (h, z) \in \mathcal B_{3, *} : h = 2 \varepsilon |\ln^{\gamma} \varepsilon| \}.
\end{equation}
Inside $\Pi$ (and also inside the zone $\tilde \Pi$ defined below in Section~\ref{ss:nonres_proof}) we have (for any $\alpha > 0$)
\begin{equation} \label{e:in-Pi}
  \varepsilon / h \le \abs{\ln^{-\gamma} \varepsilon} < 1,
  \qquad
  |\ln h| < |\ln \varepsilon|,
  \qquad
  \varepsilon h^{-1} \abs{\ln^{-\alpha} h} \le \abs{\ln^{-(\gamma + \alpha)} \varepsilon}.
\end{equation}
For a resonance $s$ denote $b_s = e^{-C_F s_2} < 1$, where $C_F$ is defined in~\eqref{e:fourier}.
Recall that when $s$ is fixed, $\hat h(z)$ denotes the resonant value of $h$ given by $\omega(\hat h(z), z) = s_2/s_1$.
Set
\begin{equation} \label{e:delta}
  \delta_s(z) =
  \sqrt{\varepsilon b_s \hat h^{-1} \ln^{-4} \hat h}
  + \varepsilon \hat h^{-1} |\ln^{-3} \hat h| \ln^2 \varepsilon.
\end{equation}
Take any large enough $C_Z$ (it should be greater than some constant $K_Z$ which will be determined in the proof of Lemma~\ref{l:nonres}).
For each resonance given by $s = (s_1, s_2)$, define its \emph{inner resonant zone}
\begin{equation} \label{e:Z-r}
  Z(s) = \Big\{ (h, z) \in \Pi: |\omega(h, z) - s_2/s_1| \le C_Z \delta_s(z) \Big\}.
\end{equation}
Set $Z_r = Z(s_r)$.
Such resonant zone has width $\sim \sqrt{\varepsilon h}$ in $h$ if $s_2 \sim 1$ and $h \gtrsim \varepsilon \ln^2 \varepsilon$.

The following lemma shows that the value of $h$ stays roughly the same between two neighboring resonances for fixed $z$.
\begin{lemma} \label{l:est-xi-r-r+1}
  Let $\xi_{r} > \xi_{r+1}$ be neighboring resonances.
  Fix $z$ and let $\hat h_r > \hat h_{r+1}$ be the corresponding values of $h$.
  Then if $\hat h_{r+1} \ge \varepsilon |\ln^\gamma \varepsilon|$, we have
  \begin{equation} \label{e:est-xi-r-r+1}
    \hat h_r < \Big( 1 + O(|\ln^{-1} \varepsilon|) \Big) \hat h_{r+1}.
  \end{equation}
\end{lemma}
\begin{proof}
  As $\omega \sim \ln^{-1} h$, $\hat h_{r+1} \ge \varepsilon |\ln^\gamma \varepsilon|$ implies $\xi_{r+1} \gtrsim |\ln^{-1} \varepsilon|$.
  As $\xi_r^{-1}$ and $\xi_{r+1}^{-1}$ are two neighboring rational numbers that can be written as $s_1/s_2$ with $|s_1| + |s_2| \le N$, $N \sim \ln^2 \varepsilon$, taking fixed $s_2 \approx N \xi_r/3$ and changing $s_1 \approx N/3$ with unit step gives the estimate
  \[
    |\xi_r^{-1} - \xi_{r+1}^{-1}| \lesssim s_2^{-1} \sim (\xi_r N)^{-1} \lesssim |\ln^{-1} \varepsilon|.
  \]
 Integrating the equality $\pdv{}{h}\/(\omega^{-1}) \sim h^{-1}$, we get
 \[
  |\ln \hat h_r - \ln \hat h_{r+1}| \sim |\xi_r^{-1} - \xi_{r+1}^{-1}| \lesssim |\ln^{-1} \varepsilon|.
 \]
 Taking exponent gives~\eqref{e:est-xi-r-r+1}.
\end{proof}

\begin{lemma} \label{l:no-intersection}
  For any $C_Z > 0$ for all sufficiently small $\varepsilon$ each point $w \in \Pi$ lies inside at most one zone $Z_r$.
\end{lemma}
\begin{proof}
  Take two resonances $s$ and $s'$.
  We have
  \[
    |s_2/s_1 - s'_2/s'_1| \ge s_1^{-1} (s'_1)^{-1} = s_2^{-1} (s_2/s_1) (s'_1)^{-1} \gtrsim s_2^{-1} \ln^{-3} \varepsilon
  \]
  and (as $b_s s_2^2 = O(1)$)
  \[
    (\delta_s / |s_2/s_1 - s'_2/s'_1|)^2 \lesssim \varepsilon h^{-1} \ln^{-4} h \ln^6 \varepsilon \lesssim \ln^{2 -\gamma} \varepsilon
  \]
  with the last inequality following from~\eqref{e:in-Pi}.
  Similarly, we have $(\delta_{s'} / |s_2/s_1 - s'_2/s'_1|)^2 \lesssim \ln^{2 -\gamma} \varepsilon$.
  Thus for $\gamma_2 > 2$ and small enough $\varepsilon$ these resonant zones are disjoint.
\end{proof}
\begin{remark}
  Take an integer $k \sim \ln \varepsilon$,
  $\frac{s_2'}{s_1'} = \frac{k-1}{k^2}$,
  $\frac{s_2}{s_1} = \frac{1}{k+1}$.
  We have
  \[
    \abs{\frac{s_2'}{s_1'} - \frac{s_2}{s_1}} = \frac{1}{k^2(k+1)} \sim \ln^{-3} \varepsilon.
  \]
  As $b_s \sim 1$ (and even $a_s \sim 1$ for $a_s$ defined later in Subsection~\ref{ss:nonres_proof}) and $\ln h \sim \omega^{-1} \sim \ln \varepsilon$, we have
  \[
    \delta_{s} \gtrsim \sqrt{\varepsilon h^{-1} \ln^{-4} \varepsilon}.
  \]
  To avoid $s_2'/s_1'$ being inside the resonant zone of $s_2/s_1$, we need $h \gtrsim \varepsilon \ln^2 \varepsilon$.
\end{remark}

Denote by $Z_{r, r+1} \subset \Pi$
the zone between two neighboring resonance zones $Z_r$ and $Z_{r+1}$, $\xi_r > \xi_{r+1}$. For a more formal definition we need to denote by $S_z \subset \Pi$ the set of all $w \in \Pi$ with given $z$. The intersection $Z_{r, r+1} \cap S_z$ is defined as the segment between $Z_r \cap S_z$ and $Z_{r+1} \cap S_z$ if both are non-empty. If one of these sets is empty, we take the segment between the other set and an endpoint of the segment $S_z$. If both these sets are empty, $Z_{r, r+1} \cap S_z$ is also empty.

\subsection{Lemma on crossing non-resonant zones}
It is convenient to use action-angle variables to describe crossing non-resonant zones.
We will denote $W = (w, \varphi, \lambda)$.
Take some initial data
$W_{init} = (w_{init}, \varphi_{init}, \lambda_{init}) \in \mathcal A_3$.
Let
$W(\lambda) = (w(\lambda), \varphi(\lambda), \lambda)$
be the solution of the perturbed system~\eqref{e:perturbed-aa} with this initial data.
Let us also denote by $\overline w(\lambda)$ some solution of the averaged system~\eqref{e:avg-w} with $\overline w(\lambda_{init})$ close to $w(\lambda_{init})$.
We will use the notation $h(\lambda) = h(w(\lambda))$, $\overline h(\lambda) = h(\overline w(\lambda))$.

Along the solution of the averaged system $\omega$ decreases due to~\eqref{e:K-omega}. Until this solution reaches $\partial \Pi$, it passes the zones in the following order: $Z_1, Z_{1, 2}, Z_2, Z_{2, 3}, \dots$. The evolution of slow variables given by the perturbed system $w(\lambda)$ passes the zones more or less in the same order, but as it oscillates, it can leave and reenter the zones. The lemma below covers the times from the first entry to $Z_{r, r+1}$ until the first entry to $Z_{r+1}$ (or reaching $\partial \Pi$).

\begin{lemma} \label{l:nonres}
  Fix $z_* \in \mathcal Z$, then there exist $\gamma_1, K_Z > 0$ such that the following holds for any $C_Z > K_Z$ for some $C_6, C_7, C_t > 0$ and any small enough $c_z, \omega_0, \varepsilon > 0$.
  Suppose at some time $\lambda_0 > \lambda_{init}$ we have
  \[
    w(\lambda_0) \in Z_{r, r+1}, \qquad
    \omega(w(\lambda_0)) < \omega_0.
  \]

  \noindent a) Then there exists $\lambda_1 > \lambda_0$ such that $w(\lambda_1)$ lies
  in $\partial \Pi$ or on the border between $Z_{r, r+1}$ and $Z_{r+1}$.

  \noindent b)
  \begin{equation}
    \varepsilon (\lambda_1 - \lambda_0) \le C_t (\xi_r - \xi_{r+1}) h(\lambda_1) \abs{\ln^3 h(\lambda_1)}.
  \end{equation}

  \noindent c) $h(\lambda) \le \frac{5}{3} h(\lambda_0)$ for $\lambda \in [\lambda_0, \lambda_1]$.

  \noindent d)
  Assume that for some $\lambda_{01} \in [\lambda_0, \lambda_1]$ we also have
  \begin{equation} \label{e:h-overline-h}
    \overline h(\lambda_{01}) > 0.5 h(\lambda_{01}),
    \qquad
    \overline h(\lambda_0) < 2 h(\lambda_0).
  \end{equation}
  Denote $b_{s_r} = e^{-C_F s_2}$, where $\xi_r = s_2/s_1$. Denote $b_{s_{r+1}}$ by the same formula with $\xi_r$ replaced by $\xi_{r+1}$.
  Then for all $\lambda \in [\lambda_0, \lambda_{01}]$ we have (in the error term below we denote $h = h(\lambda_1)$)
  \begin{equation} \label{e:nonres-est}
    \norm{w(\lambda) - \overline w(\lambda)}
    <
    e^{
      C_6 (\xi_r - \xi_{r+1})
    }
    \norm{w(\lambda_0) - \overline w(\lambda_0)}
    + C_7 (\sqrt{b_{s_r}} + \sqrt{b_{s_{r+1}}}) \sqrt{\varepsilon h} \abs{\ln h}
    + C_7 \varepsilon |\ln^{\gamma_1} \varepsilon|.
  \end{equation}
\end{lemma}

\subsection{Outer and middle resonant zones}

Suppose we are given $C_Z, C'_Z$ with $C'_Z > C_Z > 0$ and numbers $D'(s)$ for each resonance $s = (s_1, s_2)$.
Define the \emph{outer resonant zone}
\begin{equation} \label{e:outer-res}
  Z'(s) = \Big\{ (h, z) \in \Pi: |\omega(h, z) - s_2/s_1| \le C'_Z \delta_s(z) \Big\}.
\end{equation}
Fix some resonant $\hat \omega = s_2/s_1$.
Given a function $\psi(h, z)$, denote $\hat \psi(z) = \psi(\hat h(z), z)$.
Denote
\begin{equation} \label{e:alpha}
  \alpha(z) = \alpha_s(z) = \sqrt{\varepsilon / \hatw{\pdv{\omega}{I}}}
  \sim \sqrt{\varepsilon \hat h \ln^3 \hat h}.
\end{equation}
Let $D(h, z) = D_s(h, z)$ be defined by
\begin{equation}
  I(h, z) = \hat I(z) + D \alpha(z) \sqrt{\hat \omega}.
\end{equation}
The \emph{middle resonant zone} is defined in terms of $D$ as
\begin{equation} \label{e:middle-res}
  Z^m(s) = \Big\{ (h, z) \in \Pi: |D(h, z)| \le D'(s) \Big\}.
\end{equation}
Set $Z^m_r = Z^m(s_r)$ and $Z'_r = Z'(s_r)$.

\begin{lemma} \label{l:3-res-zones}
  Given $K, S_2 > 0$, there exists $K_1 > K$ such that for any $D_{p, 0} > K_1$ we can pick $C'_Z > C_Z > K$ and functions $D'_p(s, z)$ such that for small enough $\varepsilon$ for any $s$ we have
  \begin{itemize}
    \item $Z(s) \subset Z^m(s) \subset Z'(s)$,
    \item $D'_p(s, z) \le D_{p, 0}$ and $D'_p(s, z) = D_{p, 0}$ if $s_2 > S_2$,
    \item $D'_p(s, z) = c_1 D_{p, 0} \Big(
      \sqrt{b_s} + \sqrt{\varepsilon \hat h(z)^{-1}} \ln^{-1} \hat h(z) \ln^2 \varepsilon
      \Big)$ for some $c_1 > 0$.
  \end{itemize}
\end{lemma}
\begin{proof}
  Denote
  \[
    Z(s, C) = \Big\{ (h, z) \in \Pi: |\omega(h, z) - s_2/s_1| \le C \delta_s(z) \Big\}.
  \]
  This gives inner resonant zones for $C = C_Z$ and outer resonant zones for $C = C'_Z$.
  Let $D_s(z, C)$ be the width of this zone in $D$.
  It is easy to check that $\delta_s(z) = O(\ln^{-1} \varepsilon)$ if $(\hat h_s(z), z) \in \Pi$, so for fixed $z$ the values of $h$ in $Z(s, C)$ differ by $o(h)$. This implies $\pdv{\omega}{I} \sim \hat{\pdv{\omega}{I}}$ in this zone (for any $C$ this holds for small enough $\varepsilon$).

  As $\gamma \ge 2$, we have $\varepsilon \hat h^{-1} \ln^{-3} h \ln^2 \varepsilon \lesssim \sqrt{\varepsilon \hat h^{-1} \ln^{-4} h }$.
  Denote $\tilde \delta_s = \sqrt{\varepsilon h^{-1} \ln^{-4} h}$.
  We have $\delta_s \le C_1 \tilde \delta_s$ for some $C_1 > 0$.
  We have
  \begin{equation}
    D_s(z, C) \sim C \delta_s / (\hat{\pdv{\omega}{I}} \alpha \sqrt{\hat \omega})
    \sim C \delta_s / \sqrt{\varepsilon \hat{\pdv{\omega}{I}} \hat \omega}
    \sim C \frac{\delta_s}{\tilde \delta_s}.
  \end{equation}
  Take $d'(s) =  C_1^{-1} \delta_s / \tilde \delta_s \le 1$ if $s_2 > S_2$ and $d'(s) = 1$ otherwise.
  Note that $\delta_s \sim \tilde \delta_s$ if $s_2 \le S_2$ as $b_s \sim 1$ for such $s$.
  Then for some $C_2 > 1$ we have
  \begin{equation}
    D_s(z, C) \in [C C_2^{-1} d'(s) , C C_2 d'(s)] \qquad \text{for all } s.
  \end{equation}
  Note that the value of $C_2$ does not depend on $K$ and $D_{p, 0}$, it only depends on $S_2$.
  Take
  \[
    C_Z = D_{p, 0} C_2^{-1}, \qquad
    C'_Z = D_{p, 0} C_2, \qquad
    D'_p(s) = D_{p, 0}d'(s)
  \]
  and $K_1$ so large that $D_{p, 0} > K_1$ implies $C_Z > K$.
  We have
  \[
    D_s(z, C_Z) \le C_Z C_2 d'(s) = D_{p, 0} d'(s) = D'_p(s), \qquad
    D_s(z, C'_Z) \ge C'_Z C^{-1}_1 d'(s) = D_{p, 0} d'(s) = D'_p(s).
  \]
  This implies $Z_r(s) \subset Z^m_r(s) \subset Z'_r(s)$.
\end{proof}

\subsection{Lemmas on crossing resonant zones} \label{s:est-res-cross}
We will call a resonance high-numerator if $s_2 < S_2$ and low-numerator otherwise. The constant $S_2$ is picked in the proof of Lemma~\ref{l:est-weak-res} below.
Denote by $B_r(O)$ the open ball with center $O$ and radius $r$.
The lemma below covers crossing middle resonant zones of high-numerator resonances.
\begin{lemma} \label{l:est-weak-res}
  Given $z_* \in \mathcal Z$,
  for any $C_{s_1} > 0$ and $D_p > 1$
  for any large enough $S_2$
  there exist $C_\rho, C > 1$
  such that
  for any small enough $c_z, \omega_0, \varepsilon_0 > 0$
  for any $\hat \omega = s_2 / s_1 \in (0, \omega_0)$ with $s_2 > S_2$ and $s_1 < C_{s_1} \ln^2 \varepsilon$
  we have the following for any $D'_p \in (0, D_p]$.
  Take some initial condition $X_0 = (p_0, q_0, z_0, \lambda_0)$ with
  \[
    |I(X_0) - \hat I(z_0)| \le D'_p \alpha(z_0) \hat \omega^{0.5}, \qquad
    z_0 \in B_{c_z}(z_*), \qquad
    h(X_0) > C_\rho \varepsilon |\ln \varepsilon^5|.
  \]
  Denote by $X(\lambda)$ the solution of the perturbed system with this initial data. Then this solution crosses the hypersurface
  \[
    I = \hat I(z) - D'_p \alpha(z) \hat \omega^{0.5}
  \]
  at some time $\lambda_1 > \lambda_0$ with
  \[
    \varepsilon(\lambda_1 - \lambda_0)
    \le C D'_p \alpha(z_0) \hat \omega^{0.5} + 2\pi s_1 \varepsilon.
  \]
\end{lemma}

\noindent This lemma is proved in Section~\ref{ss:weak-res-proof}.

\noindent The next lemma covers crossing middle resonant zones of low-numerator resonances.
\begin{lemma} \label{l:est-res-cross}
  Given $z_* \in \mathcal Z_B$, for any
  $C_{s_1}, D_{p, 0}, \Lambda > 1$
  there exist
  \[
    D_p > D_{p, 0}, \; C_\rho, C > 1
  \]
  such that for any small enough $c_z, \omega_0, \varepsilon > 0$
  for any
  $\hat \omega = s_2 / s_1 \in (0, \omega_0)$
  with
  $s_1 < C_{s_1} \ln^2 \varepsilon$
  we have the following.
  Set
  \begin{equation}
    \alpha_* = \max_{B_{c_z}(z_*)} {\alpha(z)}.
  \end{equation}
  Then there exists a set $\mathcal E_s \subset \mathcal A_3$ with
  \[
    m(\mathcal E_s) \le C s_1 \hat \omega^{-1} \alpha_*
  \]
  such that the following holds.

  Take some initial condition $X_{init} = (p_{init}, q_{init}, z_{init}, \lambda_{init}) \in \mathcal A_3 \setminus \mathcal E_s$. Let $X(\lambda)$ be the solution of the perturbed system with this initial data. Suppose that this solution
  crosses the hypersurface
  \[
    I = \hat I(z) + D_p \alpha(z) \hat \omega^{0.5}
  \]
  at some time $\lambda_0 > \lambda_{init}$ with $\varepsilon (\lambda_0 - \lambda_{init}) \le \Lambda$, $z(\lambda_0) \in B_{c_z}(z_*)$ and
  $h(\lambda_0) > C_\rho \varepsilon |\ln \varepsilon^5|$.
  Then this solution crosses the hypersurface
  \[
    I = \hat I(z) - D_p \alpha(z) \hat \omega^{0.5}
  \]
  at some time $\lambda_1 > \lambda_0$ with
  \[
    \varepsilon(\lambda_1 - \lambda_0)
    \le C \alpha(z_0) |\ln \varepsilon| \hat \omega^{0.5}
    \le C \alpha_* |\ln \varepsilon| \hat \omega^{0.5}.
  \]
\end{lemma}
\noindent This lemma is proved in Section~\ref{ss:low-num-res}.

\begin{remark} \label{r:res-effects-decrease-formal}
  Resonances near separatrices have smaller effect on the dynamics.
  Namely, exceptional set in Lemma~\ref{l:est-res-cross} has measure $O(\sqrt{\varepsilon h_*})$ (up to some power of $\ln h_*$), where $h_* = \min_{B_{c_z}(z_*)} \hat h(z)$. Time spent in resonant zone (i.e. between two hypersurfaces with given $I$) in Lemma~\ref{l:est-res-cross} and Lemma~\ref{l:est-weak-res} is $O(\sqrt{h_* / \varepsilon})$ (up to some power of $\ln h_*$).
  The change of slow variables $h$ and $z$ while inside resonant zone is $O(\sqrt{\varepsilon h_*})$ (up to some power of $\ln h_*$).
\end{remark}
\begin{proof}
  The first two estimates (on measure of exceptional set and time spent in resonant zones) follow from $\alpha \sim \sqrt{\varepsilon h \ln^3 h}$ and Lemmas~\ref{l:est-weak-res} and~\ref{l:est-res-cross}. The estimate on the change of slow variables follows from the estimate on time spent in resonant zones, as the rate of change of slow variables is $O(\varepsilon)$.
\end{proof}

\newpage
\section{Proof of the lemma on approaching separatrices} \label{s:approach-proof}

\subsection{Picking constants and excluded set $\mathcal E$}
Considerations for non-resonant zones give us $C_{s_1}$ and some bound from below $K_Z$ on $C_Z$.
Lemma~\ref{l:est-res-cross} gives us the value of $S_2$.
Lemma~\ref{l:3-res-zones} (with $K = K_Z$ and $S_2$ fixed above) gives us $D_{p, 0}$.
Apply Lemma~\ref{l:est-res-cross} (with this $D_{p, 0}$) to get $D_p$.
Plug $K_1 = D_p$ in Lemma~\ref{l:3-res-zones} to get $C'_Z > C_Z > K_Z$ and $D'_p(s, z) \le D_p$ (and $D'_p(s, z) = D_p$ if $s$ is a low-numerator resonance).
Now resonant (inner, outer, middle) zones and non-resonant zones (zones between inner resonant zones) are defined.
Pick $c_z$, $\omega_0$ so that lemmas~\ref{l:nonres},~\ref{l:est-weak-res} and~\ref{l:est-res-cross} hold with this $\omega_0$ and with $c_z$ in these lemmas equal to $2 c_z$.
Pick $\Lambda > 0$ so that any solution of the averaged system under consideration starting with $\lambda = \lambda_0$ crosses $h=-c_h$ after time less than $\varepsilon^{-1} \Lambda$.
Let $\mathcal E_1$ be the union of excluded sets $\mathcal E_s$ for all low-numerator resonances provided by Lemma~\ref{l:est-res-cross} and $\mathcal E_2$ be the union of all middle resonant zones of all low-numerator resonances.
Define the excluded set $\mathcal E$ by $\mathcal E = \mathcal E_1 \cup \mathcal E_2$.

\subsection{Passing resonant and non-resonant zones}
Consider solutions $X(\lambda)$ of~\eqref{e:perturbed-pq} and $\overline X(\lambda)$ of~\eqref{e:avg-h-z} as in the statement of Lemma~\ref{l:approach-sep} (then $X(\lambda_{init}) \in \mathcal A_3 \setminus \mathcal E$).
Denote
\[
  w(\lambda) = (I(h(\lambda), z(\lambda)), z(\lambda)) = w(X(\lambda)), \qquad
  \overline w(\lambda) = (I(\overline h(\lambda), \overline z(\lambda)), \overline z(\lambda)) = w(\overline X(\lambda)).
\]
In this subsection we introduce the quantities $d_r$ and $d_{r, r+1}$ that measure how much $w$ and $\overline w$ can deviate from each other when $w(\lambda)$ passes through $Z^m_r$ and $Z_{r, r+1}$, respectively.
Denote
\[
  h_*(s) = \max_{z \in B_{2c_z}(z_*)} h_s(z), \qquad
  \alpha_*(s) = \max_{z \in B_{2c_z}(z_*)} \alpha_s(z).
\]
For resonant zones, Lemma~\ref{l:est-res-cross} and Lemma~\ref{l:est-weak-res} provide estimates for the time of resonant crossing, let us now estimate how much $w$ can deviate from $\overline w$ using that $\dot w$ and $\dot{\overline w}$ are bounded.
\begin{lemma} \label{l:through-res}
  For some $C_{res} > 1$ the following holds.
  Suppose that $X_{init} \in \mathcal A_3 \setminus \mathcal E$ and
  at some time $\lambda_0$ with
  \[
    \varepsilon(\lambda_0 - \lambda_{init}) \le \Lambda, \qquad
    z(\lambda_0) \in B_{2c_z}(z_*), \qquad
    \omega(\lambda_0) \le \omega_0
  \]
  $w(\lambda)$ reaches the border of $Z^m(s)$ given by
  $I = \hat I(z) + D'_p(s, z) \alpha(z)$. We assume $s_1 < C_{s_1}\ln^2 \varepsilon$.

  Then $w(\lambda)$ exits $Z^m(s)$ at some time $\lambda_1 > \lambda_0$ via the other border $I = \hat I(z) - D'_p(s, z) \alpha(z)$ and we have the estimate
  \begin{equation}
    \norm{w(\lambda_1) - \overline w(\lambda_1)} \le
    \norm{w(\lambda_0) - \overline w(\lambda_0)}
    + d(s),
  \end{equation}
  where
  \begin{equation} \label{e:d-low}
    d(s) = C_{res} \alpha_*(s)|\ln \varepsilon| \sqrt{s_1/s_2} \qquad \text{for low-numerator resonances}
  \end{equation}
  and
  \begin{equation} \label{e:d-high}
    d(s) = C_{res} \sqrt{b_s} \alpha_*(s)\sqrt{s_1/s_2}
    + C_{res} \varepsilon |\ln^5 \varepsilon| \qquad \text{for high-numerator resonances}.
  \end{equation}
\end{lemma}
\begin{proof}
  Low-numerator resonances.
  By Lemma~\ref{l:est-res-cross} such $\lambda_1$ exists and
  \[
    \varepsilon(\lambda_1 - \lambda_0) \lesssim \alpha_*(s) |\ln \varepsilon| \sqrt{s_2/s_1}.
  \]
  As $\dot w$ and $\dot{\overline w}$ are bounded by $O(\varepsilon \hat \omega^{-1}) = O(\varepsilon s_1/s_2)$ (we use Lemma~\ref{l:d-f-0-I} to estimate the rate of change of $w$), this gives the estimate~\eqref{e:d-low}.

  High-numerator resonances. We argue in the same way as for low-numerator resonances, but apply Lemma~\ref{l:est-weak-res} instead. We have (we use the formula for $D_p'(s, z)$ from Lemma~\ref{l:3-res-zones} and the estimates $\alpha(z_0) \sim \sqrt{\varepsilon \hat h(z_0) \ln^3 \hat h(z_0)}$ and $s_1 \lesssim \ln^2 \varepsilon$)
  \[
    \varepsilon (\lambda_1 - \lambda_0) \lesssim
    D_p'(s, z_0) \alpha(z(\lambda_0)) \hat \omega^{0.5} + \varepsilon s_1 \lesssim \sqrt{b \hat \omega} \alpha_* + \varepsilon |\ln^5 \varepsilon|.
  \]
\end{proof}
\noindent We will use the notation $d_r = d(s_r)$ for shorthand.

Passing non-resonant zones is described by Lemma~\ref{l:nonres}.
Denote by
\begin{equation}
  d_{r, r+1} = C_{nonres} \sqrt{\varepsilon h_*(s_r)} \abs{\ln h_*(s_r)}(\sqrt{b_{s_r}} + \sqrt{b_{s_{r+1}}})
  + C_{nonres} \varepsilon |\ln^{\gamma_1} \varepsilon|.
\end{equation}
Here we pick large enough $C_{nonres}$ so that the last term in~\eqref{e:nonres-est} is bounded by $d_{r, r+1}$ (as long as $z(w(\lambda))$ stays in $B_{2c_z}(z_*)$). This is possible, as for fixed $z$ by Lemma~\ref{l:est-xi-r-r+1} the values of $h$ in $Z_{r, r+1}$ are close to each other, so we have in the last term of~\eqref{e:nonres-est} $h \sim \hat h_{s_r}(z) \lesssim h_*(s_r)$.

\subsection{Estimating sums over resonant vectors $s$} \label{ss:est-sum}
\begin{lemma} \label{l:est-sum}
  For any $a, b \in \mathbb R$ and $c, d > 0$ we have (sums are taken for all $s = (s_1, s_2) \in \mathbb N^2$)
  \begin{align}
  \begin{split}
    &\sum_{s: s_2 \le S_2} h_*(s)^c s_1^a \lesssim 1, \qquad
    \sum_s h_*(s)^c b_s^d s_1^a s_2^b \lesssim 1, \\
    &\sum_{s: s_2 \le S_2} \alpha_*(s) s_1^a \lesssim \sqrt{\varepsilon}, \qquad
    \sum_s b_s^d \alpha_*(s) s_1^a s_2^b \lesssim \sqrt{\varepsilon}.
  \end{split}
  \end{align}
\end{lemma}
\begin{proof}
  The period $T$ depends on $h$ and $z$: $T = -a(h, z) \ln h + b(h, z)$ by Lemma~\ref{l:period-cont}.
  From this we have
  \begin{equation}
    \ln h = -a^{-1}(T - b),
    \qquad
    \ln \hat h_s(z) = -a^{-1}(2\pi \hat \omega^{-1} - b).
  \end{equation}
  Thus for some $c_1 > 0$ we have
  \begin{equation} \label{e:est-h-exp-T}
    h_*(s) \lesssim e^{- c_1 \hat \omega^{-1}}.
  \end{equation}
  If $s_2 \le S_2$, we also have
  \[
    h_*(s) \lesssim e^{-c_2 s_1},
  \]
  where $c_2 = c_1 / S_2$. This implies the upper left estimate.

  As for any $f$ we have
  \[
    h_*(s)^{c/2} (s_1/s_2)^f \sim h_*(s)^{c/2} |\ln^f h_*(s)| \lesssim 1, \qquad
    b_s^{d/2} s_2^f = e^{-d C_F s_2/2} s_2^f \lesssim 1,
  \]
  the upper right estimate follows from
  \begin{equation} \label{e:tmp-5839}
    \sum_s h_*(s)^{c/2} b_{s}^{d/2} = O(1).
  \end{equation}
  Let us now prove this estimate.
  Recall the notation $\abs{s} = s_1 + s_2$. Either $\hat \omega^{-1} = s_1/s_2 \gtrsim \sqrt{\abs{s}}$ or $s_2 \gtrsim \sqrt{\abs{s}}$.
  But $h$ decreases exponentially with the growth of $\hat \omega^{-1}$ by~\eqref{e:est-h-exp-T} and $b_s$ decreases exponentially with the growth of $s_2$. Therefore, we have $h_*(s)^{c/2} b_s^{d/2} = O(e^{-c_3 \sqrt{\abs{s}}})$ for some $c_3$, which implies the required estimate~\eqref{e:tmp-5839}.

  Finally, the lower estimates follow from the upper estimates if we take into account
  \[
    \alpha \sim \sqrt{\varepsilon h_*(s) \ln^3 h_*(s)} \sim \sqrt{\varepsilon h_*(s) s_1^3 s_2^{-3}}.
  \]
\end{proof}

By Lemma~\ref{l:est-sum} we have
$m(\mathcal E_1) \le \sum m(\mathcal E_s) = O(\sqrt \varepsilon)$.
As the width of low-numerator resonant zones in $h$ is $\sim \sqrt{\varepsilon h_s}$, the width in $I$ is $\sim \sqrt{\varepsilon h_s}s_1/s_2$ and thus $m(Z^m(s)) = O(\sqrt{\varepsilon h_s}s_1)$. By Lemma~\ref{l:est-sum} we have
$m(\mathcal E_2) \le \sum m(Z^m(s)) = O(\sqrt \varepsilon)$.
This implies
\begin{equation}
  m(\mathcal E) = O(\sqrt \varepsilon).
\end{equation}

\begin{lemma} \label{l:sum-res}
\begin{equation}
  \sum_r d_r \lesssim \sqrt{\varepsilon} |\ln \varepsilon|.
\end{equation}
\end{lemma}
\begin{proof}
  The sum over low-numerator resonances is $O(\sqrt{\varepsilon} |\ln \varepsilon|)$ by the bottom left estimate of Lemma~\ref{l:est-sum}.
  The sum of the terms $C_{res} \sqrt{b_s} \alpha_*(s)\sqrt{s_1/s_2}$ over high-numerator resonances is $O(\sqrt{\varepsilon})$ by the bottom right estimate of Lemma~\ref{l:est-sum}. Finally, as the number of resonances is bounded by some power of $\ln \varepsilon$, the sum of the terms $C_{res} \varepsilon |\ln^5 \varepsilon|$ is also $O(\sqrt{\varepsilon})$.
\end{proof}
\begin{lemma} \label{l:sum-nonres}
  \begin{equation}
    \sum_r d_{r, r+1} = O(\sqrt{\varepsilon}).
  \end{equation}
\end{lemma}
\begin{proof}
  Firstly, the sum of the terms $\varepsilon |\ln^{\gamma_1} \varepsilon|$ is $O(\sqrt{\varepsilon})$, as the total number of resonances is bounded by some power of $\ln \varepsilon$.
  Secondly, we have $h_*(r+1) \sim h_*(r)$ by Lemma~\ref{l:est-xi-r-r+1}, thus it is enough to prove
  \[
    \sum_r \sqrt{h_*(r)} \abs{\ln h_*(r)} \sqrt{b_{s_r}} = O(1).
  \]
  This estimate follows from Lemma~\ref{l:est-sum}, as $|\ln h_*(r)| \sim (s_1/s_2)|_{s=s_r}$.
\end{proof}

\subsection{End of the proof}
\begin{lemma} \label{l:nonres-to-sqrt}
  For any large enough $C > 0$ there is $\lambda_f$ such that
  $h(\lambda_f) = C \sqrt{\varepsilon}$.
  For all $\lambda \in [\lambda_0, \lambda_f]$ we have
  $h(\lambda) \ge C \sqrt{\varepsilon}$
  and
  \begin{equation} \label{e:close}
    \norm{w(\lambda) - \overline w(\lambda)} <  O(\sqrt{\varepsilon} |\ln \varepsilon|).
  \end{equation}
\end{lemma}
\begin{proof}
  Set $\lambda_f$ equal to the first moment such that $h(\lambda_f) = C \sqrt{\varepsilon}$ or
  $|h(\lambda_f) - \overline h(\lambda_f)| > 0.5 C \sqrt{\varepsilon}$ or
  $\norm{z - z_*} = 2c_z$ or $\norm{w(\lambda_f) - w(h{=}0, z_*)} = k$ or $\lambda_f = \lambda_{init} + \varepsilon^{-1} \Lambda$.
  For $\lambda < \lambda_f$ we can apply the last part of Lemma~\ref{l:nonres} and we have $\overline h \sim h$.

  Consider the trajectory $w(\lambda)$ for $\lambda \le \lambda_f$.
  We will apply Lemma~\ref{l:through-res} to cover the moments from the time $w(\lambda)$ first enters $Z^m_r$ until $w(\lambda)$ reaches the border of $Z^m_r$ (with $\omega < \hat \omega(s_r)$). We will apply Lemma~\ref{l:nonres} from this moment until $w(\lambda)$ first reaches $Z^{m+1}_r$.
  Let us renumerate the resonances in such a way that
  the point $w(\lambda_0)$ is either in non-resonant zone $Z_{12}$ or in high-numerator middle resonant zone $Z^m_1$ (it is not in low-numerator resonant zone, as low-numerator resonant zones lie in $\mathcal E$). In the latter case, by Lemma~\ref{l:through-res} this solution enters the non-resonant zone $Z_{12}$ and the difference between $w(\lambda)$ and $\overline w(\lambda)$ accumulated is $ \le d_1$.
  Thus we can assume that $w(\lambda)$ starts in $Z_{12}$ at $t=0$ with
  \[
    \delta_1 = \norm{w - \overline w} \le d_1 + O(\sqrt{\varepsilon}|\ln \varepsilon|).
  \]
  Then $w(\lambda)$ passes $Z_{12}$ and approaches the second resonance zone $Z^m_2$, on entering this zone by Lemma~\ref{l:nonres} we have
  \[
    \norm{w - \overline w}
    \le
    \delta_1
    e^{C_6(
      \xi_1 - \xi_2
    )}
    + d_{1,2}.
  \]
  On leaving $Z^m_2$ into $Z_{2, 3}$ we have
  \[
    \delta_2 = \norm{w - \overline w}
    \le
    \delta_1
    e^{C_6(
      \xi_1 - \xi_2
    )}
    + d_{1,2} + d_2.
  \]
  After passing $Z_{2, 3}$ and $Z^m_3$, on exit from $Z^m_3$ we have
  \[
    \delta_3 = \norm{w - \overline w}
    \le
    \delta_1
    e^{C_6(
      \xi_1 - \xi_3
    )}
    +
    (d_{1, 2} + d_2)
    e^{C_6(
      \xi_2 - \xi_3
    )}
    +
    d_{2, 3} + d_3.
  \]
  Continuing like this, we get the estimate
  \begin{equation} \label{e:local-0285}
    \delta_k \le (O(\sqrt{\varepsilon}|\ln \varepsilon|) + d_1 + d_{1,2} + d_2 + \dots + d_{k-1, d_k} + d_k)
    \times
    e^{C_6(
      \xi_1 - \xi_k
    )}.
  \end{equation}
  As we have $\xi_i = O(1)$, by Lemma~\ref{l:sum-res} and Lemma~\ref{l:sum-nonres} this gives us
  \begin{equation}
    \delta_k \lesssim \sqrt{\varepsilon} |\ln \varepsilon|.
  \end{equation}
  This implies that we cannot have $\norm{z(\lambda_f) - z_*} = 2c_z$, as $\norm{z(\lambda_0) - z_*} \le c_z$.
  Similarly, we see that $\norm{w(\lambda_f) - w(h{=}0, z_*)} \ne k$, as near separatrices the solution of averaged system is close to $w(h{=}0, z_*)$.
  From~\eqref{e:local-0285} and $\norm{\pdv{h}{w}} = O(\ln^{-1} h)$ (this is proved in Lemma~\ref{l:d-f-0-I}) we have
  \[
    \abs{h(\lambda) - \overline h(\lambda)} \lesssim \sqrt{\varepsilon} |\ln \varepsilon|\abs{\ln^{-1} h(\lambda)}.
  \]
  This means that for large enough $C$ we cannot have $|h(\lambda_f) - \overline h(\lambda_f)| = 0.5 C \sqrt{\varepsilon}$. This also means $\overline h(\lambda_f) > 0$, thus we cannot have $\lambda_f = \lambda_0 + \varepsilon^{-1} \Lambda$. This leaves just one possibility: $h(\lambda_f) = C \sqrt{\varepsilon}$.
\end{proof}

\begin{remark}
  By Lemma~\ref{l:d-f-0-I} we have $\norm{\pdv{h}{w}} = O(\ln^{-1} h)$. Taking this into account, in Lemma~\ref{l:nonres-to-sqrt} we have
  \[
    \abs{h(\lambda_f) - \overline h(\lambda_f)} = O(\sqrt{\varepsilon}),
    \qquad
    \norm{z(\lambda_f) - \overline z(\lambda_f)} = O(\sqrt{\varepsilon} \abs{\ln \varepsilon}).
  \]
\end{remark}

\begin{lemma} \label{l:nonres-after-sqrt}
  Assume that for some $\lambda_f$ we have
  \[
    h(\lambda_f) = C \sqrt{\varepsilon}, \qquad
    \norm{h(\lambda) - \overline h(\lambda)} = O(\sqrt{\varepsilon}), \qquad
    \norm{z(\lambda) - \overline z(\lambda)} = O(\sqrt{\varepsilon} |\ln \varepsilon|).
  \]
  Then there is
  $\lambda_\Pi > \lambda_f$ with $h(w(\lambda_\Pi)) = 2 \varepsilon \ln^\gamma \varepsilon$
  such that for all $\lambda \in [\lambda_f, \lambda_\Pi]$
  \begin{align} \label{e:nonres-after-sqrt}
  \begin{split}
    \abs{h(\lambda) - \overline h(\lambda)} < O(\sqrt{\varepsilon}), \qquad
    \norm{z(\lambda) - \overline z(\lambda)} < O(\sqrt{\varepsilon} |\ln \varepsilon|).
  \end{split}
  \end{align}
\end{lemma}
\begin{proof}
  Denote $v = (h, z)$, $v(\lambda) = v(w(\lambda))$ and
  $\overline v(\lambda) = v(\overline w(\lambda))$.
  Set $\lambda_\Pi > \lambda_f$ be the first moment such that
  $v(\lambda_\Pi) \in \partial \Pi$ or
  $\norm{z(\lambda_\Pi) - z_*} = 2c_z$
  or $\norm{w(\lambda_\Pi) - w(h{=}0, z_*)} = k$
  or $\lambda_\Pi = \lambda_{init} + \varepsilon^{-1} \Lambda$.
  Consider $v(\lambda)$ for $\lambda \in [\lambda_f, \lambda_\Pi]$.

  The solution $v(\lambda)$ subsequently passes non-resonant and resonant zones, and the behaviour in these zones is described by Lemma~\ref{l:nonres} and Lemma~\ref{l:through-res}, respectively.
  The total time $\lambda_{tot}$ between $\lambda_f$ and $\lambda_\Pi$ is split into the time $\lambda_{res}$ spent in resonant zones and the time $\lambda_{nonres}$ spent in non-resonant zones.
  By Lemma~\ref{l:sum-res} we have
  $\varepsilon \lambda_{res} \lesssim \sqrt{\varepsilon} \abs{\ln \varepsilon}$.
  By Lemma~\ref{l:nonres} we have
  $\varepsilon \lambda_{nonres} \lesssim h_{max} \ln^2 \varepsilon$, where
  $h_{max}$ is the maximum value of $h(\lambda)$ for $\lambda \in [\lambda_f, \lambda_\Pi]$. Indeed, in our domain we have $\ln h \sim \ln \varepsilon$ and the sum of all terms $\xi_r - \xi_{r+1}$ is $O(\ln^{-1} \varepsilon)$.
  As
  $h_{max} \lesssim \sqrt{\varepsilon}$,
  we have
  \begin{equation}
    \varepsilon \lambda_{tot} \lesssim \sqrt{\varepsilon} \ln^{2} \varepsilon.
  \end{equation}
  As $\norm{\dot v}, \norm{\dot{\overline v}} \lesssim \varepsilon$, for $t \in [\lambda_f, \lambda_\Pi]$ all the values of $v(\lambda)$ and $\overline v(\lambda)$ differ from each other for different $\lambda$ by at most $O(\sqrt{\varepsilon} \ln^{2} \varepsilon)$. Thus we cannot have
  \[
    \norm{z(\lambda_\Pi) - z_*} = 2c_z \text{ or }
    \norm{w(\lambda_\Pi) - w(h{=}0, z_*)} = k \text{ or }
    \lambda_\Pi = \lambda_{init} + \varepsilon^{-1} \Lambda
  \]
  (in the last case we would have $\overline h(\lambda_\Pi) = -c_h$ and this is impossible). The remaining possibility is $h(\lambda_\Pi) = 2 \varepsilon \ln^\gamma \varepsilon$.

  The estimate on the change of $z(\lambda)$ allows us to show that $h(\lambda)$ decreases exponentially with the growth of $T(v(\lambda))$.
  The period $T$ depends on $h$ and $z$: $T = -a(h, z) \ln h + b(h, z)$ by Lemma~\ref{l:period-cont}.
  From this we have
  \begin{equation}
    \ln h = -a^{-1}(T - b).
  \end{equation}
  As $h$ and $z$ change by at most $O(\sqrt{\varepsilon} \abs{\ln^{2} \varepsilon})$,
  for some constants $a_0, b_0$ we have
  \[
    a(h, z) = a_0 + O(\sqrt{\varepsilon} \abs{\ln^{2} \varepsilon}),
    \quad
    b(h, z) = b_0 + O(\sqrt{\varepsilon} \abs{\ln^{2} \varepsilon}),
  \]
  so
  \[
    h = (1 + O(\sqrt{\varepsilon} \abs{\ln^{3} \varepsilon})) e^{-a_0^{-1}(T - b_0)}.
  \]

  Now we can get a better estimate for $\lambda_{nonres}$. Denote $\omega_r = \hat \omega(s_r)$, $T_r = 2\pi / \omega_r$, let $h_r$ be the value of $h$ when entering $Z_{r, r+1}$.
  By Lemma~\ref{l:nonres} we have
  \[
    \varepsilon \lambda_{nonres}
    \lesssim
    \sum_{r > r_0} (\omega_r - \omega_{r+1}) h_r \ln^3 h_r
    =
    (2 \pi^{-1}) \sum_{r > r_0} \omega_r \omega_{r+1}(T_{r+1} - T_r) h_r \ln^3 h_r
    \sim
   \sum_{r > r_0}  h_r \ln h_r (T_{r+1} - T_r).
  \]
  It is easy to see that $T_{r+1} - T_r$ are bounded, e.g. we can take resonances $\omega = 1/n$, then $T = 2 \pi n$ (up to $T = 2 \pi N \sim \ln^2 \varepsilon$).
  As $h(T)$ decreases exponentially, we have
  $\varepsilon \lambda_{nonres} \lesssim h_{r_0 + 1} \abs{\ln h_{r_0 + 1}}
  \lesssim h_{max} \abs{\ln \varepsilon}
  \sim \sqrt{\varepsilon} \abs{\ln \varepsilon}$. So we have
  \begin{equation} \label{e:t-tot}
    \varepsilon \lambda_{tot} \lesssim \sqrt{\varepsilon} |\ln \varepsilon|.
  \end{equation}

  We have $\dot{\overline h} \lesssim \varepsilon \abs{\ln^{-1} \varepsilon}$.
  Similarly, for any $\lambda_1, \lambda_2 \in [\lambda_f, \lambda_\Pi]$ we have
  \[
    \int_{\lambda_1}^{\lambda_2} \dot h d\lambda \lesssim (\lambda_2 - \lambda_1) \varepsilon \abs{\ln^{-1} \varepsilon} + O(\varepsilon),
  \]
  as the integral of $\dot h$ during one period is $O(\varepsilon)$ and $T \sim |\ln \varepsilon|$. Hence,
  \[
    \abs{\overline h(\lambda) - h(\lambda)}
    <
    \abs{\overline h(\lambda_f) - h(\lambda_f)} + O(\varepsilon \abs{\ln^{-1} \varepsilon}) \lambda_{tot} + O(\varepsilon)
    \le
    O(\sqrt{\varepsilon}).
  \]
  As $\dot z, \dot{\overline z} = O(\varepsilon)$, we have
  \[
    \norm{\overline z(\lambda) - z(\lambda)}
    <
    \norm{\overline z(\lambda_f) - z(\lambda_f)} + O(\varepsilon) \lambda_{tot}
    \le
    O(\sqrt{\varepsilon} \abs{\ln \varepsilon}).
  \]
\end{proof}
Together Lemma~\ref{l:nonres-to-sqrt} and Lemma~\ref{l:nonres-after-sqrt} imply Lemma~\ref{l:approach-sep}.

\newpage
\section{Crossing non-resonant zone: proof} \label{ss:nonres_proof}
This section is devoted to the proof of Lemma~\ref{l:nonres}.
Let us define the function
\begin{equation}
  u(w, \varphi, \lambda) = \sum_{1 \le |m| \le N, m \in \mathbb Z^2}
  \frac{f_m e^{i(m_1\varphi + m_2 \lambda)}}{i(m_1\omega(w) + m_2)}.
\end{equation}
Given $w$, let
$\xi_r = \frac{s_2}{s_1}$
be the nearest resonance to $\omega(w)$. Let us also denote
$\Delta = |\omega(w) - \xi_r|$.
\begin{lemma} \label{l:est-u}
  Denote
  \begin{equation} \label{e:as-bs}
    a_s(w) = exp\bigg( -C_F s_2 - C_F \frac{s_1}{T(w)} \bigg),
    \qquad
    b_s = exp\big( -C_F s_2 \big),
    \qquad
    0 < a_s < b_s < 1.
  \end{equation}
  Then inside $\Pi$ we have (provided that $s_2/s_1$ is the nearest resonance to $\omega(w)$)
  \begin{align}
  \begin{split}
    &\norm{u} \lesssim a_s s_1^{-1} \Delta^{-1} + \ln^2 \varepsilon,\\
    &\norm{\pdv{u}{\lambda}} \lesssim a_s s_2 s_1^{-1} \Delta^{-1} + \ln^2 \varepsilon, \qquad
    \norm{\pdv{u}{\varphi}} \lesssim a_s \Delta^{-1}  + |\ln^3 \varepsilon|, \\
    &\norm{\pdv{u}{h}} \lesssim
    b_s s_1^{-1} h^{-1} \ln^{-2} h \; \Delta^{-1}(\Delta^{-1} + \abs{\ln h})
    + h^{-1} \ln^{-1} h \ln^4 \varepsilon, \\
    &\norm{\pdv{u}{z}} \lesssim
    b_s s_1^{-1} |\ln^{-1} h| \; \Delta^{-1}(\Delta^{-1} + \abs{\ln h})
    + \ln^4 \varepsilon. \\
  \end{split}
  \end{align}
\end{lemma}
\noindent This lemma is proved in Appendix~\ref{a:appendix}.

Take small $\tilde k > k$ (here $k$ is from Section~\ref{s:res-zones}) and denote by $\tilde{\mathcal B}_{3, *} \subset \mathcal B_3$ the set of points $(h, z)$ that satisfy $\norm{w(h, z) - w(h{=}0, z_*)} < \tilde k$. The value of $\tilde k$ is picked so that $\tilde{\mathcal B}_{3, *}$ is far from $\partial \mathcal B$.
Let us consider the sets
\begin{equation} \label{e:Pi-tilde}
  \tilde \Pi = \{ (h, z) \in \tilde{\mathcal B}_{3, *}: h \ge \varepsilon |\ln^{\gamma} \varepsilon| \},
  \qquad
  \partial \tilde \Pi = \{ (h, z) \in \tilde{\mathcal B}_{3, *}: h = \varepsilon |\ln^{\gamma} \varepsilon| \}.
\end{equation}
where $\gamma > 2$ is the same as in~\eqref{e:Pi}. We have
$\Pi \subset \tilde \Pi$.
The (large enough) constant $C_{\tilde Z} \in (0, C_Z)$ will be chosen later in Lemma~\ref{l:beta}.
For each resonance $\xi_r = s_2/s_1$ define the zone $\tilde Z_r$ by the condition
\begin{equation} \label{e:tilde-Z-r}
  \tilde Z_r = \Big\{ (\tilde h, \tilde z) \in \tilde \Pi:
  |\omega(\tilde h, \tilde z) - \xi_r|
  \le
  C_{\tilde Z} \delta_s(\tilde z) \Big\}.
\end{equation}
Denote by $\tilde Z_{r, r+1}$
the zone between two neighboring zones $\tilde Z_r$ and $\tilde Z_{r+1}$, $\xi_r > \xi_{r+1}$. We have $Z_{r, r+1} \subset \tilde Z_{r, r+1}$.
As the zones $\tilde Z$ differ from $Z$ only in the values of constants, the properties of the zones $Z$ discussed in Section~\ref{s:res-zones} also hold for the zones $\tilde Z$ for large enough $C_{\tilde Z}$.
Let us also note that in $\tilde Z_{r, r+1}$ for $s$ corresponding to the nearest resonance ($\xi_r$ or $\xi_{r+1}$) we have (recall that by Lemma~\ref{l:est-xi-r-r+1} we have $\hat h_{\xi_r}, \hat h_{\xi_{r+1}} \in [0.9 h, 1.1 h]$)
\begin{equation} \label{e:Delta}
  \Delta^{-1} \le 2 C_{\tilde Z}^{-1} \sqrt{h \ln^4 h /(\varepsilon b_s)},
  \qquad
  \Delta^{-1} \le 2 C_{\tilde Z}^{-1} \varepsilon^{-1} h \ln^{-3} h \ln^{-2} \varepsilon.
\end{equation}

\begin{lemma} \label{l:invertible}
  There is $C_u > 0$ such that for large enough $C_{\tilde Z}$ the following estimates hold
  in $\tilde Z_{r, r+1}$ for all $r$ (with $s$ corresponding to the nearest to $\omega(w)$ resonance)
  \begin{itemize}
    \item
    $\norm{\varepsilon u}
    < C_u C_{\tilde Z}^{-1} \sqrt{b_s \varepsilon h \ln^2 h} + C_u \varepsilon \ln^2 \varepsilon
    < \min(0.1 h, |\ln^{-1} \varepsilon|)$,
    \item $\norm{\pdv{u}{w}}, \norm{\pdv{u}{\varphi}}, \norm{\pdv{u}{\lambda}} \le 0.1 \varepsilon^{-1}$.
  \end{itemize}
  Moreover, the coordinate change
  \begin{equation} \label{e:varchange}
    U : \tilde w, \varphi, \lambda \mapsto \tilde w + \varepsilon u(\tilde w, \varphi, \lambda), \varphi, \lambda
  \end{equation}
  is invertible in ${(\tilde Z_{r, r+1}) \times [0, 2\pi]^2}$.
\end{lemma}
\begin{proof}
  First, let us prove that inside $\tilde Z_{r, r+1}$ we have (for $s$ corresponding to $\xi_r$ or $\xi_{r+1}$)
  \begin{equation} \label{e:s1}
    s_1^{-1} =  O(\ln^{-1} h).
  \end{equation}
  By~\eqref{e:est-xi-r-r+1} we have $\xi_r, \xi_{r+1} \sim - \ln^{-1} h$, so
  $s_1^{-1} = s_2^{-1} \xi_{r_1} = O(\ln^{-1} h)$, where $r_1 = r$ or $r_1 = r+1$.

  The estimates on $u$ and its derivatives follow from Lemma~\ref{l:est-u} and from~\eqref{e:Delta},~\eqref{e:s1} and~\eqref{e:in-Pi}. To estimate the $I$-derivative, we use $\pdv{u}{I} = \omega \pdv{u}{h}$.

  Denote $Y = (\tilde w, \varphi, \lambda)$. From the estimates on the derivatives of $u$ we have $\varepsilon \norm{\pdv{u}{Y}} < 0.5$, so the coordinate change $U$ is invertible.
\end{proof}

\begin{lemma} \label{l:uh}
  Take some $\tilde w \in \tilde Z_{r, r+1}$, $\varphi$, $\lambda$.
  Set $w = \tilde w + \varepsilon u(\tilde w, \varphi, \lambda)$.
  Let $h = h(w)$, $\tilde h = h(\tilde w)$.
  Then we have $\norm{h - \tilde h} < \min(0.25 \tilde h, |\ln^{-1} \varepsilon|)$.
\end{lemma}
\begin{proof}
  For $\alpha \in [0, 1]$ let us denote $w_\alpha = \tilde w + \alpha \varepsilon u(\tilde w, \varphi, \lambda)$, $h_\alpha = h(w_\alpha)$.
  Take largest possible $\alpha'$ such that we have $|h_\alpha - \tilde h| \le 0.25 \tilde h$ for all $\alpha \le \alpha'$. We have
  $h_{\alpha'} - \tilde h = \pdv{h}{w}\big|_\xi \alpha' u$
  with some $\xi \in [\tilde w, \tilde w + \alpha' u]$. As we have $h \in [0.75 \tilde h, 1.25 \tilde h]$ in $[\tilde w, \tilde w + \alpha' u]$, Lemma~\ref{l:d-f-0-I} gives $|h_{\alpha'} - \tilde h| = O(\ln^{-1} h) u < O(\ln^{-1} \tilde h) \min(0.1 \tilde h, |\ln^{-1} \varepsilon|)$
  close enough to the separatrices by Lemma~\ref{l:invertible}.
  Thus we actually have $\alpha' = 1$ and the estimate for $h_{\alpha'}$ gives
  $|h(w) - \tilde h| < \min(0.25 \tilde h, |\ln^{-1} \varepsilon|)$.
\end{proof}

\noindent Recall that the zones $Z_r$ depend on the constant $C_Z$.
\begin{lemma} \label{l:cover}
  Given $C_{\tilde Z}$, there is $C_{Z, 0} > 0$ such that for any $C_Z > C_{Z, 0}$ for all sufficiently small $\varepsilon$ and we have the following.
  Take the map $U$ defined by~\eqref{e:varchange}. Then for all $r$
  \[
    Z_{r, r+1} \times [0, 2\pi]^2 \subset U(\tilde Z_{r, r+1} \times [0, 2\pi]^2).
  \]
\end{lemma}
\begin{proof}
  For $w_0 \in \tilde Z_{r, r+1}$ denote $h_0 = h(w_0)$.
  Set
  \[
    R(w_0) = \varepsilon \; max_{\varphi, \lambda} \norm{u(w_0, \varphi, \lambda)},
    \qquad
    B(w_0) = \{ w \in \tilde Z_{r, r+1} : \norm{w - w_0} \le R(w_0)\}.
  \]
  Let us note that by Lemma~\ref{l:invertible} we have
  \begin{equation} \label{e:R-to-zero}
    \sup_{w_0 \in \cup_r \tilde Z_{r, r+1}} R(w_0) \to 0 \text{ for } \varepsilon \to 0.
  \end{equation}
  By Lemma~\ref{l:uh} we also have
  \begin{equation} \label{e:R-h}
    |h(w) - h_0| \le 0.25 h_0 \qquad \text{for } w \in B(w_0).
  \end{equation}

  As we have
  $Z_{r, r+1} \subset \tilde Z_{r, r+1}$ (due to $C_Z > C_{\tilde Z}$),
  it is enough to prove that for any
  $w_0 \in \partial \tilde Z_{r, r+1}$
  the ball $B(w_0)$ does not intersect $Z_{r, r+1}$.

  First, we may have $w_0 \in \partial \tilde{\mathcal B}_{3, *}$ with
  $\norm{w_0 - w(h{=}0, z_*)} = \tilde k$.
  Then for $w \in B(w_0)$ we have
  $\norm{w - w(h{=}0, z_*)} = \tilde k + o(1) > k$
  for small $\varepsilon$, thus the ball $B(w_0)$ does not intersect $\mathcal B_{3, *} \supset Z_{r, r+1}$.

  Second, we may have $w_0 \in \partial \tilde \Pi$
  with
  $h_0 = \varepsilon |\ln^{\gamma} \varepsilon|$.
  By~\eqref{e:R-h}, for any $(h, z) \in B(w_0)$ we have
  $h \le 1.5 \varepsilon |\ln^{\gamma} \varepsilon|$, so $B(w_0)$ does not intersect $\Pi \supset Z_{r, r+1}$.

  Finally, we may have $|\omega(w_0) - \xi_r| = C_{\tilde Z} \hat \delta_{\xi_r}(z_0)$ (or the same for $r{+}1$ instead of $r$, this case is treated in the same way). We will write $\delta$ instead of $\delta_{\xi_r}$ for brevity.
  By~\eqref{e:R-h} we have $h \in [0.25 h_0, 1.25 h_0]$ in $B(w_0)$.
  By Lemma~\ref{l:est-xi-r-r+1} this implies
  $\hat h(z(w)) \in [0.5 h_0, 1.5 h_0]$
  in $B(w_0)$.
  So in $\lesssim$-estimates below we will write $h$ for both $h(w)$ and $\hat h(z(w))$ for $w \in B(w_0)$.
  This also means
  $\delta(z) \in [0.5\delta(z_0), 2\delta(z_0)]$
  in $B(w_0)$.

  Denote $\Delta \omega = \abs{\omega(w) - \omega(w_0)}$. It is enough to prove that there is a constant $c_1 > 0$ (that does not depend on $r$ or $w_0$) such that $\Delta \omega < c_1 \delta(z_0)$, then we can take $C_{Z, 0} = 2(c_1 + C_{\tilde Z})$ and have
  \[
    \abs{\omega - \xi_r} < (c_1 + C_{\tilde Z}) \delta(z_0) < 2 (c_1 + C_{\tilde Z}) \delta(z) <  C_Z \delta(z)
  \]
  in $B(w_0)$. Hence, $B(w_0)$ does not intersect $Z_{r, r+1}$.

  We have $\Delta \omega \le \norm{\pdv{\omega}{w}\/(w_{int})} R(w_0)$ for some $w_{int} \in [w_0, w]$.
  We have $\norm{\pdv{\omega}{w}} \lesssim h_0^{-1} \ln^{-3} h_0$.
  We will use that
  \[
    R(w_0) \lesssim \varepsilon b_s s_1^{-1} \delta^{-1}(w_0) + \varepsilon \ln^2 \varepsilon
  \]
  by Lemma~\ref{l:est-u}.
  Hence we obtain using~\eqref{e:s1}
  \[
    \norm{\Delta \omega} \lesssim
    \varepsilon b_s h^{-1} \ln^{-4} h \; \delta^{-1}(w_0)
    + \varepsilon h^{-1} \ln^{-3} h \ln^2 \varepsilon.
  \]
  We have
  $
    \varepsilon h^{-1} \ln^{-3} h \ln^{2} \varepsilon
    \lesssim \delta(w_0).
  $
  We also have
  \[
    \varepsilon b_s h^{-1} \ln^{-4} h \; \delta^{-1}(w_0)
    \lesssim \delta(w_0),
  \]
  as
  $
    \delta(w_0) \ge \sqrt{ \varepsilon b_s h^{-1} \ln^{-4} h}
  $.
  This completes the proof.
\end{proof}

\begin{lemma} \label{l:beta}
  For some $C_\beta > 0$ and for large enough value of the constant $C_{\tilde Z}$ defined above the following holds.
  Inside $\tilde Z_{r, r+1}$ the change of variables
  $w = \tilde w + \varepsilon u(w, \varphi, \lambda)$
  takes the perturbed system~\eqref{e:perturbed-aa} to the form
  \begin{align} \label{e:dot-tilde-v}
  \begin{split}
    \dot {\tilde w} &= \varepsilon f_0(\tilde w) + \varepsilon^2 \beta, \\
    \dot \varphi &= \omega(\tilde w + \varepsilon u) + \varepsilon g(\tilde w + \varepsilon u, \varphi, \lambda, \varepsilon), \\
    \dot{\lambda} &= 1,
  \end{split}
  \end{align}
  where $\beta$ is a smooth function of $\tilde w, \varphi, \lambda$ that depends on $\varepsilon$
  and
  \begin{equation} \label{e:est-beta}
    \norm{\beta}
    = O(b_s \Delta^{-2} h^{-1} \ln^{-4} h)
    + O(b_s \Delta^{-1} h^{-1} \ln^{-4} h \ln^3 \varepsilon)
    + O(h^{-1} \ln^{-3} h \ln^5 \varepsilon)
  \end{equation}
  (here $s$ corresponds to the nearest to $\omega(\tilde w)$ resonance).
  Moreover, for small enough $\varepsilon$ we have
  \begin{equation} \label{e:est-beta-2}
    \varepsilon^2 \norm{\beta}
    \le C_\beta C_{\tilde Z}^{-2} \varepsilon
  \end{equation}
  and
  \begin{align} \label{e:decrease}
  \begin{split}
    &-1.5 \varepsilon K_h |\ln^{-1} h| < \dv{h(\tilde w)}{\lambda} < -0.5 \varepsilon K_h |\ln^{-1} h|,
    \\
    &-1.5 \varepsilon K_\omega h^{-1} |\ln^{-3} h| < \dv{\omega( \tilde{w} )}{\lambda} < -0.5 \varepsilon K_\omega h^{-1} |\ln^{-3} h|, \\
    &0.5 \abs{\dv{\omega( \tilde{w} )}{\lambda}} > C_{\tilde Z} \dv{\delta_{s_r}(z(\tilde w))}{\lambda},
  \end{split}
  \end{align}
  where $\dv{\tilde w}{\lambda}$ is given by~\eqref{e:dot-tilde-v}.
\end{lemma}
\begin{proof}
  From~\cite[Proof of Lemma $7.3$]{neishtadt2014} we have
  \[
    \varepsilon \beta = \bigg(
      \big(E + \varepsilon \pdv{u}{\tilde w} \big)^{-1} - E
      \bigg)
      (f_0(\tilde z) + \varepsilon \beta_1) + \varepsilon \beta_1
  \]
  with
  \begin{align}
  \begin{split}
    \varepsilon \beta_1 &= [f(w, \varphi, \lambda, \varepsilon) - f(w, \varphi, \lambda, 0)]
    + [f(\tilde w + \varepsilon u, \varphi, \lambda, 0) - f(\tilde w, \varphi, \lambda, 0)] \\
    &+ R_N f
    - \varepsilon \pdv{u}{\varphi} f_\varphi
    - \pdv{u}{\varphi}[\omega(\tilde w + \varepsilon u) - \omega(\tilde w)].
  \end{split}
  \end{align}
  Note that as $w \in \Pi$ we have $|\ln h| < |\ln \varepsilon|$.
  As $\omega$ is bounded, $1 \lesssim \Delta^{-1} \lesssim \Delta^{-2}$.
  For a vector function $g(w) = (g_1(w), \dots g_n(w))$ denote
  $(\pdv{g}{w})_{int} = (\pdv{g_1}{w}\/(\eta_1), \dots, \pdv{g_n}{w}\/(\eta_n))$,
  where $\eta_1, \dots, \eta_n \in [w, w + \varepsilon u]$ are some intermediate points.
  By Lemma~\ref{l:uh} the values of $h$ for the points $\eta_i$ are in $[0.5h(w), 2h(w)]$.
  We have the following estimates (using the estimates from Lemma~\ref{l:est-u}, Lemma~\ref{l:fourier-remainder},~\eqref{e:est-real} and~\eqref{e:s1}).
  \begin{align}
  \begin{split}
    &\varepsilon^{-1}\norm{ f(w, \varphi, \lambda, \varepsilon) - f(w, \varphi, \lambda, 0)}
    = O(\ln h), \\
    &\varepsilon^{-1}\norm{f(\tilde w + \varepsilon u, \varphi, \lambda, 0) - f(\tilde w, \varphi, \lambda, 0)}
    = \norm{\Big(\pdv{f}{w}\Big)_{int} u}
    = O(h^{-1} \ln^{-2} h) b_s \Delta^{-1} + O(h^{-1} \ln^{-1} h \ln^2 \varepsilon), \\
    &\varepsilon^{-1}\norm{R_N f} < 1, \\
    &\norm{\pdv{u}{\varphi} f_\varphi}
    = O(h^{-1} \ln^{-2} h) b_s \Delta^{-1} + O(h^{-1} \ln^{-2} h \ln^3 \varepsilon), \\
    & \varepsilon^{-1}\norm{\pdv{u}{\varphi}[\omega(\tilde w + \varepsilon u) - \omega(\tilde w)]}
    = \norm{\pdv{u}{\varphi} \Big(\pdv{\omega}{w}\Big)_{int} u} = \\
    &\qquad = O(h^{-1}\ln^{-4} h) \Big(b_s^2 \Delta^{-2} + b_s \Delta^{-1} |\ln^3 \varepsilon|\Big) + O(h^{-1}\ln^{-3} h \ln^5 \varepsilon).
  \end{split}
  \end{align}
  This gives (we use $|\ln \varepsilon| \gtrsim |\ln h|$ to simplify the expression below)
  \begin{equation}
    \norm{\beta_1}
    \le
    O(h^{-1} \ln^{-4} h) b_s^2 \Delta^{-2}
    + O(h^{-1} \ln^{-4} h \ln^3 \varepsilon)b_s \Delta^{-1}
    + O(h^{-1}\ln^{-3} h \ln^5 \varepsilon).
  \end{equation}
  As $\norm{\varepsilon \pdv{u}{w}} \le 0.5$ and $\norm{f_{w, 0}} = O(1)$,
  we have the estimate
  \[
    \norm{\beta} \lesssim \norm{\beta_1} + \norm{f_{0}} \norm{\pdv{u}{\tilde w}}
  \]
  As we have the estimate
  \[
    \norm{\pdv{u}{\tilde w}} \lesssim b_s h^{-1} \ln^{-4} h (\Delta^{-2} + \abs{\ln h} \Delta^{-1})
    + h^{-1} \ln^{-2} h \ln^4 \varepsilon,
  \]
  by Lemma~\ref{l:est-u} and~\eqref{e:s1} (we use $\pdv{u}{I} = \omega \pdv{u}{h}$), this implies~\eqref{e:est-beta}.

  From~\eqref{e:est-beta} and~\eqref{e:Delta} we have~\eqref{e:est-beta-2}, the estimate for the second and third terms of~\eqref{e:est-beta} uses~\eqref{e:in-Pi}. Here we use that $\gamma > 2$ and so
  $|\ln \varepsilon|^{2-\gamma} < C_{\tilde Z}^{-2}$
  for small $\varepsilon$.

  As $\pdv{h}{w} = O(\ln^{-1} h)$ by Lemma~\ref{l:d-f-0-I}, for large $C_{\tilde Z}$ we have $\norm{\pdv{h}{w}} C_{\beta} C_{\tilde Z}^{-2} < 0.5 K_h |\ln^{-1} h|$. Then~\eqref{e:est-beta-2} and~\eqref{e:K-h} imply the first part of~\eqref{e:decrease}.
  For some $c_1 > 0$ we have $\norm{\pdv{\omega}{w}} < c_1 h^{-1} \ln^{-3} h$ in $\tilde {\mathcal B}_{3, *}$. For large $C_{\tilde Z}$ we have $C_{\beta} C_{\tilde Z}^{-2} c_1 < 0.5 K_\omega$. This (together with~\eqref{e:est-beta-2} and~\eqref{e:K-omega}) means that the second part of~\eqref{e:decrease} also holds.

  We have $\pdv{\hat h_{\xi_r}}{z} = O(h \ln h)$ by~\eqref{e:est}.
  This implies $\pdv{\delta_{s_r}}{z} = O(\sqrt{\varepsilon/h} \ln^{3} \varepsilon)$.
  As $\dv{z}{\lambda} = O(\varepsilon)$, this means $\dv{\delta_{s_r}}{\lambda} = O(\varepsilon^{3/2} h^{-1/2} \ln^3 \varepsilon)$.
  By the estimate on $\dv{\omega}{\lambda}$ this implies the last part of~\eqref{e:decrease}.

\end{proof}

\begin{lemma} \label{l:h-const}
  There is a constant $c_1 > 0$ such that the following holds.
  Assume that a solution
  $(\tilde w(\lambda), \varphi(\lambda), \lambda(\lambda))$
  of~\eqref{e:dot-tilde-v} stays inside $\tilde Z_{r, r+1}$
  for $\lambda \in [\lambda_1, \lambda_2]$.
  Denote $\tilde h(\lambda) = h(\tilde w(\lambda))$.
  Then we have $\tilde h(\lambda_2) > c_1 \tilde h(\lambda_1)$.
\end{lemma}
\begin{proof}
  By~\eqref{e:decrease} both $h$ and $\omega$ decrease along our solution and
  $\dv{\omega}{h} \sim h^{-1} \ln^{-2} h$. This gives $\dv{}{h}\/(\omega^{-1}) \sim h^{-1}$.
  Afterwards this lemma is proved like Lemma~\ref{l:est-xi-r-r+1}, but we integrate the derivative $\dv{}{h}\/(\omega^{-1})$ along our solution instead of the partial derivative for fixed $z$.
\end{proof}

\noindent
Let us say that a domain $D \subset \mathbb R^m$ is $L$-approximately convex ($L \ge 1$), if any two points $w_1, w_2 \in D$ can be connected by a piecewise linear path with length at most $L \norm{w_1 - w_2}$ that lies in $D$.
In such domain any vector-function $\psi(w)$ satisfies the estimate
\[
  \norm{\psi(w_1) - \psi(w_2)} \le L \; \max_D \norm{\pdv{\psi}{w}} \; \norm{w_1 - w_2}.
\]
We will need the following lemma. It is well-known for convex domains (with $L=1$) and generalises straightforwardly for $L$-approximately convex domains, using the estimate above.
\begin{lemma} \label{l:ode-est}
  Consider two ODEs
  \begin{equation}
    \dot w_1 = a(w_1),
    \qquad
    \dot w_2 = a(w_2) + b(\lambda)
  \end{equation}
  defined in some $L$-convex domain $D$.
  Consider two solutions $w_1(\lambda), w_2(\lambda)$ with
  \[
    \norm{w_1(\lambda_0) - w_2(\lambda_0)} < \delta
  \]
  that exist and stay in $D$ up to the moment $T$.
  Assume that in $D$ we have the estimate $L \norm{\pdv{a}{w}} \le A$.
  Then for any $\lambda \in [\lambda_0, T)$ we have the estimate
  \begin{equation} \label{e:ode-est}
    \norm{w_2(\lambda) - w_1(\lambda)}
    \le e^{A(\lambda-\lambda_0)}
    \bigg(
      \delta + \int_{\tau = \lambda_0}^\lambda \norm{b(\tau)} d\tau
    \bigg).
  \end{equation}
\end{lemma}

\begin{lemma} \label{l:L-conv}
  There exists $L \ge 1$ such that for all $h_0 > 0$ the domain
  $\{w \in \tilde{\mathcal B}_{3, *} : h(w) > h_0\}$ is $L$-approximately convex.
\end{lemma}
\begin{proof}
  Let us build a path connecting $w_1$ with $w_2$. We assume $I(w_1) \le I(w_2)$.
  Let $I_{max}$ and $I_{min}$ be the maximum and minimum of $I(h_0, z)$, where $z$ lies in the segment connecting $z(w_1)$ and $z(w_2)$.
  If $I(w_1), I(w_2) > I_{max}$, we can connect $w_1$ and $w_2$ by a segment.
  Otherwise, we have $I(w_1) \le I_{max}$.
  Then we connect by a segment $w_1$ with $(I_{max}, z(w_1))$, then $(I_{max}, z(w_1))$ with $(I_{max}, z(w_2))$ and, finally, $(I_{max}, z(w_2))$ with $w_2$. The length of this path is bounded by
  $(I_{max} - I_{min}) + \norm{z_1 - z_2} + (I(w_2) - I(w_1))$ if $I(w_2) > I_{max}$ and by $2(I_{max} - I_{min}) + \norm{z_1 - z_2}$ otherwise.

  Let us prove that for some $L_0$ we have
  \begin{equation} \label{e:local-1055}
    I_{max} - I_{min} \le L_0 \norm{z_1 - z_2}.
  \end{equation}
  As $\pdv{I}{h} = \omega^{-1}$, for any $z', z'' \in [z(w_1), z(w_2)]$ we have
  \[
    |I(h_0, z') - I(h_0, z'')| \le |I(0, z') - I(0, z'')| + \in\lambda_0^{h_0} \norm{\pdv{}{z} \omega^{-1}} \norm{z' - z''} dh_1 \lesssim \norm{z' - z''}.
  \]
  This gives~\eqref{e:local-1055}. So our domain is approximately convex with $L = 2L_0 + 2$.
\end{proof}

\begin{lemma} \label{l:nonres-tilde}
  There exist $C_6, \tilde C_t > 0$ and $\gamma_1 \in \mathbb R$ such that the following statement holds. Assume that for some $\lambda_0 > 0$ and $\tilde w_0$ we have
  \[
    \tilde w_0 \in \tilde Z_{r, r+1}.
  \]
  a) Then for any $\varphi_0 \in [0, 2\pi]$ there exists $\lambda_1 > \lambda_0$ such that the solution $\tilde X(\lambda)$ of~\eqref{e:dot-tilde-v} with
  $\tilde X(\lambda_0) = (\tilde w_0, \varphi_0, \lambda_0)$
  is defined for all
  $\lambda \in [\lambda_0, \lambda_1]$.
  This solution satisfies
  $\tilde w(\lambda) \in \tilde Z_{r, r+1}$
  and $\tilde w(\lambda_1)$ lies on the boundary of
    $\tilde Z_{r, r+1}$:
    in $\partial \tilde \Pi$ or on the border with $\tilde Z_{r+1}$.

  \noindent b)
  Denote $\tilde h(\lambda) = h(\tilde w(\lambda))$.
  Then
  \begin{equation} \label{e:time-nonres-zone-tilde}
    \varepsilon(\lambda_1 - \lambda_0) \le \tilde C_t (\xi_r - \xi_{r+1}) \tilde h(\lambda_1) |\ln^3 \tilde h(\lambda_1)|.
  \end{equation}

  \noindent c)
  Denote $\overline h(\lambda) = h(\overline w(\lambda))$. Assume that for some $\lambda_{01} \in [\lambda_0, \lambda_1]$ we have
  \begin{equation} \label{l:ratio-h-overline-h}
    \overline h(\lambda_{01}) > 0.25 \tilde h(\lambda_{01}),
    \qquad
    \overline h(\lambda_0) < 3 \tilde h(\lambda_0).
  \end{equation}
  Then for all $\lambda \in [\lambda_0, \lambda_{01}]$ we have (we use the notation $h = h(\lambda_1)$ in the formula below)
  \begin{equation} \label{e:v-tilde-approx}
    \norm{\tilde w(\lambda) - \overline w(\lambda)} <
    e^{
      C_6 (\xi_r - \xi_{r+1})
    }
    \norm{\tilde w(\lambda_0) - \overline w(\lambda_0)}
    + O(\sqrt{\varepsilon h} \abs{\ln h})(\sqrt{b_{s_r}} + \sqrt{b_{s_{r+1}}})
    + O(\varepsilon |\ln^{\gamma_1} \varepsilon|).
  \end{equation}
\end{lemma}
\begin{proof}
  By~\eqref{e:decrease} the value of $\abs{\omega - \xi_r} - C_{\tilde Z_r} \delta_{s_r}$ increases with the time, so the solution $\tilde w(\lambda)$ does not cross the border of
  $\tilde Z_r$ and $\tilde Z_{r, r+1}$. This proves the first statement of the lemma.

  As $\omega$ decreases with the time, we can take $\omega$ as an independent variable. From~\eqref{e:decrease} we have
  \begin{equation} \label{e:d-omega-t}
    \dv{\lambda}{\omega} \sim \varepsilon^{-1} h \ln^3 h,
  \end{equation}
  so
  \begin{equation} \label{e:t1-t2}
    \lambda_1 - \lambda_0 = O(\varepsilon^{-1} h \ln^3 h)(\omega(\lambda_0) - \omega(\lambda_1)).
  \end{equation}
  As $\omega \in [\xi_{r+1}, \xi_r]$, this gives~\eqref{e:time-nonres-zone-tilde}.

  Let us now prove the estimate for $\norm{\tilde w(\lambda) - \overline w(\lambda)}$.
  Note that by Lemma~\ref{l:h-const} and~\eqref{l:ratio-h-overline-h}
  the values of $h(\lambda)$ for different $\lambda \in [\lambda_0, \lambda_{01}]$ differ at most by a constant factor between themselves and also with $h(\overline w(\lambda))$. Thus we will simply write $h$ in $O$-estimates (for definiteness, set $h = h(\lambda_1)$). Let us use Lemma~\ref{l:ode-est} with $D = \{ w \in G : h(w) > c_1 h \}$ for $c_1 > 0$ chosen to cover all the values of $h$ discussed above. This domain is $L$-convex for $L$ chosen in Lemma~\ref{l:L-conv}.
  We have
  $a = (\varepsilon f_{I, 0}, \varepsilon f_{z, 0})$.
  By Lemma~\ref{l:d-f-0-I} we have
  $L \norm{\pdv{a}{w}} \le A \sim \varepsilon h^{-1} \abs{\ln^{-3} h}$.
  By~\eqref{e:time-nonres-zone-tilde} for some $C_6 > 0$ we have
  \[
    A(\lambda-\lambda_0) \le C_6 (\xi_r - \xi_{r+1}).
  \]
  We have $b(\lambda) = \varepsilon^2 \beta$. For some $\gamma_1$ we get from~\eqref{e:est-beta} the following estimate for $\beta$
  \begin{equation}
    \norm{\beta}
    = O(b_s \Delta^{-2} h^{-1} \ln^{-4} h)
    + O(h^{-1} |\ln^{\gamma_1 - 4} \varepsilon|)(b_s \Delta^{-1} + 1).
  \end{equation}
  Denote
  $\omega_0 = \omega(\lambda_0)$, $\omega_1 = \omega(\lambda_1)$, we have $\omega_0 > \omega_1$.
  By~\eqref{e:d-omega-t} we have
  \[
    \int_{\tau = \lambda_0}^{\lambda} \norm{b(\tau)} d\tau
    \le \int_{\tau = \lambda_0}^{\lambda_1} \norm{b(\tau)} d\tau
    \lesssim \varepsilon h \abs{\ln^3 h} \int_{\omega_1}^{\omega_0} \norm{\beta} d\omega.
  \]
  We will use~\eqref{e:Delta} in the estimates for the integrals below.
  Clearly,
  \[
    \int_{\omega_1}^{\omega_0} b_s \Delta^{-2} d \omega
    \le b_{s_r} \Delta^{-1}(\omega_0) + b_{s_{r+1}} \Delta^{-1}(\omega_1)
    = O(\sqrt{h \abs{\ln^4 h}/\varepsilon}) (\sqrt{b_{s_r}} + \sqrt{b_{s_{r+1}}})
  \]
  and
  \[
    \int_{\omega_1}^{\omega_0} (b_s \Delta^{-1} + 1) d \omega
    \le (|\ln \Delta(\omega_0)| + |\ln \Delta(\omega_1)|) + (\omega_0 - \omega_1)
    = O(\ln \varepsilon).
  \]
  By the estimate on $\norm{\beta}$ this yields
  \[
    \int_{\omega_1}^{\omega_0} \norm{\beta} d\omega
    = O(\varepsilon^{-1/2} h^{-1/2} \abs{\ln^{-2} h}) (\sqrt{b_{s_r}} + \sqrt{b_{s_{r+1}}})
    + O(h^{-1} |\ln^{\gamma_1 - 3} \varepsilon|),
  \]
  \[
    \int_{\omega_1}^{\omega_0} \abs{b(\tau)} d\tau
    = O(\varepsilon^{1/2} h^{1/2} \abs{\ln h})(\sqrt{b_{s_r}} + \sqrt{b_{s_{r+1}}})
    + O(\varepsilon |\ln^{\gamma_1} \varepsilon|).
  \]
  Finally, in Lemma~\ref{l:ode-est} we take $\delta = \norm{w(\lambda_0) - \overline w(\lambda_0)}$. Then this lemma gives the estimate~\eqref{e:v-tilde-approx}.
\end{proof}

\begin{proof}[Proof of Lemma~\ref{l:nonres}]
  Let us start with fixing the values of the constants $C_{\tilde Z}, K_Z$.
  Take $C_{\tilde Z}$ as needed by Lemma~\ref{l:beta}. Then pick $K_Z$ so that Lemma~\ref{l:cover} holds for all $C_Z > K_Z$.

  Consider the coordinate change $U$ given by~\eqref{e:varchange}.
  By Lemma~\ref{l:cover} we have
  \[
    Z_{r, r+1} \times [0, 2\pi]^2 \subset U(\tilde Z_{r, r+1} \times [0, 2\pi]^2).
  \]
  Define $\tilde w(\lambda)$ by the formula
  \[
    (\tilde w(\lambda), \varphi(\lambda), \lambda(\lambda)) = U^{-1}(w(\lambda), \varphi(\lambda), \lambda(\lambda))\in \tilde Z_{r, r+1} \times [0, 2\pi]^2.
  \]
  Set
  \[
    h(\lambda) = h(w(\lambda)), \qquad \tilde h(\lambda) = h(\tilde w(\lambda)), \qquad \overline h(\lambda) = h(\overline w(\lambda)).
  \]
  By Lemma~\ref{l:uh} we have
  \begin{equation} \label{e:h-tilde-h-2-3}
    h(\lambda) \in [0.5 \tilde h(\lambda), 1.5 \tilde h(\lambda)], \qquad \tilde h(\lambda) \in [(2/3)h(\lambda), 2h(\lambda)].
  \end{equation}
  Let us apply Lemma~\ref{l:nonres-tilde}, this lemma gives some moment $\lambda_1$ that we denote by $\tilde \lambda_1$ to avoid the conflict with $\lambda_1$ from the current lemma.
  From Lemma~\ref{l:cover} and the continuity of $w(\lambda)$ if $\tilde w(\tilde \lambda_1)$ is in $\partial \tilde \Pi$ or on the border between $\tilde Z_r$ and $\tilde Z_{r+1}$, for some $\lambda_1 < \tilde \lambda_1$ the point $w(\lambda_1)$ is in $\partial \Pi$ or on the border between $Z_r$ and $Z_{r+1}$.
  We have proved that $\lambda_1$ exists.

   The estimate on $\lambda_1 - \lambda_0$ follows from the estimate on $\tilde \lambda_1 - \lambda_0$ provided by~\eqref{e:time-nonres-zone-tilde} and from~\eqref{e:h-tilde-h-2-3}.

   The estimate $h(\lambda) \le \frac{5}{3} h(\lambda_0)$ follows from Lemma~\ref{l:uh}, as $\tilde h(\lambda)$ decreases: we have
   \[
    h(\lambda) \le \frac 5 4 \tilde h(\lambda) \le \frac 5 4 \tilde h(\lambda_0) \le
    \frac 5 4 \times \frac 4 3 h(\lambda_0).
   \]

  By~\eqref{e:h-overline-h} and~\eqref{e:h-tilde-h-2-3} we have
  \[
    \overline h(\lambda_{01}) \ge 0.5 h(\lambda_{01}) \ge 0.25 \tilde h(\lambda_{01}),
    \qquad
    \overline h(\lambda_0) \le 2 h(\lambda_0) \le 3 \tilde h(\lambda_0).
  \]
  Hence, we may apply the last part of Lemma~\ref{l:nonres-tilde} and obtain an estimate for $\norm{\tilde w(\lambda) - \overline w(\lambda)}$ for all $\lambda \in [\lambda_0, \lambda_{01}]$.
  By the estimate on $\varepsilon u$ in Lemma~\ref{l:invertible} we have
  \begin{equation}
    \abs{\norm{\tilde w(\lambda) - \overline w(\lambda)} - \norm{w(\lambda) - \overline w(\lambda)}}
    \lesssim  (\sqrt{b_{s_{r}}} + \sqrt{b_{s_{r+1}}})\sqrt{\varepsilon \tilde h(\lambda) \ln^2 \tilde h(\lambda)} + \varepsilon \ln^2 \varepsilon.
  \end{equation}
  Hence, the estimate from Lemma~\ref{l:nonres-tilde} means that for all $\lambda \in [\lambda_0, \lambda_{01}]$
  \begin{equation}
    \norm{w(\lambda) - \overline w(\lambda)} <
    e^{
      C_6 (\xi_r - \xi_{r+1})
    } \norm{w(\lambda_0) - \overline w(\lambda_0)}
    + O(\sqrt{b_{s_r}} + \sqrt{b_{s_{r+1}}}) \sqrt{\varepsilon h} \abs{\ln h}
    + O(\varepsilon |\ln^{\gamma_1} \varepsilon|).
  \end{equation}
  We have proved the estimate~\eqref{e:nonres-est}.
\end{proof}

\newpage
\section{Auxiliary system describing resonance crossing} \label{s:aux}
\subsection{Transition to auxiliary system: statement of lemma}
In this subsection we state lemmas on the auxilliary system describing passage through resonances. These lemmas will be proved in the rest of the current section.
Denote
\begin{equation} \label{e:alphabet}
  \alpha(z) = \sqrt{\varepsilon / \hatw{\pdv{\omega}{I}}}
  \sim \sqrt{\varepsilon \hat h \ln^3 \hat h},
  \qquad
  \beta(z) = \sqrt{\varepsilon \hatw{\pdv{\omega}{I}}}
  \sim \sqrt{\varepsilon \hat h^{-1} \ln^{-3} \hat h}.
\end{equation}
We have $\varepsilon = \alpha \beta$.
For given resonance $s = (s_1, s_2)$ denote $\hat \omega = s_2/s_1$ and let $\hat h(z)$ be given by $\omega(\hat h, z) = \hat \omega$.

\begin{lemma} \label{l:d-z}
  There exists $d_Z > 0$ such that for any $\hat \omega = s_2/s_1$ and $z_0$ such that
  $(\hat h(z_0), z_0) \in \mathcal B$ for any $z \in \mathbb C^n$ with $\norm{z - z_0} \le \hat \omega d_Z$ we have
  \begin{equation} \label{e:hat-h-bound}
    |\hat h(z) - \hat h(z_0)| < \hat h(z_0)/10, \qquad
    |\alpha(z) - \alpha(z_0)| < \alpha(z_0)/10, \qquad
    |\beta(z) - \beta(z_0)| < \beta(z_0)/10.
  \end{equation}
\end{lemma}

\begin{lemma} \label{l:aux}
  For any $C_{s_1} > 0$ for any large enough $C_\gamma > 0$
  there exist $C, C_\rho, c_Z > 0$ such that for any small enough $\varepsilon > 0$, any $z_0$, and any resonance $s = (s_1, s_2)$ with
  \[
    |s_1| < C_{s_1} \ln^2 \varepsilon, \qquad
    (\hat h(z_0), z_0) \in \mathcal B, \qquad
    \hat h(z_0) > C_\rho \varepsilon |\ln^5 \varepsilon|
  \]
  after a coordinate change $(p, q) \to (P, Q)$ depending on $\lambda$ (with period $2 \pi s_1$) and $z$,
  and the time and coordinate change given by
  \begin{equation} \label{e:tau-Z}
      \dv{\tau}{\lambda} = \hat \omega^{-0.5} \beta(z), \qquad
      Z = \hat \omega^{-1}(z - z_0)
  \end{equation}
  the perturbed system~\eqref{e:init} can be rewritten as
  \begin{align} \label{e:res-passing}
  \begin{split}
    P' &= - F_s(Q, z_0 + \hat \omega Z) - \hat \omega^{-2} s_1 \beta \pdv{H_7}{Q}\/(P, Q, Z) + \alpha u_P(P, Q, Z, \tau, \varepsilon), \\
    Q' &= P + \hat \omega^{-2} s_1 \beta \pdv{H_7}{P}\/(P, Q, Z) + \alpha \hat \omega^{-0.5} u_Q(P, Q, Z, \tau, \varepsilon), \\
    Z' &= \alpha \hat \omega^{-0.5} u_z(P, Q, Z, \tau, \varepsilon), \\
    \tau' &= 1
  \end{split}
  \end{align}
  in the domain
  \begin{equation} \label{e:D-res-passing}
    \mathcal D = \Big\{Z, P, Q, \tau \in \mathbb R^{n+3} :
    \norm{Z} < c_Z, \;
    \abs{P} < \frac{\hat \omega^{-0.5}}{4}, \;
    \abs{Q} < \frac{\hat \omega^{-1} C_\gamma}{4}
    \Big\}.
  \end{equation}
  The values of $\alpha$ and $\beta$ in~\eqref{e:res-passing} are taken at
  $z = z_0 + \hat \omega Z$.
  The function $F_s(Q, z)$ is described by Lemma~\ref{l:lim} below.
  We have the estimates
  \[
    \norm{H_7}_{C^2}, |u_Q|, |u_z| < C, \qquad |u_P| < C \hat \omega^{-1}
  \]
  and
  \begin{equation} \label{e:p-q-J-gamma}
    \abs{\alpha^{-1}(I - \hat I) - \sqrt{\hat \omega} P(J, \gamma, z, t)} < C |s_1 \hat \omega^{-1} \beta|,
    \qquad
    \abs{\gamma - \hat \omega Q(J, \gamma, z, t)} < C |s_1 \beta|,
  \end{equation}
  where $\gamma = \varphi - (s_2/s_1) \lambda$.
\end{lemma}

\begin{lemma} \label{l:weak-res-aux}
  Denote $D(P, Q, z, \lambda)=\alpha^{-1}\hat \omega^{-0.5}(I - \hat I)$, then this function is $2\pi s_1$-periodic in $\lambda$ and we have
  \begin{equation} \label{e:d-D-minus-P}
    |D-P|, \; \abs{\pdv{D}{Q}} < C |s_1 \hat \omega^{-1.5} \beta|, \qquad
    \abs{\pdv{D}{P} - 1} < C |s_1 \hat \omega^{-1} \beta|, \qquad
    \norm{\pdv{D}{z}} < C |s_1 \hat \omega^{-2.5} \beta|.
  \end{equation}
\end{lemma}

\begin{lemma} \label{l:lim} \;
  \begin{itemize}
    \item $F_s(Q, z)$ is $2\pi$-periodic in $Q$.
    \item $F_s$ can be continued to
    \begin{equation}
      \mathcal D_F = \big\{
        Q, z \in [0, 2\pi] \times \mathbb C^n: |z-z_0| < 0.5c
      \big\},
    \end{equation}
    where $c$ is the constant from Lemma~\ref{l:cont}, it does not depend on $s$.
    The function $F_s$ is uniquely determined by $s$ and the perturbed system in the action-angle variables~\eqref{e:init}.
    \item For any $\delta_1 > 0$ there exist $S$, $\delta_2 > 0$ such that for any $s_2, s_1$ with  $s_2 > S$ and $s_2/s_1 < \delta_2$ we have
    \[
      \norm{F_s - \Theta_3(z)}_{C^1} < \delta_1, \qquad
      \norm{\pdv{F_s}{z} - \pdv{\Theta_3}{z}}_{C^1} < \delta_1
      \qquad
      \text{in } \mathcal D_F
    \]
    (recall that $\Theta_3$ is given by~\eqref{e:Theta-3}).
    \item For each $s_2$ there exist finite set $\mathcal F_{s_2}$ such that we have
    \[
      \min_{\tilde F \in \mathcal F_{s_2}} \Big(
        \norm{F_s - \tilde F}_{C^1} + \norm{\pdv{F_s}{z} - \pdv{\tilde F}{z}}_{C^1}
      \Big)
      \; \xrightarrow[s_1 \to \infty] \; 0
      \qquad \text{in } \mathcal D_F.
    \]
  \end{itemize}
\end{lemma}

\subsection{Rescaling the action}
In this subsection we start the proof of Lemma~\ref{l:aux}, this proof continues till Subsection~\ref{ss:end-proof}. Lemma~\ref{l:d-z} is obtained as a byproduct in this subsection.

First, we use Lemma~\ref{l:cont} to continue the perturbed system in the complex domain $\mathcal D_0$ given by~\eqref{e:D-0}, we will use the notation $c_{cont}$ from~\eqref{e:D-0}. As $\omega$ is bounded in $\mathcal B$, for small enough $d_Z$ we have $\hat \omega d_Z < c_{cont}$ for all $s$.
By~\eqref{e:est}
\begin{equation} %
  \pdv{\hat h}{z} = O(h \ln h)
  \qquad \text{in } \mathcal D_0.
\end{equation}
For any $K>10$, reducing $d_Z$ if needed, we get
\begin{equation}
  |\hat h(z) - \hat h(z_0)| < \hat h(z_0)/K.
  \qquad \text{if\;\;} |z - z_0| < \hat \omega d_Z
\end{equation}
This clearly implies the first estimate in~\eqref{e:hat-h-bound}, other two estimates also follow by~\eqref{e:alphabet}.
This proves Lemma~\ref{l:d-z}.

Let us continue the proof of Lemma~\ref{l:aux}, we will assume below that $c_Z$ in this lemma satisfies $8c_w \le d_Z$.
We will use the notation $\hat \psi(z) = \psi(z, \hat h(z))$.
Let us replace the energy variable $h$ by the rescaled action $J$ that will be defined shortly.
First, note that the action $I(h, z)$ can be continued to $\mathcal D_0$ by the formula $I = \int_{\varphi=0}^{2\pi} p \pdv{q}{\varphi} d\varphi$, where $p(h, z, \varphi)$ and $q(h, z, \varphi)$ are defined in $\mathcal D_0$ by Lemma~\ref{l:cont}.
As by~\eqref{e:est} we have $\pdv{\omega}{I} \ne 0$ in $\mathcal D_0$, $\sqrt{\hatw{\pdv{\omega}{I}}}$
is uniquely continued from the real square root).
Hence, $\alpha$ and $\beta$ are correctly defined by~\eqref{e:alphabet} even for complex $z$.
Let us define (in $\mathcal D_0$) the rescaled action $J$ by the formula
\[
  J = \alpha^{-1} (I - \hat I),
\]
then
\[
  I = \hat I + \alpha J, \qquad
  \pdv{}{{I}} = \alpha^{-1} \pdv{}{J}.
\]

Denote by $\pdv{}{z}_J$ the $z$-derivative for fixed $J, \varphi, \lambda$.
Denote by $f_J = \pdv{J}{h} f_h + \pdv{J}{z} f_z$ the $J$-component of the vector field $f$, here $\pdv{J}{h} = \alpha^{-1}\omega^{-1}$.
Denote
\begin{equation}
  \Div_J f = \pdv{f_z}{z}_J + \pdv{f_\varphi}{\varphi} + \pdv{f_J}{J}.
\end{equation}
\begin{lemma} \label{l:est-div}
  $\Div_J f = O(1) + a(z) f_z$ with $a(z) = O(\hat \omega^{-1})$ in $\mathcal D_0$.
  Thus, $\Div_J f = O(\hat \omega^{-1})$.
\end{lemma}
\begin{proof}
  We need the following formula \href{http://www.owlnet.rice.edu/~fjones/chap15.pdf}{[see page~15 of http://owlnet.rice.edu/~fjones/chap15.pdf]} for the divergence in curvilinear coordinates.
  Let $\tilde x_i$ be curvilinear coordinates and $x_i$ be cartesian coordinates, let $\tilde F_i$ and $F_i$ be components of a vector field $F$ in these coordinates. Let $D$ be the Jacobian of the map $T$ given by $x = T(\tilde x)$. Then
  \begin{equation}
    \sum \pdv{F_i}{x_i} = \sum \pdv{\tilde F_i}{\tilde x_i}
    + D^{-1} \sum \tilde F_i \pdv{D}{\tilde x_i}.
  \end{equation}

  Let us apply this formula to coordinate systems $\tilde x = (J, \varphi, z)$ and $x = (p, q, z)$ (for fixed $\lambda$). The map $(I, \varphi) \mapsto (p, q)$ is volume-preserving for any fixed value of $z$, so for fixed $z$ the map $(J, \varphi) \mapsto (p, q)$ has Jacobian equal to $\alpha(z)$. Thus, we have $D = \alpha(z)$.
  \begin{equation}
    \pdv{f_p}{p} + \pdv{f_q}{q} + \pdv{f_z}{z}_{p, q} = \Div_J f
    + \alpha(z)^{-1} f_z \pdv{\alpha}{z}.
  \end{equation}
  Using that $\pdv{\hat h}{z} = O(h \ln h)$ by~\eqref{e:est}, we get $\pdv{}{z} \hat{\pdv{\omega}{h}} = \pdv[2]{\omega}{h}{z} + \pdv[2]{\omega}{h} \pdv{\hat h}{z} = O(h^{-1} \ln^{-1} h)$ by~\eqref{e:est}.
  As we can write $\alpha = \varepsilon^{1/2} \hat \omega^{-1/2} \hat{\pdv{\omega}{h}}^{-1/2}$, where only the last multiplier depends on $z$, we have
  \begin{equation} \label{e:d-alpha-z}
    \pdv{\alpha}{z} \sim
    \varepsilon^{1/2} \hat \omega^{-1/2} \hat{\pdv{\omega}{h}}^{-3/2} \pdv{}{z} \hat{\pdv{\omega}{h}}
    = \alpha O(\ln h).
  \end{equation}
  Thus $\Div_J f = O(1) + a(z) f_z$ with $a(z) = \alpha^{-1}\pdv{\alpha}{z}= O(\ln h) = O(\hat \omega^{-1})$.
\end{proof}

Let us prove some estimates that will be used later.
We have
$\pdv{\hat I}{z} = \hat{\pdv{I}{z}} + \hat{\pdv{I}{h}} \pdv{\hat h}{z} = \hat{\pdv{I}{z}} + O(h \ln^2 h)$.
As $\pdv[2]{I}{z}{h} \sim \pdv{T}{z} = O(\ln h)$,
we have
$\pdv{}{z}\/ (I - \hat I) = O(h-\hat h) \ln h + O(h \ln^2 h) = O(h \ln^2 h)$ in $\mathcal D_0$.
Thus we have in $\mathcal D_0$ (as $I - \hat I = O(h \ln h)$)
\begin{align}
\begin{split} \label{e:d-J-z}
  \pdv{J}{z}
  &= \alpha^{-1} \pdv{}{z}\/(I - \hat I) - \alpha^{-2} \pdv{\alpha}{z}\/(I - \hat I)
  = O(h \ln^2 h) \alpha^{-1} + O(\ln h) \alpha^{-1}(I - \hat I) = \\
  &= O(h \ln^2 h) \alpha^{-1} = O(\ln^{-1} h) \beta^{-1}.
\end{split}
\end{align}

Denote
\begin{align}
\begin{split}
  \mathcal D_{1, 0} &= \Big\{z, J, \varphi, \lambda \in \mathbb C^{n+3} :
  \abs{z - z_0} < 4 c_w \hat \omega, \;
  \abs{J} < c_{cont, J} \beta^{-1} \ln^{-2} \hat h, \;
  \abs{\Im \varphi} < c_{cont} \hat \omega, \;
  \abs{\Im \lambda} < c_{cont}
  \Big\}, \\
  \mathcal D_1 &= \Big\{z, J, \varphi, \lambda \in \mathbb C^{n+3} :
  \abs{z - z_0} < 2 c_w \hat \omega, \;
  \abs{J} < \frac{c_{cont, J}}{2} \beta^{-1} \ln^{-2} \hat h, \;
  \abs{\Im \varphi} < \frac{c_{cont}}{2} \hat \omega, \;
  \abs{\Im \lambda} < \frac{c_{cont}}{2}
  \Big\}.
\end{split}
\end{align}
Here the constant $c_{cont, J} > 0$ is chosen so that
the image of $\mathcal D_{1, 0}$ under the map $J, z, \varphi, \lambda \mapsto h, z, \varphi, \lambda$ lies in $\mathcal D_0$.
We will use the domain $\mathcal D_1$ below, the domain $\mathcal D_{1, 0}$ is needed to obtain Cauchy estimates in $\mathcal D_1$.
The width of these domains in $J$ is $O(\alpha^{-1} \hat h \ln \hat h) = O(\sqrt{\varepsilon^{-1} \hat h\ln^{-1} \hat h}) = O(\beta^{-1} \ln^{-2} \hat h)$.

We have $\pdv{f_z}{z}_J = O(\ln h)$ in $\mathcal D_1$ by Cauchy formula.
However, we have weaker estimates $O_*(h^{-1})$ for $\pdv{f_\varphi}{\varphi}$ and $\pdv{f_J}{J}$ (this is not used later, so we skip the proof). Let us set
\begin{equation} \label{e:def-H1}
  H_1(J, z, \varphi, \lambda) = - \int_{\psi = 0}^\varphi (
    f_J(J, z, \psi, \lambda, \varepsilon{=}0) - \langle f_J(J, z, \varphi, \lambda, \varepsilon{=}0) \rangle_\varphi
  ) d\psi.
\end{equation}
We have $f_J = \pdv{J}{h} f_h|_{\varepsilon=0} + \pdv{J}{z} f_z|_{\varepsilon=0} = O(\ln h) \alpha^{-1} + O(\ln^{-1} h) \beta^{-1} = O(\ln h) \alpha^{-1}$. By the estimates above
we have
$\int_{\psi = 0}^\varphi \pdv{J}{z}f_z d\psi = O(\beta^{-1}\ln^{-1} h) = O(\alpha^{-1})$.
As $\pdv{J}{h} = \omega^{-1}\alpha^{-1}$, we have
$\int_{0}^\varphi \pdv{J}{h}f_h d\psi = \alpha^{-1} \int_0^{\omega^{-1} \varphi} f_h dt = O(\alpha^{-1})$
by Lemma~\ref{l:int-psi}.
Thus we have $\int_0^{\varphi} f_J d \psi = O(\alpha^{-1})$ for any $\varphi$.
This implies
\begin{equation} \label{e:est-avg-f-J}
  \langle f_J \rangle_\varphi = O(\alpha^{-1}), \qquad
  H_1 = O(\alpha^{-1}) \qquad
  \text{in }\mathcal D_1.
\end{equation}
We can compute
\begin{align}
\begin{split}
  \pdv{H_1}{J}
  &= - \int_{\psi = 0}^\varphi \Big(
    \pdv{f_J}{J} - \pdv{}{J} \langle f_J \rangle_\varphi
  \Big) d\psi =\\
  &=
  f_\varphi - f_\varphi|_{\varphi=0} +
  \int_{\psi = 0}^\varphi \Big(
    \pdv{f_z}{z}_J - \Div_J f
    + \pdv{}{J} \langle f_J \rangle_\varphi
  \Big) d\psi.
\end{split}
\end{align}
Hence, we have
\begin{equation}
  f_J = - \pdv{H_1}{\varphi} + \tilde f_J(J, z, \lambda) + \check f_J(J, z, \varphi, \lambda, \varepsilon), \qquad
  f_\varphi = \pdv{H_1}{J} + \tilde f_\varphi(J, z, \varphi, \lambda, \varepsilon),
\end{equation}
where
\begin{align} \label{e:g}
\begin{split}
  \tilde f_J &= \langle f_J|_{\varepsilon=0} \rangle_\varphi, \\
  \tilde f_\varphi &= \bigg( f_\varphi|_{\varphi=0}
  - \int_{\psi = 0}^\varphi \Big(
    \pdv{f_z}{z}_J - \Div_J f
    + \pdv{}{J} \langle f_J \rangle_\varphi
  \Big) d\psi \bigg)\bigg|_{\varepsilon=0}
  + (f_\varphi - f_\varphi|_{\varepsilon=0}), \\
  \check f_J &= f_J - f_J|_{\varepsilon=0}.
\end{split}
\end{align}
By~\eqref{e:est-avg-f-J} we have $\tilde f_J = O(\alpha^{-1}).$
We also have
\begin{equation} \label{e:d-avg-fJ-J}
  \pdv{}{J} \langle f_J \rangle_\varphi = \langle \Div_J f - \pdv{f_z}{z}_J \rangle_\varphi = O(\ln h).
\end{equation}
Note that $\tilde f_\varphi$ is clearly $2\pi$-periodic. Let us also set
$\tilde f_z = f_z$.
\begin{lemma} \label{l:diff-eps-0}
  \begin{equation}
    f_J|_{\varepsilon=0} - f_J = O(\ln h) \beta,
    \qquad
    f_\varphi|_{\varepsilon=0} - f_\varphi = O(\varepsilon h^{-1} \ln^{-2} h).
  \end{equation}
\end{lemma}
\begin{proof}
  Denote $\vec{f} = (f_p, f_q, f_z)$. The maps $\vec{f} \mapsto f_J$ and $\vec{f} \mapsto f_\varphi$ are linear. As $\vec{f}|_{\varepsilon=0} - \vec{f} = O(\varepsilon)$, the lemma follows from the estimates $f_J = O(\ln h) \alpha$ and $f_\varphi = O(h^{-1} \ln^{-2} h)$ that hold for any $\vec{f}$.
\end{proof}
As $\varphi = 0$ is given by a transversal to one of the separatrices that is separated from the saddle, the time $t=\omega^{-1} \varphi$ is a smooth function of the coordinates $p, q, z$ near $\varphi = 0$ (we treat values of $\varphi$ near $2\pi$ as negative values near $0$). Therefore, from $\varphi = \omega t$ we get
\[
  \pdv{\varphi}{q}\Big|_{\varphi = 0} = \omega \pdv{t}{q}\Big|_{t = 0} = O(\ln^{-1} h),
  \;\;
  \pdv{\varphi}{p}\Big|_{\varphi = 0} = \omega \pdv{t}{p}\Big|_{t = 0} = O(\ln^{-1} h),
  \;\;
  \pdv{\varphi}{z}\Big|_{\varphi = 0} = \omega \pdv{t}{z}\Big|_{t = 0} = O(\ln^{-1} h).
\]
Hence,
$f_\varphi|_{\varphi = 0} = O(\ln^{-1} h)$.
We have by Lemma~\ref{l:est-div} (we use the notation $a(z)$ from this lemma, with $a(z) = O(\ln h)$)
\[
  \int_{\psi=0}^\varphi \Div_J f - \langle \Div_J f \rangle_\varphi d\psi
  = O(1) +
  a(z)\int_{\psi=0}^\varphi (f_z - f_z(C)) -  \langle f_z - f_z(C) \rangle_\varphi d\psi.
\]
This is $O(1)$, as the integral above can be estimated by Lemma~\ref{l:int-psi}.
By~\eqref{e:g} this implies $\tilde f_\varphi = O(1)$ (note that $f_\varphi|_{\varepsilon=0} - f_\varphi = O(1)$ by Lemma~\ref{l:diff-eps-0}, as $h \gtrsim \varepsilon \ln^5 \varepsilon$).

The system~\eqref{e:init} rewrites in $\mathcal D_1$ (we use the notation $\omega(J, z) = \omega(h(J, z), z)$) as
\begin{align}
\begin{split} \label{e:init-J}
  \dot J &= -\varepsilon \pdv{H_1}{\varphi} + \varepsilon \tilde f_J(J, z, \lambda) + \varepsilon \check f_J(J, z, \varphi, \lambda, \varepsilon), \\
  \dot z &= \varepsilon \tilde f_z(J, z, \varphi, \lambda, \varepsilon), \\
  \dot \varphi &= \omega(J, z) + \varepsilon \pdv{H_1}{J} + \varepsilon \tilde f_\varphi(J, z, \varphi, \lambda, \varepsilon), \\
  \dot \lambda &= 1.
\end{split}
\end{align}

\subsection{Transition to resonant phase} \label{s:res-phase}
It will be convenient to use new angle variables $\gamma, \mu$ given by
\begin{equation} \label{e:mu-gamma}
  \gamma = \varphi - (s_2/s_1) \lambda, \; \mu = \lambda/s_1;
  \qquad
  \varphi = \gamma + s_2 \mu, \; \lambda = s_1 \mu.
\end{equation}
Note that $\pdv{}{\varphi} = \pdv{}{\gamma}$, as both these derivatives are taken for fixed $\lambda$ and $\mu$. After the coordinate change
$(J, z, \varphi, \lambda) \mapsto (J, z, \gamma, \mu)$
and the time change $\psi' = \dv{\psi}{\mu} = s_1 \dot \psi$
the system~\eqref{e:init-J} rewrites as
\begin{align}
\begin{split}
  J' &= - s_1 \varepsilon \pdv{H_1}{\gamma} + s_1 \varepsilon \tilde f_J + s_1 \varepsilon \check f_J, \\
  z' &= s_1 \varepsilon \tilde f_z, \\
  \gamma' &= s_1 \omega - s_2 + s_1 \varepsilon \pdv{H_1}{J} + s_1 \varepsilon \tilde f_\varphi, \\
  \mu' &= 1.
\end{split}
\end{align}
The coefficients of this system are $2\pi$-periodic in $\gamma$ and $\mu$ and also are invariant under the translation
$(\gamma, \mu) \mapsto (\gamma - 2\pi \frac{s_2}{s_1}, \mu + \frac{2\pi}{s_1})$.
Set
\begin{align}
\begin{split}
  H_2(J, z, \gamma, \mu) = \alpha \Big(
    H_1(J, z, \gamma+s_2\mu, s_1\mu) \; - \; \gamma \tilde f_J(0, z, s_1\mu)
  \Big)
  + \beta^{-1} \int_0^J (\omega(\tilde J, z) - s_2/s_1) d \tilde J.
\end{split}
\end{align}
Note that $H_2$ is $2\pi$-periodic in $\mu$.
Set
\begin{equation}
  g_J(J, z, \gamma, \mu, \varepsilon) = \tilde f_J(J, z, s_1 \mu) - \tilde f_J(0, z, s_1 \mu) + \check f_J,
  \qquad
  g_z(J, z, \gamma, \mu, \varepsilon) = \tilde f_z,
  \qquad g_\varphi(J, z, \gamma, \mu, \varepsilon) = \tilde f_\varphi.
\end{equation}
Then our system rewrites as
\begin{align}
\begin{split} \label{e:before-averaging}
  J' &= - s_1\beta \pdv{H_2}{\gamma}
  + s_1 \varepsilon g_J, \\
  z' &= s_1 \varepsilon g_z, \\
  \gamma' &= s_1\beta \pdv{H_2}{J}
    + s_1 \varepsilon g_\varphi, \\
  \mu' &= 1.
\end{split}
\end{align}

We will consider the Hamiltonian part of this system (i.e. without the $g_{*}$ terms) in the domain
\begin{equation} \label{e:D2}
  \mathcal D_2 = \Big\{z, J, \gamma, \mu \in \mathbb C^{n+3} :
  \abs{z-z_0} < 2 c_w \hat \omega, \;
  \abs{J} < 1, \;
  \abs{\Im \gamma} < \frac{c_{cont}}{4} \hat \omega, \;
  \abs{\Re \gamma} < C_\gamma, \;
  \abs{\Im \mu} < \frac{c_{cont}}{4 s_1}
  \Big\},
\end{equation}
where $C_\gamma > 0$ should satisfy $C_\gamma \gg \hat \omega c_{cont}$. Note that we have $\beta^{-1} \ln^{-3} h > 1$ for small enough $\varepsilon$ given $h > \varepsilon |\ln \varepsilon|^\rho$ for $\rho > 1$. Thus, the image of this domain under the map $J, z, \gamma, \mu \mapsto J, z, \varphi, \lambda$ lies inside $\mathcal D_1$.
Using the estimates~\eqref{e:est} on $\pdv{\omega}{I}$, $\pdv[2]{\omega}{I}$, we can compute
$\beta^{-1} (\omega(I{=}\hat I + \alpha J) - s_2/s_1)
= J + O(\beta \ln^2 h)$.
As $H_1, \tilde f_J = O(\alpha^{-1})$ by~\eqref{e:est-avg-f-J}, we have
$H_2 = O(1)$ in $\mathcal D_2$ (also for $\rho > 1$).
We also have $g_J = O(\ln h)$; $g_\varphi, g_z = O(1)$ in the real part of $\mathcal D_2$. Indeed, the estimates on $\tilde f_\varphi$ and $\tilde f_z$ were obtained above and the estimate on $g_J$ follows from~\eqref{e:d-avg-fJ-J} and Lemma~\ref{l:diff-eps-0}.

As we have $\alpha H_1 = O(1)$ in $\mathcal D_1$ by~\eqref{e:est-avg-f-J}, we get Cauchy estimate $\pdv{\alpha H_1}{J} = O(\beta  \ln^2 h)$ in $\mathcal D_2$.
This allows us to separate the main part of $H_2$: in $\mathcal D_2$ we have
\begin{align}
\begin{split} \label{e:H-pendulum}
  &H_2 = J^2/2 + H_{2, 0}(z, \gamma, \mu) + \beta \ln^2 h \; H_{2, 1}(J, z, \gamma, \mu),
  \qquad \text{where} \\
  &H_{2, 0} = \alpha \hat H_1(z, \gamma + s_2\mu, s_1\mu) - \gamma \alpha \hat {\tilde f}_J(z, s_1\mu),
  \qquad
  H_{2, 1} = O(1).
\end{split}
\end{align}

\subsection{Averaging over time}
Let us state the following lemma that follows from~\cite{neishtadt1984separation}.
\begin{lemma} \label{l:averaging}
  Consider a Hamiltonian system with the Hamiltonian $\varepsilon H(p, q, t)$ periodically depending on time $t$ (with the period $2\pi$) and slow variables $p, q$:
  \begin{align} \label{e:lem-avg}
  \begin{split}
    \dot q  &= \varepsilon \pdv{H}{p}\/(p, q, t), \\
    \dot p  &= -\varepsilon \pdv{H}{q}\/(p, q, t), \\
    \dot t  &= 1.
  \end{split}
  \end{align}
  Assume that the Hamiltonian $\varepsilon H(p, q, t)$ is defined for $(p, q)$ in a complex neighborhood $U_\delta$ of some real domain $U \subset \mathbb R^2$ of width $\delta > 0$:
  $U_\delta = U + \{ p, q \in \mathbb C^2; \abs{p}, \abs{q} \le \delta \}$
  and $t \in S^1_t = \mathbb R / 2 \pi \mathbb Z$
  and we have $\abs{H} < C_H$ in $U_\delta \times S^1_t$. We assume $H$ to be analytic in $p, q$ and continuous in $t$.

  Then there is $C > 0$ depending only on $C_H$ and $\delta$ such that for all $\varepsilon \in (0, C)$
  there are new canonical variables $\tilde p(p, q, t), \; \tilde q(p, q, t)$ with
  \[
    \abs{\tilde p(p, q, t) - p} \le C^{-1} \varepsilon,
    \qquad
    \abs{\tilde q(p, q, t) - q} \le C^{-1} \varepsilon
  \]
  such that in these coordinates our system is defined in $V = U_{0.5 \delta} \times S^1_t$ and is given by the Hamiltonian
  \[
    \varepsilon \overline H(\tilde p, \tilde q) + \varepsilon \Delta H(\tilde p, \tilde q, t)
  \]
  with
  \[
    \norm{\overline H(\tilde p, \tilde q)
    - \langle H(p, q, t) \rangle_t|_{p=\tilde p, q = \tilde q}}_{V}
    \le C^{-1} \varepsilon,
    \qquad
    \norm{\Delta H}_{V} \le \exp(-C \varepsilon^{-1}).
  \]

\end{lemma}
\noindent Let us obtain explicit dependence of the estimates in this lemma on the width of the complex domain where the system is defined.
\begin{corollary} \label{c:averaging}
  Consider the system~\eqref{e:lem-avg} with the Hamiltonian $\varepsilon H$ periodically depending on time $t$ (with the period $2\pi$) and slow variables $p, q$.
  Assume that the Hamiltonian $\varepsilon H(p, q, t)$ is defined for $(p, q)$ in a complex neighborhood $U_{\delta_p, \delta_q}$ of some real domain $U \subset \mathbb R^1$ of width $\delta_p, \delta_p \in (0, 1)$ in $p$ and $q$, respectively:
  \[
    U_{\delta_p, \delta_q} = U + \{ p, q \in \mathbb C^2; \abs{p} < \delta_p, \abs{q} < \delta_q \}
  \]
  and $t \in S^1_t = \mathbb R / 2 \pi \mathbb Z$ and we have $\abs{H} < C_H$ in $U_{\delta_p, \delta_q} \times S^1_t$.
  We assume $H$ to be analytic in $p$, $q$ and continuous in $t$,

  Then there is $C > 0$ depending only on $C_H$ such that for all $\varepsilon \in (0, C \delta_p \delta_q)$
  there are new canonical variables $\tilde p(p, q, t), \; \tilde q(p, q, t)$ with
  \[
    \abs{\tilde p(p, q, t) - p} \le C^{-1} \varepsilon \delta_q^{-1},
    \qquad
    \abs{\tilde q(p, q, t) - q} \le C^{-1} \varepsilon \delta_p^{-1}
  \]
  such that in these coordinates our system is defined in $V = U_{0.5 \delta_p, 0.5 \delta_q} \times S^1_t$ and is given by the Hamiltonian
  \[
    \varepsilon \overline H(\tilde p, \tilde q) + \varepsilon \Delta H(\tilde p, \tilde q, t)
  \]
  with
  \[
    \norm{\overline H(\tilde p, \tilde q)
    - \langle H(p, q, t) \rangle_t|_{p=\tilde p, q = \tilde q}}_{V}
    \le C^{-1} \varepsilon \delta_p^{-1} \delta_q^{-1},
    \qquad
    \norm{\Delta H}_{V} \le \exp(-C \delta_p \delta_q \varepsilon^{-1}).
  \]

\end{corollary}
\begin{proof}
  Let us make a coordinate change $p = \delta_p p'$, $q = \delta_q q'$. The Hamiltonian in the new coordinates is $\varepsilon (\delta_p \delta_q)^{-1} H$. Denote $\varepsilon' = \varepsilon (\delta_p \delta_q)^{-1}$, then the new system is given by the Hamiltonian $\varepsilon' H$. This system is analytic for $(p', q')$ in a complex neighborhood $U'_{1}$ of some real domain $U'$ with width $1$ in both $p'$ and $q'$. The domain $U'$ is large, but the constant $C$ in Lemma~\ref{l:averaging} does not depend on $U'$ and $U'_\delta$, it only depends on the width $\delta$.
  Lemma~\ref{l:averaging} gives us new coordinates $\tilde p', \tilde q'$ with
  $\abs{\tilde p' - p'}, \abs{\tilde q' - q'} \le C^{-1} \varepsilon'$.
  In the coordinates $\tilde p', \tilde q'$ the system is given by the Hamiltonian
  $\varepsilon' \overline H'(\tilde p', \tilde q') + \varepsilon' \Delta H'(\tilde p', \tilde q', t)$ with
  \[
    \norm{\overline H'(\tilde p', \tilde q')
    - \langle H(p, q, t) \rangle_t|_{p=\delta_p\tilde p', q = \delta_q \tilde q'}}_{U'_{0.5}\times S^1_t}
    \le C^{-1} \varepsilon \delta_p^{-1} \delta_q^{-1};
    \;
    \norm{\Delta H'}_{U'_{0.5}\times S^1_t} \le \exp(-C \delta_p \delta_q \varepsilon^{-1}).
  \]
  Set $\tilde p = \delta_p \tilde p'$, $\tilde q = \delta_q \tilde q'$ and let $\overline H(\tilde p, \tilde q)$, $\Delta H(\tilde p, \tilde q)$ be $\overline H'$ and $\Delta H'$ written in these coordinates.
  It is easy to check that the estimates above imply the estimates in the statement of this corollary.
\end{proof}

\subsection{After averaging} \label{ss:end-proof}
Let us recall the system~\eqref{e:before-averaging}
\begin{align}
\begin{split}
  J' &= - s_1\beta \pdv{H_2}{\gamma}
  + s_1 \varepsilon g_J, \\
  z' &= s_1 \varepsilon g_z, \\
  \gamma' &= s_1\beta \pdv{H_2}{J}
    + s_1 \varepsilon g_\varphi, \\
  \mu' &= 1.
\end{split}
\end{align}
The Hamiltonian part of this system is defined in
\begin{equation}
  \mathcal D_2 = \Big\{z, J, \gamma, \mu \in \mathbb C^{n+3} :
  \abs{z-z_0} < 2 c_Z \hat \omega, \;
  \abs{J} < 1, \;
  \abs{\Im \gamma} < \frac{c_{cont}}{4} \hat \omega, \;
  \abs{\Re \gamma} < C_\gamma, \;
  \abs{\Im \mu} < \frac{c_{cont}}{4 s_1}
  \Big\},
\end{equation}
and the whole system is defined in the real part of this complex domain.
Let us apply Lemma~\ref{c:averaging} to the Hamiltonian part of~\eqref{e:before-averaging}, with $t = \mu$,
$U = (-0.5, 0.5)_J \times (-C_\gamma + \frac{c_{cont} \hat \omega}{4}, C_\gamma - \frac{c_{cont} \hat \omega}{4})_\gamma$,
$\delta_J = 0.5$, $\delta_\gamma = \frac{c_{cont} \hat \omega}{4}$ and $\varepsilon_1 = s_1 \beta$. Here we denote by $\varepsilon_1$ the $\varepsilon$ variable used in Corollary~\ref{c:averaging} to distinguish it from $\varepsilon$ in~\eqref{e:init}. This corollary gives new coordinates that we denote
$\tilde P, \tilde Q$.
Corollary~\ref{c:averaging} is applied separately for different values of $z$, but it is easy to check that the construction in~\cite{neishtadt1984separation} gives $\tilde P(J, \gamma, z, \mu)$ and $\tilde Q(J, \gamma, z, \mu)$ that are analytic in $z$.
Let us make a scale transformation
\begin{equation}
  Q = \tilde Q / \hat \omega, \; P = \tilde P / \sqrt{\hat \omega}; \qquad
  \tilde Q = \hat \omega Q, \; \tilde P = \sqrt{\hat \omega} P,
\end{equation}
this will simplify the main part of our system that will be written later.
By Corollary~\ref{c:averaging} we have
\begin{equation} \label{e:pq-close-Jgamma}
  \abs{J - \sqrt{\hat \omega} P(J, \gamma, z, t)} = O(s_1 \hat \omega^{-1} \beta),
  \qquad
  \abs{\gamma - \hat \omega Q(J, \gamma, z, t)} = O(s_1 \beta).
\end{equation}
In $\tilde P, \tilde Q$ coordinates the system~\eqref{e:before-averaging} without $g$ is given by the Hamiltonian
\begin{align}
\begin{split}
  &\tilde H(\tilde P, \tilde Q, z, \mu) = s_1\beta \tilde H_3(\tilde P, \tilde Q, z) + s_1^2 \hat \omega^{-1} \beta^2 \tilde H_4(\tilde P, \tilde Q, z)
    + s_1 \beta \exp(-C \hat \omega s_1^{-1} \beta^{-1}) \tilde H_5(\tilde P, \tilde Q, z, \mu), \\
  &\tilde H_3(\tilde P, \tilde Q, z)
  = \langle H_2(J, \gamma, z, \mu) \rangle_\mu|_{J = \tilde P, \gamma = \tilde Q}, \\
  &\tilde H_3, \tilde H_4, \tilde H_5 = O(1). \\
\end{split}
\end{align}
Denote by $H(P, Q, z, \mu)$ the new Hamiltonian after scaling, we have
$H = \hat \omega^{-1.5} \tilde H$.
Hence, in the $P, Q$ coordinates this rewrites as
\begin{align}
\begin{split}
  &H = \hat \omega^{-1.5} s_1\beta H_3( P,  Q, z) + s_1^2 \hat \omega^{-2.5} \beta^2  H_4( P,  Q, z)
    + \hat \omega^{-1.5} s_1 \beta \exp(-C \hat \omega s_1^{-1} \beta^{-1})
    H_5( P,  Q, z, \mu), \\
  &H_3(P, Q)
    = \langle H_2(J, \gamma, \mu) \rangle_\mu|_{J = \sqrt{\hat \omega}P, \gamma = \hat \omega Q}, \\
  & H_3,  H_4,  H_5 = O(1). \\
\end{split}
\end{align}
This system is defined in the domain
(we reduce this domain a bit to have a shorter formula, taking into account that
 $C_\gamma \gg c_{cont} \hat \omega$)
\begin{equation} \label{e:D3}
  \mathcal D_3 = \Big\{z, P, Q, \mu \in \mathbb C^{n+2} \times \mathbb R :
  \abs{z - z_0} < 2 c_Z \hat \omega, \;
  \abs{P} < \frac{\hat \omega^{-0.5}}{2}, \;
  \abs{\Re Q} < \frac{C_\gamma \hat \omega^{-1}}{2}, \;
  \abs{\Im Q} < \frac{c_{cont}}{8} %
  \Big\}.
\end{equation}
Let us also consider real domain
\begin{equation} \label{e:D4}
  \mathcal D_4 = \Big\{z, P, Q, \mu \in \mathbb R^{n+3} :
  \abs{z - z_0} < c_Z \hat \omega, \;
  \abs{P} < \frac{\hat \omega^{-0.5}}{4}, \;
  \abs{Q} < \frac{C_\gamma \hat \omega^{-1}}{4}
  \Big\}
\end{equation}
and the same domain rewritten using $Z=\hat \omega^{-1}(z-z_0)$ instead of $z$:
\begin{equation} \label{e:D4-w}
  \mathcal D = \Big\{Z, P, Q, \mu \in \mathbb R^{n+3} :
  \abs{Z} < c_Z, \;
  \abs{P} < \frac{\hat \omega^{-0.5}}{4}, \;
  \abs{Q} < \frac{C_\gamma \hat \omega^{-1}}{4}
  \Big\}.
\end{equation}
We have Cauchy estimates valid in $\mathcal D$
\begin{equation} \label{e:Cauchy-after-avg}
  \norm{H_i}_{C^2} = O(1)
  \text{\qquad for } i=3,4,5.
\end{equation}
Denote
$C_{H_5} = \hat \omega^{-1.5} s_1 \beta \exp(-C \hat \omega s_1^{-1} \beta^{-1})$.
For large enough $C_\rho$ we have
\begin{equation} \label{e:cH5}
  C s_1 \hat \omega^{-1} \beta < |\ln^{-1} \varepsilon| / 4
\end{equation}
and thus
 $C_{H_5} = O(\varepsilon^3)$.
By~\eqref{e:Cauchy-after-avg} this means that the corresponding terms in Hamiltonian equations are also $O(\varepsilon^3)$.

Let us now include the terms appearing after we reintroduce the terms $s_1 \varepsilon g$ in~\eqref{e:before-averaging} rewritten in the new coordinates.
Denote
\begin{align}
\begin{split}
  u_P &=
    \pdv{\tilde P}{J} g_J
    + \pdv{\tilde P}{\gamma} g_\varphi
    + \pdv{\tilde P}{z} g_z
    - C_{H_5} s_1^{-1} \varepsilon^{-1} \hat \omega^{0.5} \pdv{H_5}{Q}, \\
  u_Q &=
    \pdv{\tilde Q}{J} g_J
    + \pdv{\tilde Q}{\gamma} g_\varphi
    + \pdv{\tilde Q}{z} g_z
    + C_{H_5} s_1^{-1} \varepsilon^{-1} \hat \omega \pdv{H_5}{P}, \\
  u_z &= g_z.
\end{split}
\end{align}
From the estimates $g_J = O(\ln h)$, $g_\varphi, g_z = O(1)$ we get $u_P = O(\ln h)$, $u_Q, u_z = O(1)$
(we use that $\pdv{(\tilde P - J)}{x} = O(1)$, $\pdv{(\tilde Q - \gamma)}{x} = O(\ln^{-1} \varepsilon)$ for $x=J, \gamma, z$ in $\mathcal D_4$ by the Cauchy formula and~\eqref{e:pq-close-Jgamma} and~\eqref{e:cH5}).
Now the system~\eqref{e:before-averaging} rewrites as
\begin{align}
\begin{split} \label{e:after-averaging}
  P' &= - \hat \omega^{-1.5} s_1 \beta \pdv{H_3}{Q}\/(P, Q, z) - \omega^{-2.5} s_1^2 \beta^2 \pdv{H_4}{Q}\/(P, Q, z) + s_1 \varepsilon \hat \omega^{-0.5} u_P(P, Q, z, \mu), \\
  Q' &= \hat \omega^{-1.5} s_1 \beta \pdv{H_3}{P}\/(P, Q, z) + \hat \omega^{-2.5} s_1^2 \beta^2 \pdv{H_4}{P}\/(P, Q, z) + s_1 \varepsilon \hat \omega^{-1} u_Q(P, Q, z, \mu), \\
  z' &= s_1 \varepsilon u_z(P, Q, z, \mu), \\
  \mu' &= 1.
\end{split}
\end{align}
After the time change
$\dv{\tau}{\mu} = s_1 \hat \omega^{-0.5} \beta(z)$,
$\dv{\tau}{\lambda} = \hat \omega^{-0.5} \beta(z)$
we obtain the system (we recycle $'$ to denote also the derivative with respect to the new time $\tau$)
\begin{align}
\begin{split}
  P' &= - \hat \omega^{-1} \pdv{H_3}{Q}\/(P, Q, z) - \hat \omega^{-2} s_1 \beta \pdv{H_4}{Q}\/(P, Q, z) + \alpha u_P(P, Q, z, \tau), \\
  Q' &= \hat \omega^{-1} \pdv{H_3}{P}\/(P, Q, z) + \hat \omega^{-2} s_1 \beta \pdv{H_4}{P}\/(P, Q, z) + \alpha \hat \omega^{-0.5} u_Q(P, Q, z, \tau), \\
  z' &= \alpha \hat \omega^{0.5} u_z(P, Q, z, \tau), \\
  \tau' &= 1.
\end{split}
\end{align}

Let us now use~\eqref{e:H-pendulum} to separate the main part of this system.
We replace
\[
  H_3 = \langle H_2(J, \gamma, z, \mu) \rangle_\mu|_{J = \sqrt{\hat \omega}P, \gamma = \hat \omega Q}
\]
with its main part (corresponding to $J^2/2 + H_{2, 0}(z, \gamma, \mu)$ from~\eqref{e:H-pendulum}) that we denote $H_6$ and add the remainder to $H_4$ (this sum is denoted $H_7$).
We have in $\mathcal D_3$
\begin{equation} \label{e:H56}
  H_6(P, Q) = \hat \omega \frac{P^2}{2} + \langle H_{2, 0} \rangle_\mu|_{\gamma= \hat \omega Q},
  \qquad
  H_7(P, Q) = H_4(P, Q)
  + O(s_2^{-1}) \langle H_{2, 1} \rangle_\mu|_{J = \sqrt{\hat \omega}P, \gamma=\hat \omega Q} = O(1).
\end{equation}
Denote
\begin{equation} \label{e:def-Fs}
  F_s = \hat \omega^{-1} \pdv{H_6}{Q}.
\end{equation}
Then the system above rewrites as
\begin{align} \label{e:sys-p-q}
\begin{split}
  P' &= - F_s(Q, z) - \hat \omega^{-2} s_1 \beta \pdv{H_7}{Q}\/(P, Q, z) + \alpha u_P(P, Q, z, \tau), \\
  Q' &= P + \hat \omega^{-2} s_1 \beta \pdv{H_7}{P}\/(P, Q, z) + \alpha \hat \omega^{-0.5} u_Q(P, Q, z, \tau), \\
  z' &= \alpha \hat \omega^{0.5} u_z(P, Q, z, \tau), \\
  \tau' &= 1.
\end{split}
\end{align}
This system is defined in the domain $\mathcal D_4$;
using $Z$ instead of $z$ gives the system~\eqref{e:res-passing} defined in $\mathcal D$.
As $H_7 = O(1)$ in $\mathcal D_3$, in $\mathcal D$ we also have $\norm{H_7}_{C_2} = O(1)$ by Cauchy formula (moreover, $\pdv{H_7}{P} = O(\sqrt{\hat \omega})$).
Thus we have in $\mathcal D$:
\[
  \norm{H_7}_{C^2}, \norm{u_Q}_{C}, \norm{u_z}_{C} = O(1),
  \qquad
  \norm{u_P}_{C} = O(\ln h).
\]
This completes the proof of Lemma~\ref{l:aux}. \qed

\begin{proof}[Proof of Lemma~\ref{l:weak-res-aux}]
  As $D$ depends $2\pi$-periodically on $\mu$, it depends $2\pi s_1$-periodically on $\lambda$.
  As $D = J/\sqrt{\hat \omega}$, by~\eqref{e:pq-close-Jgamma} we have
  \[
    \abs{D - P} = O(s_1 \hat \omega^{-1.5} \beta) \text{\qquad in } \mathcal D_3 .
  \]
  The estimates of Lemma~\ref{l:weak-res-aux} in $\mathcal D_4 \subset \mathcal D_3$ follow by Cauchy formula.
\end{proof}

\subsection{Main part of the Hamiltonian}
In this section we prove Lemma~\ref{l:lim}.
By~\eqref{e:def-Fs},~\eqref{e:H56},~\eqref{e:H-pendulum} and~\eqref{e:mu-gamma} we have
\begin{equation}
  -F_s = - \hat \omega^{-1} \pdv{H_6}{Q}
  = - \alpha \pdv{}{\gamma} \langle
    \hat H_1(z, \gamma+s_2\mu, s_1\mu)
    - \gamma \hat{\tilde f}_J(z, s_1\mu)
  \rangle_\mu|_{\gamma = \hat \omega Q}.
\end{equation}
As $\pdv{}{\gamma}$ and $\langle \cdot \rangle_\mu$ commute, $\pdv{}{\gamma}\big|_{\mu=const} = \pdv{}{\varphi}\big|_{\lambda=const}$ and $\tilde f_J = \langle f_J \rangle_\varphi$ by~\eqref{e:g}, this rewrites as
\begin{equation}
  -F_s = \alpha \Big\langle
    \Big(
      - \pdv{}{\varphi} \hat H_1(z, \varphi, \lambda)
      + \langle \hat f_J(z, \varphi, \lambda) \rangle_\varphi
    \Big) \Big|_{\varphi = s_2\mu + \hat \omega Q, \lambda = s_1\mu}
  \Big\rangle_\mu.
\end{equation}
From~\eqref{e:def-H1}
we obtain
\begin{equation}
  -F_s = \alpha \langle
    \hat f_J(z, \varphi {=} s_2\mu {+} \hat \omega Q, \lambda {=} s_1\mu, \varepsilon{=}0)
  \rangle_\mu
  = \alpha \langle \hat f_J(z, \varphi, \lambda {=} \hat \omega^{-1} \varphi {-} Q, \varepsilon{=}0) \rangle_{\varphi \in [0, 2s_2 \pi]}.
\end{equation}
Let us again use the notation $t = \hat \omega^{-1} \varphi = (s_1/s_2) \varphi$ for the time for the unperturbed system.
We have $\alpha \hat f_J =  \hat \omega^{-1} \hat f_h + \alpha \hat{\pdv{J}{z}} \hat f_z$.
Denote (we reuse the notation $g$ already used in Section~\ref{s:res-phase}, as $g$ from Section~\ref{s:res-phase} is not mentioned in the current section)
$g = \big(\hat f_h + \hat \omega \alpha \hat{\pdv{J}{z}} \hat f_z \big)|_{\varepsilon=0}$.
We have $\alpha \hat f_J = \hat \omega^{-1} g$.
We have the estimate~\eqref{e:d-J-z} $\pdv{J}{z} = O(h \ln^2 h \; \alpha^{-1})$, thus
\begin{equation} \label{e:est-g}
  g = \hat{f_h} + O(\hat h \ln \hat h) \hat{f_z} \qquad\text{in } \mathcal D_0.
\end{equation}
We will write $g(z, t, \lambda) = g(z, \varphi {=} \omega t, \lambda)$.
Let us rewrite
\begin{equation} \label{e:F-is-periodic}
  -F_s = \hat \omega^{-1} \langle g(z, \varphi, \lambda {=} \hat \omega^{-1} \varphi {-} Q) \rangle_{\varphi \in [0, 2s_2 \pi]} =
  (2\pi s_2)^{-1} \int_{t=0}^{s_2 T} g(z, t, \lambda {=} t {-} Q)dt.
\end{equation}
Let us denote by $l_0$ and $l_1$ the separatrices, let $l_0$ correspond to $\varphi \approx 0$ and $l_1$ to $\varphi \approx \pi$.
Let us split the phase curve of the unperturbed system for given $h, z$ into $2$ pieces $\hat l_{0}, \; \hat l_{1}$ close to the separatrices. We cut the phase curve by the line $y=x$ (cf. fig.~\ref{f:transversals}).
Let us define the coordinates $t_0=t,\; t_1=t - 0.5 T$ on $\hat l_0$ and $\hat l_1$, respectively. These coordinates are defined up to adding $i T, \; i \in \mathbb Z$ and are given by the time passed after crossing the transversals $\varphi=0$ and $\varphi=\pi$, respectively.
One may check that for $h \to 0$ these transversals approach some limit points on the separatrices, so the coordinates $t_0, t_1$ can be continued to the separatrices themselves.
Note that $i T = 2\pi \frac{i s_1}{s_2}$.
We can split the integral above as
\begin{align}
\begin{split}
  -F_s &=
  (2\pi s_2)^{-1} \int_{\hat l_{0}} \sum_{i=0}^{s_2-1} g\Big(
    z, t_0{=}t_0, \lambda {=} t_0 {+} 2 \pi \frac{i s_1}{s_2} {-} Q
  \Big)dt_0 \; +\\
  &+
  (2\pi s_2)^{-1} \int_{\hat l_{1}} \sum_{i=0}^{s_2-1} g\Big(
    z, t_1{=}t_1, \lambda {=} t_1 {+}  \pi \frac{(2 i + 1)s_1}{s_2} {-} Q
  \Big)dt_1.
\end{split}
\end{align}
As $s_1$ and $s_2$ are coprime, we have $\{ i s_1/ s_2 \bmod 1\}_{i=0}^{s_2-1} = \{ i/ s_2 \bmod 1\}_{i=0}^{s_2-1}$ and this rewrites as
\begin{align}
\begin{split} \label{e:potential}
  F_s =
  &- (2\pi)^{-1} \int_{\hat l_{0}} \Big\langle g \Big(
    z, t_0{=}t_0, \lambda {=} t_0 {+} 2 \pi \frac{i}{s_2} {-} Q
  \Big) \Big\rangle_{i=0, \dots, s_2-1} dt_0 \\
  &-
  (2 \pi)^{-1} \int_{\hat l_{1}} \Big\langle g \Big(
    z, t_1{=}t_1, \lambda {=} t_1 {+} \frac{\pi s_1}{s_2} {+} 2 \pi \frac{i}{s_2} {-} Q
  \Big) \Big\rangle_{i=0, \dots, s_2-1} dt_1.
\end{split}
\end{align}
Taking the derivative of the expression above yields
\begin{align}
\begin{split} \label{e:mp-derivative}
  \pdv{F_s}{Q}
  &=(2 \pi)^{-1} \int_{\hat l_{0}} \Big\langle \pdv{g}{\lambda}\/\Big(
    z, t_0{=}t_0, \lambda {=} t_0 {+} 2 \pi \frac{i}{s_2} {-} Q
  \Big) \Big\rangle_{i=0, \dots, s_2-1} dt_0 \;\\
  &+(2 \pi)^{-1}\int_{\hat l_{1}} \Big\langle \pdv{g}{\lambda}\/\Big(
    z, t_1{=}t_1, \lambda {=} t_1 {+} \frac{\pi s_1}{s_2} {+} 2 \pi \frac{i}{s_2} {-} Q
  \Big) \Big\rangle_{i=0, \dots, s_2-1} dt_1.
\end{split}
\end{align}
Note that as $\omega(\hat h, z) = s_2/s_1$, we have $\hat h \to 0$ for $s_2/s_1 \to 0$.
For fixed value $t$ of $t_0$ or $t_1$ by~\eqref{e:est-g} and using that the coordinates $h, t$ do not have singularities on the separatrices (i.e. $p(h, t)$ and $q(h, t)$ are smooth) we have
\begin{align}
\begin{split} \label{e:approx-g}
  g(z, t_0{=}t, \lambda)
  &= f_h(h{=}0, z, t_0{=}t, \lambda) + O(\hat h \ln \hat h), \\
  g(z, t_1{=}t, \lambda)
  &= f_h(h{=}0, z, t_1{=}t, \lambda) + O(\hat h \ln \hat h). \\
\end{split}
\end{align}
As $f_h(C) = 0$, the values of
$\max_\lambda \abs{f_h(h{=}0, z, t_0{=}t, \lambda)}$
and
$\max_\lambda \abs{f_h(h{=}0, z, t_1{=}t, \lambda)}$
exponentially decrease when $\abs{t} \to \infty$. Hence, the formulas above imply $C^0$-convergergence in
\begin{align}
\begin{split} \label{e:U-lim}
  F_s \; \xrightarrow[s_2/s_1 \to 0]{C^1} \;
  &-(2\pi)^{-1}\int_{l_0} \Big\langle f_h\Big(
    h{=}0, z, t_0{=}t_0, \lambda{=}t_0 {-} Q {+} 2\pi \frac{i}{s_2}, \varepsilon{=}0\Big) \Big\rangle_{i=0, \dots, s_2-1} dt_0 \\
  &-(2\pi)^{-1}\int_{l_1} \Big\langle f_h\Big(
    h{=}0, z, t_1{=}t_1, \lambda{=}t_1 {-} Q {+} 2\pi \frac{s_1 \bmod 2s_2}{2s_2} {+} 2\pi \frac{i}{s_2}, \varepsilon{=}0 \Big) \Big\rangle_{i=0, \dots, s_2-1} dt_1
\end{split}
\end{align}
in the domain
\begin{equation}
  \mathcal D_{F, 0} = \big\{
    Q, z \in [0, 2\pi] \times \mathbb C^n: \norm{z-z_0} < c
  \big\}.
\end{equation}
We can check that $-\pdv{F_s}{Q}$ converges to the $Q$-derivative of the right-hand side in the same way, using~\eqref{e:mp-derivative}. Note that $\pdv{f_h}{\lambda}\/(C) = 0$ and, similarly to~\eqref{e:est-g}, we have $\pdv{g}{\lambda} = \pdv{\hat f_h}{\lambda} + O(h \ln h) \pdv{\hat{f_z}}{\lambda}$.
Finally, by Cauchy formula $\pdv{F_s}{z}$ and $\pdv[2]{F_s}{z}$ converge to $\pdv{}{z}$ and $\pdv[2]{}{z}$, respectively, of the right-hand side of~\eqref{e:U-lim} in $\mathcal D_F$. Similarly, $\pdv{F_s}{z}{Q}$ converges to the $\pdv{}{z}{Q}$ of the right-hand side of~\eqref{e:U-lim}. This shows that $\pdv{F_s}{z}$ converges in $C^1$ to the $z$-derivative of the right-hand side of~\eqref{e:U-lim}.

\begin{proof}[Proof of Lemma~\ref{l:lim}]
  The first and the second parts of Lemma~\ref{l:lim} follows from~\eqref{e:F-is-periodic}. Periodicity follows from the fact that $g(\lambda)$ is $2\pi$-periodic.
  As the right-hand side of~\eqref{e:F-is-periodic} is defined in $\mathcal D_0$, we can continue $F_s$ in the domain $\mathcal D_{F, 0}$.

  By the estimate for the error of the trapezoidal integration rule we have
  \begin{equation}
    \Big\langle f_h\Big(
      h{=}0, z, t_0{=}t_0, \lambda {=} t_0 {+} 2 \pi \frac{i}{s_2} {-} Q
    \Big) \Big\rangle_{i=0}^{s_2-1}
    =
    \langle f_h(h{=}0, z, t_0{=}t_0, \lambda) \rangle_{\lambda}
    + O(s_2^{-2}).
  \end{equation}
Together with~\eqref{e:U-lim} this implies
$\norm{F_s - \Theta_3}_{C^0} \to 0$ in $\mathcal D_{F, 0}$ for $s_2 \to \infty, s_2/s_1 \to 0$.
By Cauchy formula this also means
$\norm{\pdv{F_s}{z} - \pdv{\Theta_3}{z}}_{C^0} \to 0$
and
$\norm{\pdv[2]{F_s}{z} - \pdv[2]{\Theta_3}{z}}_{C^0} \to 0$
in $\mathcal D_{F}$.
We have
\begin{equation}
  \pdv{}{Q}\/\Big\langle f_h\Big(
    h{=}0, z, t_0{=}t_0, \lambda {=} t_0 {+} 2 \pi \frac{i}{s_2} {-} Q
  \Big) \Big\rangle_{i=0}^{s_2-1}
  =
  - \Big\langle \pdv{f_h}{\lambda}\Big(
    h{=}0, z, t_0{=}t_0, \lambda {=} t_0 {+} 2 \pi \frac{i}{s_2} {-} Q
  \Big) \Big\rangle_{i=0}^{s_2-1}.
\end{equation}
Applying trapezoidal rule argument again yields
$\norm{\pdv{F_s}{Q}}_{C^0} \to 0$ in $\mathcal D_{F, 0}$
for $s_2 \to \infty, s_2/s_1 \to 0$.
Cauchy formula implies $\norm{\pdv{F_s}{Q}{z}}_{C^0} \to 0$ in $\mathcal D_F$, this proves the third part of the lemma.

  The last part of the lemma follows from~\eqref{e:U-lim} and $C^1$-convergence of $z$-derivatives in this formula established above. As the right-hand sides of~\eqref{e:U-lim} depend on $s_1 \bmod 2 s_2$ and not on $s_1$, there is a finite set of possible right-hand sides and we take this set as $\mathcal F_{s_2}$.
\end{proof}

\newpage
\section{Crossing resonant zones: proofs} \label{s:res-zones-proofs}
\subsection{High-numerator resonances: proof} \label{ss:weak-res-proof}
\begin{proof}[Proof of Lemma~\ref{l:est-weak-res}]
  Let us apply Lemma~\ref{l:aux} with $z_0$ as in the statement of lemma and $C_\gamma \ge 8\pi$.
  By~\eqref{e:alphabet} and the bound on $s_1$ from the statement of lemma we have
  \begin{equation} \label{e:loc-0267}
    \hat \omega^{-2} s_1 \beta(z_0) \lesssim C_\rho^{-0.5}.
  \end{equation}
  This means that for large enough $C_\rho$ the terms containing $H_7$ in~\eqref{e:res-passing} are small. The terms containing $u$ are also small for small enough $\varepsilon$.
  For large enough $S_2$ there exists $\delta > 0$ such that $F_s > \delta$ if $s_2 > S_2$ (by Lemma~\ref{l:lim}).
  This means $P' < - 0.5 \delta$ in~\eqref{e:res-passing}.

  Lemma~\ref{l:aux} gives the domain $\mathcal D$.
  We need~\eqref{e:hat-h-bound} to hold in this domain. This holds if $c_Z < d_Z$, where $c_Z$ is from Lemma~\ref{l:aux} and $d_Z$ is from Lemma~\ref{l:d-z}. We reduce $c_Z$ if needed so that $c_Z < d_Z$.
  Denote $D(I, z)$ by the equation $I = \hat I(z) + D \alpha(z) \hat \omega^{0.5}$.
  This also gives the function $D(P, Q, Z, \tau)$ defined in $\mathcal D$.
  By~\eqref{e:p-q-J-gamma} we have
  \begin{equation} \label{e:D-P}
    |D - P| = O(s_1\hat \omega^{-1.5} \beta) \lesssim \sqrt{\hat \omega} \qquad \text{in } \mathcal D.
  \end{equation}

  Denote by
  $\mathcal D' \subset \mathcal D$ the subdomain given by the additional restriction
  $|P| < D_p + 1$.
  Denote $Y = (P, Q, Z, \tau)$ and consider a solution $Y(\tau)$ of~\eqref{e:res-passing} with the initial condition $Y_0$, $Y_0 = (P_0, Q_0, Z_0, \tau_0)$, obtained from $X_0$ after the coordinate change of Lemma~\ref{l:aux}.
  We can add $2\pi k$ to $\varphi(X_0)$ so that $\gamma = \varphi - (s_2/s_1) \lambda \in [-\pi, \pi]$, then by~\eqref{e:p-q-J-gamma} and~\eqref{e:loc-0267} we have $Q_0 \in [-1.1\pi \hat \omega^{-1}, 1.1\pi \hat \omega^{-1}]$.
  By this estimate on $Q_0$ and~\eqref{e:D-P} $Y_0 \in \mathcal D'$ for small enough $\omega_0$.

  As $P' < 0$, the solution $Y(\tau)$ cannot leave $\mathcal D'$ through $P = D_p + 1$.
  Time $\tau$ required to leave $\mathcal D'$ without reaching the hypersurface $D = -D_p-1$ (i.e., via $\norm{Z} = c_Z$ or $|Q| = \hat \omega^{-1} C_\gamma / 4$) is
  $\gtrsim \hat \omega^{-1}$ by~\eqref{e:res-passing}.
  On the other hand, after time $\tau \lesssim D_p$ we will have $P < -D_p - 1$, thus (for small enough $\omega_0$) the solution $Y(\tau)$ leaves $\mathcal D'$ through $P = -D_p - 1$.  Denote the time $\lambda$ when this happens by $\lambda_{out}$.
  At the moment $\lambda_{out}$ we have $P = -D_p - 1$ and thus $D < -D_p \le -D_p'$.
  Hence, $\lambda_1 < \lambda_{out}$.

  Let us now obtain estimate for the time passed before crossing $D = -D_p'$, i.e. for $\lambda_1 - \lambda_0$. Such estimate can be obtained from~\eqref{e:D-P}, but we need better estimate.
  Given $\lambda_2 \in [\lambda_0, \lambda_{out} - 2\pi s_1]$, set $\lambda_3 = \lambda_2 + 2\pi s_1$.
  Denote by $\Delta \tau$ the time $\tau$ between $\lambda_2$ and $\lambda_3$.
  We have $\Delta \tau \sim s_1 \hat \omega^{-0.5} \beta(z_0)$.
  Let us compare $D(\lambda_3)$ with $D(\lambda_2)$.
  We have from~\eqref{e:res-passing}
  \[
    P(\lambda_3) \le  P(\lambda_2) - 0.5 \delta \Delta \tau,
    \;\;
    |P(\lambda_3) - P(\lambda_2)| + |Q(\lambda_3) - Q(\lambda_2)| \lesssim \Delta \tau,
    \;\;
    \norm{z(\lambda_3) - z(\lambda_2)} \lesssim \alpha(z_0) \hat \omega^{-1.5} \Delta \tau.
  \]
  By~\eqref{e:loc-0267} and~\eqref{e:d-D-minus-P} we have $\pdv{D}{P} > 0.9$ for small enough $\omega_0$.
  Together with the estimates above and~\eqref{e:d-D-minus-P}, this implies
  \[
    D(X(\lambda_3)) \le D(X(\lambda_2)) -
    \Delta \tau \Big(
      0.5 \delta + O(s_1 \hat \omega^{-1.5} \beta(z_0))
      + O(s_1 \hat \omega^{-4} \varepsilon)
    \Big) \le
    D(X(\lambda_2)) - 0.25 \delta \Delta \tau.
  \]
  As $\delta \sim 1$, this gives the estimate for $\lambda_1 - \lambda_0$ from the statement of the current lemma:
  \[
    \lambda_1 - \lambda_0 - 2\pi s_1 \lesssim s_1 D'_p/\Delta \tau
    \sim D'_p \hat \omega^{0.5} \beta(z_0)^{-1}
    \sim D'_p \hat \omega^{0.5} \varepsilon^{-1} \alpha(z_0).
  \]
\end{proof}

\subsection{Lemma on model system} \label{s:lem-model}
In this subsection we state a general lemma that will be later applied to study resonance crossing described by~\eqref{e:res-passing}.
Consider a system
\begin{align} \label{e:rc-perturbed}
\begin{split}
  p' &= -\pdv{H}{q}\/(p, q, w) + v_p(p, q, w, \tau), \\
  q' &= \pdv{H}{p}\/(p, q, w) + v_q(p, q, w, \tau), \\
  w' &= v_w(p, q, w, \tau), \\
  \tau' &= 1,
\end{split}
\end{align}
where
\begin{equation}
  H = H_0 + \Delta H(p, q, w), \qquad
  H_0 = p^2/2 + V(q),
\end{equation}
and
\begin{equation} \label{e:dH-v}
  \norm{\Delta H}_{C^2} < \varepsilon_1,
  \qquad
  \norm{v_p}_{C^0}, \norm{v_q}_{C^0}, \norm{v_w}_{C^0} < \varepsilon_2
\end{equation}
and
\begin{equation} \label{e:V-c}
    V(q) = V_{c} q + V_{per}(q),
\end{equation}
where $V_c > 0$ is a constant and $V_{per}(q)$ is $2\pi$-periodic.
We will call the system~\eqref{e:rc-perturbed} the \emph{perturbed system}.
We will call this system without the $v$ terms the \emph{intermediate system}. It is an autonomous Hamiltonian system with the Hamiltonian $H(p, q, w)$.
We will call the Hamiltonian system given by $H_0$ the \emph{unperturbed system}.
One can find analysis of the unperturbed system in~\cite[Section~9.2]{neishtadt2014}.
We will assume that the function $V$ satisfies Condition $B'$ introduced in Subsection~\ref{s:cond-B}.
Note that the saddles of the unperturbed system correspond to the local maxima of $V$.
This also holds for the intermediate system if $\varepsilon_1$ is small enough.

\begin{lemma} \label{l:near-sep}
  Fix $V(q)$ as above. Then for any small enough $\varepsilon_1 > 0$ and any large enough $D_{p, 1}, D_{p, 2} > 0$ with $D_{p, 1} < D_{p, 2}$
  for any large enough (compared with $D_{p, 2}$) $C_p$ for any large enough (compared with $C_p$) $C_q$ there exists $C > 1$  such that for any $c_Z > 0$ and any $\varepsilon_2 > 0$ with
  \begin{equation}
    \varepsilon_2 < \varepsilon_1, \qquad \varepsilon_2 (1+|\ln \varepsilon_2|) < 0.5 c_Z C^{-1}
  \end{equation}
  the following holds.

  Consider the unperturbed system in the domain
  \begin{equation} \label{e:D-near-sep}
    \mathcal D = \big\{ w, p, q, \tau \in \mathbb R^{n+3}: \;
      \norm{w} \le c_Z, \; p \in [-C_p, C_p], \; q \in [-C_q, C_q]
    \big\}.
  \end{equation}
  For any $\Delta H$, $v_p$, $v_q$, $v_w$ that satisfy~\eqref{e:dH-v} in this domain with $\varepsilon_1$, $\varepsilon_2$ fixed above, also consider the intermediate and the perturbed systems in this domain.
  Let
  \[
    X_0 = (p_0, q_0, w_0, \tau_0)
  \]
  denote some initial data with
  \[
    p_0 \in [D_{p, 1}, D_{p, 2}], \qquad
    q_0 \in [-\pi, \pi], \qquad
    \norm{w_0} < 0.5 c_Z.
  \]
  If $V$ has no local maxima, set $\Delta h_0 = 1$.
  Otherwise, let $C_i(w)$ be a saddle of the intermediate system such that $H(C_i(w_0), w_0)$ is as close as possible to $H(p_0, q_0, w_0)$ and set
  \[
    \Delta h_0 = |H(p_0, q_0, w_0) - H(C_i(w_0), w_0)|.
  \]
  We will assume
  \[
    \Delta h_0 > 2 C \varepsilon_2.
  \]
  Let $X(\tau) = (p(\tau), q(\tau), w(\tau), \tau)$ denote the solution of the perturbed system with initial data $X_0$.
  Then there exists $\tau_1 > \tau_0$ such that %
  \begin{equation}
    p(\tau_1) = -C_p, \qquad
    \tau_1 - \tau_0 < C (1 + |\ln \Delta h_0|)
  \end{equation}
  and for any $\tau \in [\tau_0, \tau_1]$ we have
  \begin{equation}
    X(\tau) \in \mathcal D, \qquad
    \abs{H(X(\tau)) - H(X(\tau_0))} < C \varepsilon_2, \qquad
    \norm{w(\tau) - w(\tau_0)} < C \varepsilon_2 (1 + |\ln \Delta h_0|).
  \end{equation}
\end{lemma}

\subsection{Proof of the lemma on model system}
Let us state a lemma that will be used to prove Lemma~\ref{l:near-sep}.
Denote by $U$ some neighborhood of $(0, 0, 0)$ with $\diam U < 1$.
Consider the Hamiltonian system given in $U$ by the Hamiltonian $H = p^2/2 - a q^2 + \Delta H(p, q, w)$, $a > 0$
and its perturbation given by a vector field $u$:
\begin{align} %
\begin{split}
  \dot p &= 2a q - \pdv{\Delta H}{q}\/(p, q, w) + u_p(p, q, w, t), \\
  \dot q &= p + \pdv{\Delta H}{p}\/(p, q, w) + u_q(p, q, w, t), \\
  \dot w &= u_z(p, q, w, t).
\end{split}
\end{align}
Assume that in $U$ we have
\begin{equation} \label{e:4a}
  \norm{\Delta H}_{C^2} < \varepsilon_1,
  \qquad
  \norm{u_p}_{C^0}, \norm{u_q}_{C^0}, \norm{u_w}_{C^0} < \varepsilon_2.
\end{equation}
The constants $\varepsilon_1$ and $\varepsilon_2$ are assumed to be small enough compared to $\min(a, a^{-1})$.
The Hamiltonian system has a saddle $C(w) = (p_C(w), q_C(w))$ with $p_C, q_C = O(\varepsilon_1)$.
Assume $C(w) \in U$ for all values of $w$ encountered in $U$. Denote $h(p, q, w) = H(p, q, w) - H(C(w), w)$.
\begin{lemma} \label{l:perturb-saddle}
  There exists $C_1 > 0$ such that for any
  $\Delta H, u_p, u_q, u_w$ as above
  the following holds. Given any initial data $(p_0, q_0, w_0, t_0)$ in $U$ such that
  \[
    \abs{h(p_0, q_0, w_0)} \ge C_1 \varepsilon_2,
  \]
  consider the solution $p(t), q(t), w(t), t$ of the perturbed system starting at $p_0, q_0, w_0, t_0$.
  Then this solution exits $U$ at some time $t_1 > t_0$ and
  for any $t \in [t_0, t_1]$ we have the estimates
  \begin{equation} %
    \abs{h(p(t), q(t), w(t)) - h(p_0, q_0, w_0)} < 0.5 C_1 \varepsilon_2,
    \qquad
    t_1 - t_0 < C_1(1 + \abs{\ln \abs{h(p_0, q_0, w_0)}}).
  \end{equation}
\end{lemma}
\begin{proof}
  We will assume $H(p_0, q_0, w_0) > H(p_C(w_0), q_C(w_0), w_0)$, the proof is similar when the opposite inequality holds.
  All $O$-estimates in this proof will be uniform in $p_0$, $q_0$, $w_0$, $\Delta H$, $u_p$, $u_q$, $u_w$.
  Denote $\tilde p = p - p_C(w)$, $\tilde q = q - q_C(w)$.
  Let us denote by $\tilde H$ the Hamiltonian $H - H(p_C(w), q_C(w), w)$ rewritten in the shifted coordinates. We can write $\tilde H = \tilde p^2/2 - a \tilde q^2 + \psi(\tilde p, \tilde q, w)$, where $\psi = O(\varepsilon_1)$.
  As at the points $(0, 0, w)$ we have
  $\tilde H = \pdv{\tilde H}{\tilde p} = \pdv{\tilde H}{\tilde q} = 0$,
  $\psi$ does not contain constant or linear terms with respect to $\tilde p, \tilde q$. Hence, we have
  \begin{equation} \label{e:ps-psi}
    \tilde H = \tilde p^2/2 - a \tilde q^2 + \psi(\tilde p, \tilde q, w),
    \qquad
    \psi, \pdv{\psi}{w} = \varepsilon_1 O(\tilde p^2 + \tilde q^2),
    \;
    \pdv{\psi}{\tilde p}, \pdv{\psi}{\tilde q} = \varepsilon_1 O(\abs{\tilde p} + \abs{\tilde q}).
  \end{equation}
  Denote $\tilde p(t) = p(t) - p_C(w(t))$, $\tilde q(t) = q(t) - q_C(w(t))$ and $h(t) = \tilde H(\tilde p(t), \tilde q(t), w(t))$.
  Let us take the largest $t_2 > t_0$ such that for any $t \in [t_0, t_2)$ we have
  \begin{equation} \label{e:ps-ind}
    (p(t), q(t), w(t)) \in U, \qquad h(t) > \varepsilon_2.
  \end{equation}
  Then for any $t \in [t_0, t_2]$ at the point $(\tilde p, \tilde q, w) = (\tilde p(t), \tilde q(t), w(t))$ we have
  \[
    h = h(t) = (1 + O(\varepsilon_1)) \tilde p^2/2 - (a+O(\varepsilon_1)) \tilde q^2.
  \]
  Hence,
  \[
    \abs{\tilde p} \ge \sqrt{a \tilde q^2 + h} \ge c_1 \max(\abs{\tilde q}, \sqrt{\varepsilon_2})
  \]
  for some $c_1 \in (0, 1)$.
  Hence, $\tilde p(t)$ has the same sign for all $t \in [t_0, t_2]$. Without loss of generality we will assume it to be positive.
  We have
  $\dot{\tilde q} = \pdv{\tilde H}{\tilde p} + \varepsilon_2 (u_q - \dv{q_C}{w}u_w) = \tilde p + \varepsilon_1 O(\abs{\tilde p} + \abs{\tilde q}) + O(\varepsilon_2)$.
  Therefore,
  \begin{equation} \label{e:est-dot-q}
    \dot{\tilde q} \ge 0.25 c_1 (\abs{\tilde p} + \abs{\tilde q}),
    \qquad
    \dot{\tilde q} \ge 0.5\sqrt{a \tilde q^2 + h}.
  \end{equation}

  By~\eqref{e:ps-psi} we have
  $\abs{\pdv{\tilde H}{\tilde p}}, \abs{\pdv{\tilde H}{\tilde q}} \le O(1) (\abs{\tilde p} + \abs{\tilde q})$
  and thus
  $\dot{\tilde q}^{-1} \pdv{\tilde H}{\tilde p}, \dot{\tilde q}^{-1} \pdv{\tilde H}{\tilde q} = O(1)$.
  As $\dot{\tilde q} > 0$, we can use $\tilde q$ as an independent variable instead of $t$. Denote by $'$ the derivative with respect to $\tilde q$. We have
  \[
    h' = \dot{\tilde q}^{-1} \pdv{\tilde H}{q}\/\Big(u_q - \dv{C_q}{w} u_w\Big)
       + \dot{\tilde q}^{-1} \pdv{\tilde H}{p}\/\Big(u_p - \dv{C_p}{w} u_w\Big)
       + \pdv{\tilde H}{w} u_w
       = O(\varepsilon_2).
  \]
  Thus for $t \in [t_0, t_2]$ we have the estimate $\abs{h(t) - h(t_0)} = O(\varepsilon_2)$ and this estimate does not depend on $t_2$. Denote $h_0 = h(t_0)$. We have
  \[
    h_0 = \tilde H(\tilde p_0, \tilde q_0, w_0)
    = H(p_0, q_0, w_0) - H(C_p(w_0), C_q(w_0), w_0) \ge C_1 \varepsilon_2.
  \]
  Hence, for large enough $C_1$ we have $h(t_2) > 0.5 h_0 > 0$ and thus the solution exists $U$ at the time $t_2$. This means that $t_1$ exists and $t_1 = t_2$. This also proves the estimate for the change of $H$.

  Using~\eqref{e:est-dot-q}, we can estimate (the details are given below)
  \[
    t_1 - t_0 = \int_{q(t_0)}^{q(t_1)} \frac{d \tilde q}{\dot{\tilde q}} \le 2 \int (a \tilde q^2 + 0.5 h_0)^{-0.5} d\tilde q = O(\abs{\ln h_0} + 1).
  \]
  The estimate for the integral above can be obtained by splitting it into two parts, with $\abs{\tilde q} \le \sqrt{h_0}$ and with $\abs{\tilde q} > \sqrt{h_0}$.
  When $\abs{\tilde q} \le \sqrt{h_0}$, we use $(a \tilde q^2 + 0.5 h_0)^{-0.5} = O(h_0^{-0.5})$, so this part is $O(1)$.
  When $\abs{\tilde q} > \sqrt{h_0}$, we use $(a \tilde q^2 + 0.5 h_0)^{-0.5} = O(\tilde q^{-1})$, so this part is $O(\abs{\ln h_0} + 1)$.
  This gives the estimate for $t_1 - t_0$ and thus completes the proof of the lemma.
\end{proof}

\begin{proof}[Proof of Lemma~\ref{l:near-sep}]
  We will consider the case where the unperturbed system has saddles, the other case is much simpler.
  Let us first consider the unperturbed system in
  \[
    \big\{ w, p, q, \tau \in \mathbb R^{n+3}: \;
      \norm{w} \le c_Z
    \big\}.
  \]
  Note that for $p > 0$ the value of the function $L = p^2/2 + V_{per}(q) = H_0 - V_c q$ decreases along the solutions of the unperturbed system, $\dot L = -p V_c < 0$. For small enough $\varepsilon_1$ (this implies that $\Delta H$ and $v$ are small) and $p > 1$
  we also have $\dot L < 0$ along the solutions of intermediate and perturbed systems.
  The function $L$ is $2 \pi$-periodic in $q$ and its contour lines such that $p > 1$ on the whole contour line provide transversals to solutions of all three systems.

  Let us take a contour line $\mathcal{D}_p$ of $L$ such that $p > 1$ on this line, let $D_{p, 0}$ be the maximum of $p$ on this line.
  We assume $D_{p, 1} > D_{p, 0}$.
  Given $D_{p, 1}$ and $D_{p, 2}$, take $C_p > D_{p, 2}$ such that there is a contour line of $L$ with $p \in (D_{p, 2}, C_p)$ on the whole contour line. This guarantees that a solution (of any of the three systems) starting with $p \le D_{p, 2}$ does not cross the line $p=C_p$. Let us now restrict to the domain
  \[
    \mathcal D_+ = \big\{ w, p, q, \tau \in \mathbb R^{n+3}: \;
      \norm{w} \le c_Z, \; |p| \le C_p
    \big\}.
  \]

  By~\eqref{e:dH-v} the set
  \[
    \mathcal H(w, \Delta H) = \big\{H(p, q, w): \; (p, q, w, \tau) \in \mathcal D_+, \; q \in [-\pi, \pi] \big\}
  \]
  is bounded (uniformly in $w$ and $\Delta H$).
  We will assume $\varepsilon_1 < 0.25 V_c$. Then we have
  \begin{equation} \label{e:H-inc}
    H(p, q+2\pi) \ge H(p, q) + \pi V_c.
  \end{equation}
  As $\mathcal H$ is bounded, this means that the values of $|q|$ on the contour lines given by $H(p, q, w) = h_0 + h_1$, $h_0 \in \mathcal H$, $h_1 \in [-1, 1]$ is bounded by some constant (uniformly in $w$ and $\Delta H$), take $C_q$ equal to this constant.
  This choice of $C_q$ gives the following property that will be used later: if we start in $\mathcal D_+$ with $q \in [-\pi, \pi]$, and the value of $H$ changes by at most $1$ (while still being in $\mathcal D_+$), the value of $q$ stays in $(-C_q, C_q)$.

  Fixed points of the unperturbed system correspond to extrema of $V$, with maxima of $V$ corresponding to saddles and minima of $V$ corresponding to centers.
  For each saddle $C_{V, i} = (0, q_{C, i})$ of the unperturbed system let us fix a small neighborhood $U_i$ in the space with coordinates $p, q, w$ (these neighborhoods will respect $2\pi$-periodicity of the unperturbed system, i.e. neighborhoods of saddles that differ by $2\pi k$ will be shifts of each other; this means that we only need to construct such neighborhoods for saddles with $q_{C, i} \in [0, 2\pi]$) in the following way. First, let us fix preliminary neighborhoods $U_{i, 0} = U_{i, 0}^{p, q} \times \{w: \norm{w} < c_Z\}$ such that there is $c_{sep} > 0$ such that the values of $H_0$ inside different $U_{i, 0}$ are separated by at least $c_{sep}$.
  Now for fixed $i$ let us define $U_i$.
  Let $q_i = q - q_{C, i}$, we can write $V = V(q_{C, i}) - a_i q_i^2 + O(q_i^3)$.
  We will now use Lemma~\ref{l:perturb-saddle} together with the notation defined there, we add tilde to expressions from this lemma to distinguish them.
  Let us apply this lemma to $\tilde q = q_i$, $\tilde p = p$, $\tilde w = w$, $\tilde a = a_i$ and $\tilde U = U_{i, 0}$, it gives us $\tilde \varepsilon_1$ such that the lemma can be applied if in $U_{i, 0}$ we have
  $\norm{\widetilde{\Delta H}}_{C_2} < \tilde \varepsilon_1$.
  Now take $U_i \subset U_{i, 0}$ such that inside $U_i$ we have $\norm{V - V(q_{C, i}) + a_i q_i^2}_{C^2} < 0.5 \tilde \varepsilon_1$ and assume $\varepsilon_1 < 0.5 \tilde \varepsilon_1$ (note that this restriction on $\varepsilon_1$ depends only on $V$, this will also hold for further restrictions on $\varepsilon_1$).
  Denote by $C_i(w)$ the saddle of the intermediate system near $C_{V, i}$.
  If needed, let us futher increase $\varepsilon_1$ so that we have $C_i \in 0.5 U_i$ for any $\Delta H$ that satisfies~\eqref{e:dH-v}, here $0.5 U_i$ is the image of $U_i$ under $w$-dependent homothety with center $C_i(w)$ and ratio $0.5$.
  Take
  \[
    \widetilde{\Delta H} = \Delta H + V - V(q_{C, i}) + a_i q_i^2
  \]
  and $\tilde u =v$, we have $H_0 + \Delta H = p^2/2 - a_i q_i^2 + \widetilde{\Delta H} + const$ in $U_i$. As $\widetilde{\Delta H}$ satisfies~\eqref{e:4a} in $U_i$ (with $\tilde \varepsilon_1$ instead of $\varepsilon_1$), we can apply Lemma~\ref{l:perturb-saddle} to describe the movement inside $U_i$ for small enough $\varepsilon_2$. The words "small enough" here give another upper bound on $\varepsilon_1$.

  By the choice of $U_i$ the values of $H_0$ in different sets $U_i$ are separated by at least $c_{sep}$.
  For small enough $\varepsilon_1$ the values of $H$ in different sets $U_i$ are separated by at least $c_{sep}/2$. Take small enough $h_b \in (0, c_{sep}/8)$ such that for any $i, \Delta H, w$ the set $|H - H(C_i(w), w)| \le 2h_b$ intersects $\partial U_i$ by four disjoints sets near intersections of separatrices of $C_i$ with $\partial U_i$.

  Let us suppose we are given some initial data.
  Denote
  \begin{equation}
    h(p, q, w) = H(p, q, w) - H(C_i(w), w).
  \end{equation}
  We will assume $h_0 < h_b$, the case $h_0 \ge h_b$ (meaning the initial condition is far from separatrices of the intermediate system) is easier and we omit it.
  We consider the case $h_0 > 0$, the proof is similar when $h_0 < 0$.
  Let $Z$ be the stripe given by
  \begin{equation}
    Z = \big\{(p, q, w, \tau) \in \mathcal D: \;
      h \in [0.5 h_0, 2h_b]
    \big\}
  \end{equation}
  (cf. Figure~\ref{f:zones}).
  \begin{figure}[ht]
    \centering
    \includegraphics[width=0.7\textwidth]{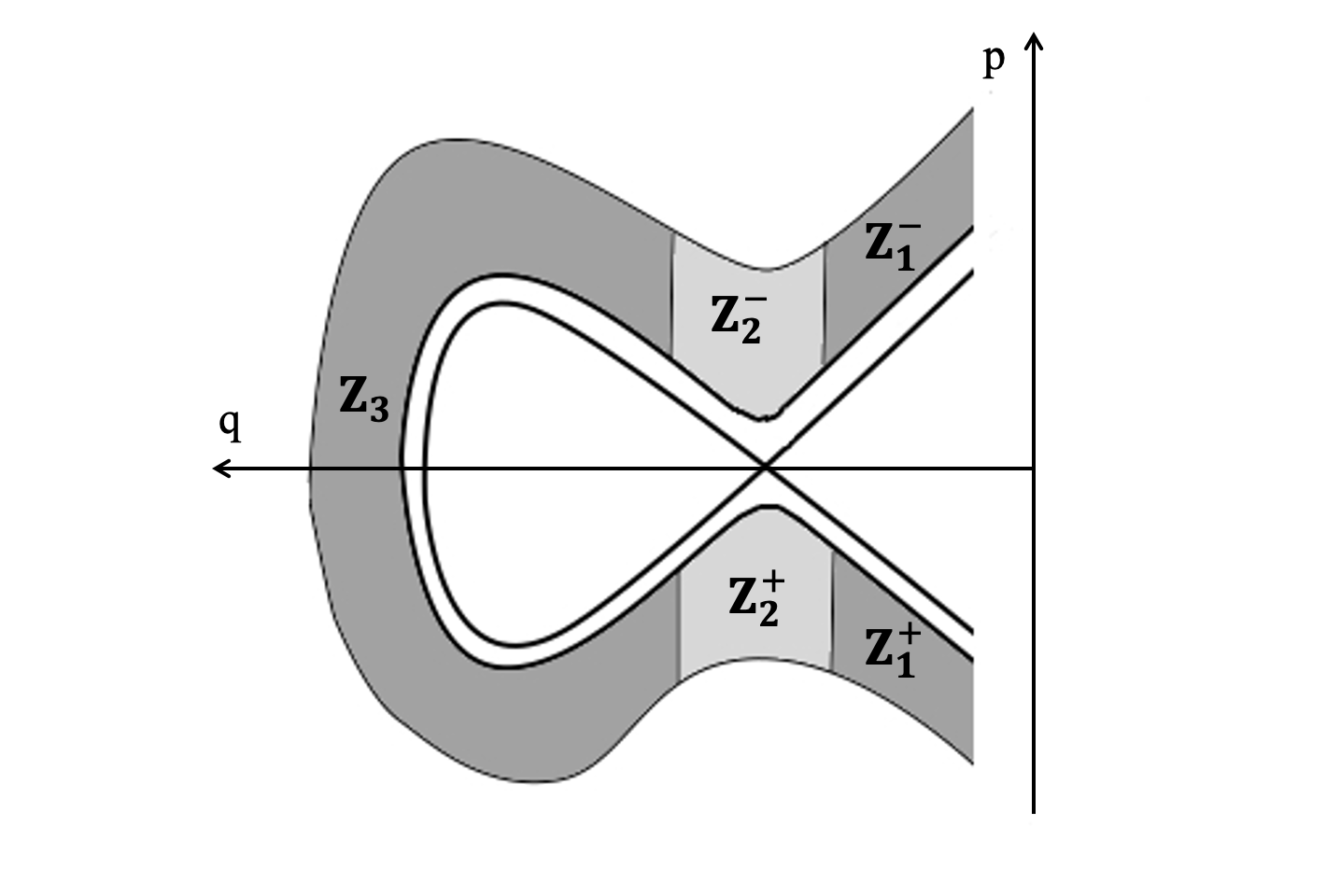}
    \caption{The zones $Z_i$.}
    \label{f:zones}
  \end{figure}
  For small enough $h_b$ we have $q \in (-C_q, C_q)$ in $Z$.
  Also, as $2 h_b < c_{sep}/4$, the zone $Z$ does not intersect $U_j$ for $j \ne i$.

  We can split $Z$ into a union of three zones $Z_1^-, Z_1^+, Z_3$ far from $C_i$ and two zones $Z_2^-, Z_2^+$ inside $U_i$. Solutions of the intermediate system visit these zones in the following order: $Z_1^-$, $Z_2^-$, $Z_3$, $Z_2^+$, $Z_1^+$. Thus, for small $\varepsilon_2$ solutions of the perturbed system starting in each zone can only leave this zone into the next zone or through the boundary of $Z$, but not into the previous zone.

  Take one of the zones $Z_1^-, Z_3, Z_1^+$.
  It is easy to show that solutions of the perturbed system starting at any point inside this zone exit it after time $O(1)$ passes.
  As $h' = O(\varepsilon_2)$ and $w' = O(\varepsilon_2)$, both $h$ and $w$ change by at most $O(\varepsilon_2)$ while passing this zone.

  Now take one of the zones $Z_2^-, Z_2^+$ inside $U_i$.
  When selecting $U_i$ above, we have checked that Lemma~\ref{l:perturb-saddle} can be applied to describe the solutions of the perturbed system inside $U_i$. By this lemma any orbit starting in our zone leaves this zone after time $O(\abs{\ln h_0} + 1)$ passes and $h$ changes by $O(\varepsilon_2)$. From the estimate on the time spent in this zone we conclude that $w$ changes by at most $O(\abs{\ln h_0} + 1) \varepsilon_2$.

  Thus the total time spent in $Z$ and changes in $h$ and $w$ before leaving $Z$ are bounded by  $O(|\ln h_0| + 1)$, $O(\varepsilon_2)$ and $O(\abs{\ln h_0} + 1) \varepsilon_2$, respectively. Take $\tau_1$ to be the moment when the solution leaves $Z$. Taking $C$ much greater than the (uniform) constants in these $O$-estimates, we obtain the estimates for $\tau_1 - \tau_0$ and the change in $h$ and $w$ from the statement of the lemma.
  For large enough $C$ we also have $h(\tau) \in [h_0 - C \varepsilon_2, h_0 + C \varepsilon_2] \subset (0.5 h_0, 2 h_b)$ and $\norm{w(\tau)-w_0} \le  C(|\ln \varepsilon_2| + 1) \varepsilon_2 < 0.5c_Z$ (due to the inequalities on $\varepsilon_2$ from the statement of the lemma) for large enough $C$.
  Recall that we have $q(\tau) \in (-C_q, C_q)$ in $Z$.
  This implies that our solution can only leave $Z$ by crossing one of the lines $p = C_p$, $p=-C_p$. However, the choice of $C_p$ above prohibits crossing $p = C_p$, thus the solution crosses $p = -C_p$.
  This completes the proof.
\end{proof}

\subsection{Applying the lemma on model system}
In this section we apply Lemma~\ref{l:near-sep} to system~\eqref{e:res-passing} (with $P, Q, Z$ in ~\eqref{e:res-passing} corresponding to $p, q, w$ in~\eqref{e:rc-perturbed} in the same order) in such a way that the estimates of this lemma will be uniform for all resonances close enough to separatrices.
We also take care of the fact that the main part of the Hamiltonian in~\eqref{e:res-passing} depends on $Z$ by taking its value at $Z=0$ as the main part and considering the difference as part of the Hamiltonian perturbation.

We will use the notation $V_F$ introduced in Subsection~\ref{s:cond-B}.
Without the $u$ terms the system~\eqref{e:res-passing} is Hamiltonian with the Hamiltonian
\begin{equation} \label{e:H-without-u}
  H_r(P, Q, Z) = P^2/2 + V(Q, z_0+\hat \omega Z) + \hat \omega^{-2} s_1 \beta H_7,
\end{equation}
where $V = V_{F_s}$.
\begin{lemma} \label{l:near-sep-2}
  Given $z_* \in \mathcal Z_B$,
  for any $C_{s_1}$ for any large enough $D_{p, 2} > D_{p, 1} > 0$ there exist
  \[
    C_p > D_{p, 2}, \qquad
    C_\rho, C, C_q > 1, \qquad
    c_z, c_Z, \omega_0, \varepsilon_0 > 0
  \]
  such that
  for any $z_0$ with $\norm{z_0 - z_*} < c_z$
  for any
  $\hat \omega = s_2 / s_1 \in (0, \omega_0)$,
  and $\varepsilon < \varepsilon_0$
  with
  \[
    |s_1| < C_{s_1} \ln^2 \varepsilon, \qquad
    \hat h(z_0) > C_\rho \varepsilon |\ln^5 \varepsilon|
  \]
  we have the following.

  $1.$ We can apply Lemma~\ref{l:aux} and the system~\eqref{e:res-passing} given by this lemma is defined in the domain
  \begin{equation}
    \mathcal D = \big\{ (Z, P, Q, \tau) \in \mathbb R^{n+3}: \;
      \norm{Z} < c_Z, \; |P| < C_p, \; |Q| < \pi \hat \omega^{-1} + C_q
    \},
  \end{equation}
  where $Z = \hat \omega^{-1}(z - z_0)$.

  $2.$
  Set $\varepsilon_2 = C \alpha(z_0) \hat \omega^{-1}$.
  Let
  \[
    X_0 = (P_0, Q_0, Z_0, \tau_0)
  \]
  denote some initial data with
  \[
    P_0 \in [D_{p, 1}, D_{p, 2}], \qquad
    Q_0 \in [-\hat \omega^{-1} \pi - 1, \hat \omega^{-1} \pi+1], \qquad
    \norm{Z_0} < c_Z/2.
  \]
  If $V_{F_s}(Q, z_*)$ (here $F_s$ is from~\eqref{e:res-passing}) has no local maxima, set $\Delta h_{r, 0} = 1$.
  Otherwise, let $C_i(Z)$ be a saddle of the intermediate system such that $H(C_i(Z_0), Z_0)$ is as close as possible to $H(p_0, q_0, Z_0)$ and set
  \[
    \Delta h_{r, 0} = |H_r(p_0, q_0, Z_0) - H_r(C_i(Z_0), Z_0)|.
  \]
  We assume
  \[
    \Delta h_{r, 0} > 2 C \varepsilon_2.
  \]
  Let $X(\tau) = (P(\tau)$, $Q(\tau)$, $Z(\tau), \tau)$ denote the solution of~\eqref{e:res-passing} with initial data $X_0$.
  Denote $h_r(\tau) = H_r(P(\tau), Q(\tau), Z(\tau))$.
  Then there exists $\tau_1 > \tau_0$ such
  \begin{equation}
    P(\tau_1) = -C_p, \qquad
    \tau_1 - \tau_0 < C (1 + |\ln \Delta h_{r, 0}|) \\
  \end{equation}
  and for $\tau \in [\tau_0, \tau_1]$ we have
  \begin{equation}
    X(\tau) \in \mathcal D, \qquad
    \abs{H_r(X(\tau)) - H_r(X(\tau_0))} < C \varepsilon_2, \qquad
    \norm{Z(\tau) - Z_0} < C \varepsilon_2 (1 + |\ln \Delta h_{r, 0}|).
  \end{equation}
\end{lemma}

\begin{proof}
  We apply Lemma~\ref{l:aux} with $C_\gamma \ge 8\pi$. This lemma gives us the values of the constants (we add $\tilde \cdot$ to the names of this constants)
  \[
    \tilde c_Z, \tilde C, \tilde C_\rho.
  \]
  Fix $c_Z = \tilde c_Z$ and assume $C_\rho \ge \tilde C_\rho$.
  We apply Lemma~\ref{l:aux} with $z_0$ as in the current lemma.
  For small enough $\omega_0$ the conditions on $s_1$ and $z_0$ in Lemma~\ref{l:aux} are satisfied and the domain $\mathcal D$ lies in the domain provided by Lemma~\ref{l:aux}.

  Let us now build a finite set $\mathcal F$ such that elements of this set approximate the functions $F_s$ for any $s$ as in the statement of lemma.
  Denote $V_\Theta = V_F$, where $F = \Theta_3(z)$.
  Let us apply Lemma~\ref{l:near-sep} to $V_\Theta(Q, z_*)$, denote by $\varepsilon_{1, 0}$ the largest value of $\varepsilon_1$ allowed by this lemma.
  Apply Lemma~\ref{l:lim} with $z_0$ in that lemma equal to $z_*$ and with
  $\delta_1 = 0.5 \min(\varepsilon_{1, 0}, \min_{\mathcal Z_B} \Theta_3(z))$.
  This lemma gives $S$, $c$ and a set $\mathcal F_{s_2}$ for each $s_2 \le S$, take
  \[
    \mathcal F = \{\Theta_3\} \cup (\cup_{s_2 \le S} \mathcal F_{s_2}), \qquad
    \mathcal F_* = \{F(Q, z_*): \; F \in \mathcal F\}.
  \]

  Let us now determine the constants used to define the domain $\mathcal D$.
  We will assume $c_z < 0.5 c$ and $\omega_0 c_Z < 0.5 c$, then estimates of Lemma~\ref{l:lim} are valid for $z = z_0 + \hat \omega Z$ with $\norm{Z} < c_Z$.
  Let $\varepsilon_{1, 1}$ be so small that Lemma~\ref{l:near-sep} can be applied to all functions $V_{F_i}, \; F_i \in \mathcal F_*$ with $\varepsilon_1 = \varepsilon_{1, 1}$. Clearly, $\varepsilon_{1, 1} \le \varepsilon_{1, 0}$.
  Let us apply Lemma~\ref{l:near-sep} to the functions $V_{F_i}, \; F_i \in \mathcal F_*$ with $\varepsilon_1 = \varepsilon_{1, 0}$ for $V_\Theta$ and $\varepsilon_1 = \varepsilon_{1, 1}$ for other $V_{F_i}$. We can take the same values of $C_p, C_q$ for all $V$, fix these constants.

  Let us now prove that for small enough $c_z$ and $\omega_0$ for any $z_0$ and $s$ we have
  \begin{equation} \label{e:F-z}
    \norm{V_{F_s}(Q, z_0 + \hat \omega Z) - V_{F_s}(Q, z_*)}_{C^2(Q, Z), Q \in [-C_q, C_q]} < 0.25 \varepsilon_{1, 1},
  \end{equation}
  (here and thereafter we consider all $Z$ with $\norm{Z} < c_Z$ when taking the norm).
  We use that $\norm{F_s}_{C^1(Q, z)}$ and $\norm{\pdv{F_s}{z}}_{C^1(Q, z)}$ are bounded by Lemma~\ref{l:lim}.
  Denote
  \[
    z = z_0 + \hat \omega Z, \qquad
    U(Q, Z) = V_{F_s}(Q, z) - V_{F_s}(Q, z_*).
  \]
  We have (in the formula below $z_{int} \in [z, z_*]$ denotes some intermediate $z$)
  \[
    U = \int_0^Q F_s(\tilde Q, z) - F_s(\tilde Q, z_*) d\tilde Q
    = \int_0^Q \pdv{F_s}{z}\/(\tilde Q, z_{int}(\tilde Q)) d\tilde Q = O(z-z_*) = O(c_z + \hat \omega c_Z).
  \]
  We check in the same way that $\pdv{U}{Q}, \pdv[2]{U}{Q} = O(c_z + \hat \omega c_Z)$.
  We have
  \[
    \pdv{U}{Z} = \hat \omega \pdv{}{z} V_{F_s}(Q, z) =
    \hat \omega \int_0^Q \pdv{F_s}{z}\/(\tilde Q, z) d\tilde Q
    = O(\hat \omega).
  \]
  We check in the same way that $\pdv{U}{Q}{Z}, \pdv[2]{U}{Z} = O(\hat \omega)$, this proves~\eqref{e:F-z}.

  Decompose $Q_0 = Q_{0, 1} + \Delta Q$ with $\Delta Q = 2\pi k$, $Q_{0, 1} \in [-\pi, \pi]$.
  To apply Lemma~\ref{l:near-sep}, we use shifted variable $\tilde Q$ defined by
  $Q = \tilde Q + \Delta Q$.
  Set $\tilde V_F(\tilde Q) = V_F(Q)|_{Q = \tilde Q}$.
  Given $F_s$ from~\eqref{e:res-passing}, we take $\tilde V_{F_i}(\tilde Q, z_*)$ as $V(\tilde Q)$ in Lemma~\ref{l:near-sep}, where $F_i$ is the closest to $F_s$ element of $\mathcal F$ and set in Lemma~\ref{l:near-sep}
  \begin{align}
  \begin{split}
    \Delta H(P, \tilde Q, Z)
    &= \tilde V_{F_s}(\tilde Q, z_0 + \hat \omega Z)
    - \tilde V_{F_i}(\tilde Q, z_*)
    + \hat \omega^{-2}s_1 \beta H_7(P, \tilde Q + \Delta Q, Z), \\
    (v_p, v_q, v_z) &= (\alpha u_P, \alpha \hat \omega^{-0.5} u_Q, \alpha \hat \omega^{-0.5} u_z)|_{Q = \tilde Q+\Delta Q}.
  \end{split}
  \end{align}
  We have
  \begin{align}
  \begin{split}
    \norm{\Delta H}_{C^2(\tilde Q, Z), \tilde Q \in [-C_q, C_q]}
    &\le
    \norm{V_{F_s}(Q, z) - V_{F_s}(Q, z_*)}_{C^2(Q, Z), Q \in [-C_q, C_q]} \\
    &+
    \norm{V_{F_i}(Q, z_*) - V_{F_s}(Q, z_*)}_{C^2(Q), Q \in [-C_q, C_q]} \\
    &+
    \norm{\hat \omega^{-2}s_1 \beta H_7}_{C^2(P, Q, Z), Q \in [\Delta Q-C_q, \Delta Q+C_q]}.
  \end{split}
  \end{align}
  The first term is bounded by $0.25 \varepsilon_{1, 1}$ due to~\eqref{e:F-z}.
  The second term is bounded by $0.5 \varepsilon_{1, 0}$ for $V_{F_i} = V_\Theta$ due to item~3 of Lemma~\ref{l:lim}, and by $0.5 \varepsilon_{1, 1}$  for other $V_{F_i}$ when $\omega_0$ is small enough due to item~4 of the same lemma.
  Let us now show that for large enough $C_\rho$ the third term is bounded by $0.25 \varepsilon_{1, 1}$.
  Indeed, by~\eqref{e:hat-h-bound} we have $|\beta(z_0+\hat \omega Z)| \sim \beta(z_0)$ for complex $Z$ with $\norm{Z} < 8 c_Z$, by Cauchy formula this means $\pdv{\beta}{Z}, \pdv[2]{\beta}{Z} = O(\beta(z_0))$. Together with the bound on $\norm{H_7}$ from Lemma~\ref{l:aux} this implies
  \[
    \norm{\hat \omega^{-2}s_1 \beta H_7}_{C^2} = O(\hat \omega^{-2} s_1 \beta(z_0)) =
    O(\sqrt{\varepsilon \hat h^{-1} \hat \omega^{-1} s_1^2})
    = O(\sqrt{\varepsilon \hat h^{-1} C_{s_1}^2 |\ln^5 \varepsilon|}) =
    O(\sqrt{C_{s_1}^2 / C_\rho}).
  \]

  Thus, $\norm{\Delta H}_{C^2(Q, Z)}$ is bounded by $\varepsilon_{1, 0}$ if $V_{F_i} = V_\Theta$ and by $\varepsilon_{1, 1}$ otherwise.
  This means that we can apply Lemma~\ref{l:near-sep} with these $\Delta H$ and $v$.
  Recall that we took $\varepsilon_2 = C \alpha(z_0) \hat \omega^{-1}$. Clearly, for large enough $C$ we have $\norm{v}_{C^0} < \varepsilon_2$. We have
  \[
    \varepsilon_2 = O(\sqrt{\varepsilon \hat h(z_0) |\ln^5 \hat h|}).
  \]
  It is easy to see that $\varepsilon_2$ satisfies the conditions in Lemma~\ref{l:near-sep} for small enough $\varepsilon_0$ and large enough $C_\rho$.
  We can finally apply Lemma~\ref{l:near-sep}, we assume that $C$ from the current lemma is greater than the constant $C$ from Lemma~\ref{l:near-sep}. Then Lemma~\ref{l:near-sep} gives us the estimates stated in the current lemma.
\end{proof}

\subsection{Low-numerator resonances: proof} \label{ss:low-num-res}
In this subsection we prove Lemma~\ref{l:est-res-cross}, thus estimating the measure of initial conditions captured into resonances and time spent in resonance zones for initial conditions that are not captured. To prove this lemma we first prove auxiliary Lemma~\ref{l:strong-res}.
Recall that $m$ denotes the Lebesgue measure on $\mathcal A$ and $B_r(O)$ denotes the ball with center $O$ and radius $r$. Denote by $m_z$ the Lebesgue measure on $\mathbb R^n$.
\begin{lemma} \label{l:strong-res}
  Given $z_* \in \mathcal Z_B$, for any
  $C_{s_1}, D_{p, 0} > 0$
  there exist
  \[
    c_z, \omega_0, \varepsilon_0 > 0, \; D_p > D_{p, 0}, \; C, C_\rho, C_B > 1
  \]
  such that
  for any
  $\hat \omega = s_2 / s_1 \in (0, \omega_0)$,
  and $\varepsilon < \varepsilon_0$ with
  $|s_1| < C_{s_1} \ln^2 \varepsilon$
  there exist a finite collection of balls $B_i$ with centers
  $z_i \in B_{c_z}(z_*)$
  and equal radii $\hat \omega c_Z$ such that
  \begin{equation} \label{e:B-i}
    \cup_{i} B_{\hat \omega c_Z / 2}(z_i) \supset B_{c_z}(z_*),
    \qquad \sum_i m_z(B_i) \le C_B
  \end{equation}
  and a collection of sets
  \[
    \mathcal V_i \subset \mathbb R^2_{p, q} \times B_i \times [0, 2\pi]_\lambda,
  \]
  such that for any $i$ we have
  \[
    m(\mathcal V_i) \le C \alpha_i^2 s_1 \hat \omega^{-0.5} m_z(B_i)
  \]
  (we denote $\alpha_i = \alpha(z_i)$).

  These collections satisfy the following:
  for any $i$ for any initial data
  $p_0, q_0, z_0, \lambda_0$ with
  \begin{equation} \label{e:init-in-B-i}
    I(p_0, q_0, z_0) = \hat I(z_0) + D_p \alpha(z_0) \hat \omega^{0.5}, \qquad
    z_0 \in B_{\hat \omega c_Z/2}(z_i), \qquad
    h(p_0, q_0, z_0) > C_\rho \varepsilon |\ln \varepsilon^5|
  \end{equation}
  at least one of the following holds.
  \begin{itemize}
    \item The solution of the perturbed system~\eqref{e:perturbed-pq} with this initial data crosses the hypersurface $I = \hat I(z) - D_p \alpha(z) \hat \omega^{0.5}$ at some time $\lambda_1 > \lambda_0$ with
    \[
      \varepsilon(\lambda_1 - \lambda_0) \le  C |\ln \varepsilon| \hat \omega^{0.5} \alpha_i.
    \]
    \item There exists $\lambda_1 > \lambda_0$ with
    \[
      \varepsilon(\lambda_1 - \lambda_0) \ge C^{-1} \hat \omega^{0.5} \alpha_i
    \]
    such that this solution remains in $\mathcal V_i$ for $\lambda \in [\lambda_0, \lambda_1]$.
  \end{itemize}
\end{lemma}
\begin{proof}
  Recall that we define $D(I, z)$ by the equality
  \[
    I = \hat I(z) + D \alpha(z) \hat \omega^{0.5}.
  \]
  The values of $D_{p, 1}$ in Lemma~\ref{l:near-sep-2} should be large enough, we denote by $D'_{p, 0}$ a constant such that this lemma can be applied for $D_{p, 1} > D'_{p, 0}$.
  The coordinate change of Lemma~\ref{l:aux} (we do not specify the details yet, this will be done later) gives the variables $P, Q$.
  By~\eqref{e:p-q-J-gamma} we have $D = P + O(s_1 \hat \omega^{-1.5} \beta)$. We have
  \[
    s_1 \hat \omega^{-1.5} \beta
    \lesssim \sqrt{\varepsilon h^{-1} \ln^4 \varepsilon} \lesssim |\ln \varepsilon|^{-0.5} \ll 1.
  \]
  Thus we can take $D_p > D_{p, 0}$ and $D_{p, 2}, D_{p, 1} > D'_{p, 0}$ with $D_{p, 2} > D_{p, 1} + 3$ such that
  \begin{itemize}
    \item $P(I, z, \varphi, \lambda) \in [D_{p, 1} + 1, D_{p, 2} - 1]$ for any $I, z, \varphi, \lambda$ with $D(I, z)=D_p$
    \item $P(I, z, \varphi, \lambda) = -D_{p, 2}$ implies $D(I, z) < -D_p$.
  \end{itemize}

  Let us now apply Lemma~\ref{l:near-sep-2} with $D_{p, 1}$ and $D_{p, 2}$ chosen above, it provides us with constants (we add $\tilde \cdot$ to constants provided by this lemma)
  \[
    \tilde C_p, \tilde C_\rho, \tilde C, \tilde C_q, \tilde c_z, \tilde c_Z, \tilde \omega_0, \tilde \varepsilon_0.
  \]
  Set $c_z = \tilde c_z$, $c_Z = \tilde c_Z$, $C_\rho = 2 \tilde C_\rho$.

  The existence of $C_B$ and $\{B_i\}$ is obvious, we can place the centers $z_i$ in the nodes of a hyper-cubic lattice (with step $\sim \hat \omega c_Z$).
  Let us now fix some $i$. Assume that condition~\eqref{e:init-in-B-i} is satisfied for some initial condition, else we have nothing to prove.
  Let us apply Lemma~\ref{l:near-sep-2} with $z_0$ (from that lemma) equal to $z_i$.
  The bound on $\hat h(z_i)$ required by Lemma~\ref{l:near-sep-2} is satisfied by~\eqref{e:hat-h-bound},~\eqref{e:init-in-B-i} and our choice of $C_\rho$.

  We will use the coordinates $P, Q$ provided by Lemma~\ref{l:aux} applied as a part of statement of Lemma~\ref{l:near-sep-2}.
  Set
  \[
    \varepsilon_{2, i} = 2 \tilde C \alpha_i \hat \omega^{-1},
  \]
  it is greater that the value of $\varepsilon_2$ in Lemma~\ref{l:near-sep-2} ($\tilde C \alpha(z_0) \hat \omega^{-1}$).
  We have $\ln \varepsilon_{2, i} \sim \ln \varepsilon$. Set
  \[
    \mathcal {\tilde V}_{0, i} = \Big\{ (Z, P, Q, \tau) \in \mathcal D_i : \norm{Z} < c_Z, \;
    P \in [D_{p, 1}, D_{p, 2}], \; Q \in [-\pi \hat \omega^{-1} - \tilde C_q, \pi \hat \omega^{-1} + \tilde C_q]  \Big\},
  \]
  where $\mathcal D_i$ denotes the domain $\mathcal D$ from Lemma~\ref{l:near-sep-2}.
  Define $\mathcal{\tilde V}_i \subset \mathcal{\tilde V}_{0, i}$ by additional condition that the value of the Hamiltonian $H_{r, i}$ defined by~\eqref{e:H-without-u} is at least $2 \tilde C \varepsilon_{2, i}$-far from its values at the saddles of the Hamiltonian system given by $H_{r, i}$.
  Let
  $\mathcal{V}_i$
  be the preimage of
  $\mathcal{\tilde V}_i$
  under the coordinate change of Lemma~\ref{l:aux}.

  Take an initial condition satisfying~\eqref{e:init-in-B-i}.
  Rewriting our initial condition in the new chart gives $P_0, Q_0, Z_0, \tau_0$ (there is some freedom, as we can add $2\pi k$ to $\varphi_0 = \varphi(p_0, q_0)$, this will be resolved later). By the choice of $D_p$ we have $P_0 \in [D_{p, 1}+1, D_{p, 2}-1]$. We can also achieve
  $Q_0 \in [-\pi \hat \omega^{-1} - 1, \pi \hat \omega^{-1} + 1]$. Indeed, by~\eqref{e:p-q-J-gamma} it is sufficient to have $\gamma(\varphi_0, \lambda_0) \in [-\pi, \pi]$
  and this can be accomplished by adding $2\pi k$ to $\varphi_0$.

  We study the evolution of this initial condition using the system~\eqref{e:res-passing}.
  As $P', Q', Z' = O(1)$ in~\eqref{e:res-passing}, solutions starting in $\mathcal{\tilde V}_{0, i}$ with $\norm{Z_0} < 0.5 c_Z$ and $P \in [D_{p, 1} + 1, D_{p, 2} - 1]$ spend in $\mathcal{\tilde V}_{0, i}$ time $\tau \gtrsim 1$, thus time $\lambda \gtrsim \varepsilon^{-1} \alpha \hat \omega^{0.5}$.
  Indeed, by~\eqref{e:tau-Z} we have
  \[
    \Delta \lambda = \varepsilon^{-1} \alpha \hat \omega^{0.5} \Delta \tau.
  \]
  If during the crossing of $\mathcal{\tilde V}_{0, i}$ the solution is not in $\mathcal{\tilde V}_i$ (at just one time), by Lemma~\ref{l:near-sep-2} this solution reaches $P = - \tilde C_p < -D_{p, 2}$ at some point.
  By the choice of $D_p$ this implies $D \le -D_p$, by continuity at some point before this solution crosses $D = -D_p$, i.e. $I = \hat I(z) - D_p \alpha(z) \hat \omega^{0.5}$.
  The estimate on the time $\lambda$ before this moment follows from the estimate $\tau \lesssim |\ln \varepsilon_2| \lesssim |\ln \varepsilon|$ on the time $\tau$ before this moment provided by Lemma~\ref{l:near-sep-2}.

  Let us consider sections $\tilde {\mathcal V}_{i, z}$ of $\tilde {\mathcal V}_i$ and sections ${\mathcal V}_{i, z}$ of ${\mathcal V}_i$ with fixed $z$.
  We will also consider sections $\tilde {\mathcal V}_{i, z, \tau}$ of $\tilde {\mathcal V}_i$ with fixed $z$ and $\tau$.
  In $\mathcal{\tilde V}_0$ we have $\pdv{H_r}{P} \approx P > D_{p, 1}$ and thus $\pdv{H_r}{P} > 0.5 D_{p, 1}$. This means that the measures of intersections of $\mathcal{\tilde V}_{i, z, \tau}$ with each segment where only $P$ varies are $O(\varepsilon_{2, i}) n_{sadd}$, where $n_{sadd}$ is the number of saddles we need to consider.
  When $Q$ is fixed, $n_{sadd} = O(1)$.
  Indeed, for fixed $Q$ the range of values of $H_r$ is bounded and so, as (we use notation $V_c$ from~\eqref{e:V-c}) $V_c > 0$, the number of saddles where the value of $H_r$ is close to that range is also bounded.
  Thus, as the range of values of $Q$ is
  $\sim \hat \omega^{-1}$,
  the measures of $\mathcal{\tilde V}_{i, z, \tau} \subset \mathbb R^2_{P, Q}$ are
  $O(\varepsilon_{2, i} \hat \omega^{-1}) = O(\alpha_i \hat \omega^{-2})$.

  To estimate the measure of the set $\mathcal V_{i, z}$, we follow the construction of the coordinate change of Lemma~\ref{l:aux} (see the proof of this lemma) in the reverse order.
  The $z$-dependent time change $\mu \to \tau$ does not change the measures of $\mathcal{\tilde V}_{i, z, \tau}$ (for the preimages of these sets $\mu$ is fixed instead of $\tau$). As our system is $2\pi$-periodic in $\mu$, integrating by $\mu$ shows that measures of the preimages of $\mathcal{\tilde V}_{i, z}$ in
  $\mathbb R^2_{P, Q} \times \mathbb [0, 2\pi]_\mu$
  are $O(\alpha_i \hat \omega^{-2})$.
  Before the scaling $\tilde P, \tilde Q \to P, Q$
  the measures of the preimages of $\mathcal{\tilde V}_{i, z}$ in $\mathbb R^2_{\tilde P, \tilde Q} \times [0, 2\pi]_\mu$ are $O(\alpha_i \hat \omega^{-0.5})$.
  The coordinate change of Corollary~\ref{c:averaging} is volume-preserving, so the measure of the preimages of $\mathcal{\tilde V}_{i, z}$ in
  $\mathbb R^2_{J, \gamma} \times [0, 2\pi]_\mu$
  is again $O(\alpha_i \hat \omega^{-0.5})$.
  We project $\gamma$ from $\mathbb R$ to $[0, 2\pi]$, this does not increase the measures.
  Before the change of angular variables the measures in
  $\mathbb R_J \times [0, 2\pi]^2_{\varphi, \lambda}$
  are $O(\alpha_i s_1 \hat \omega^{-0.5})$.
  Before the scaling of $I$ the measures in
  $\mathbb R_I \times [0, 2\pi]^2_{\varphi, \lambda}$
  are $O(\alpha_i^2 s_1 \hat \omega^{-0.5})$.

  Finally, moving from $I, \varphi$-variables back to $p, q$ (this preserves the volume), we get that the measure of the preimages of
  $\mathcal{\tilde V}_{i, z}$
  in
  $\mathbb R^2_{p, q} \times [0, 2\pi]_\lambda$
  (i.e. $\mathcal V_{i, z}$)
  is $O(\alpha_i^2 s_1 \hat \omega^{-0.5})$.
  Taking union of these preimages for $z \in B_i$ gives the estimate on the measure of $\mathcal V_i$. This completes the proof of the lemma.
\end{proof}

\begin{proof}[Proof of Lemma~\ref{l:est-res-cross}]
  Let us apply Lemma~\ref{l:strong-res} to our $z_*$, $C_{s_1}$, $D_{p, 0}$, it gives us
  \[
    c_z, \omega_0, \varepsilon_0, \; D_p > D_{p, 0}.
  \]
  Denote by $M_i$ the set of all initial conditions $X_{init}$ such that the second alternative of Lemma~\ref{l:strong-res} holds with this $i$ for some $X(\lambda_0)$ with $\lambda_0$ as above, let $m_i$ be the measure of $M_i$. Set $\mathcal E_s = \cup_i M_i$. Then the statement of the lemma holds if $X_{init} \not \in \mathcal E_s$, it only remains to estimate the measure of $\mathcal E_s$.

  Denote by $g^\lambda$ the flow given by the perturbed system.
  Let $\lambda_i \sim \varepsilon^{-1} \alpha_i \hat \omega^{0.5}$ denote the estimate from below on the time spent in $\mathcal V_i$ claimed in the second alternative. Set
  \[
    K_1 = \lceil \varepsilon^{-1} \Lambda/\lambda_i \rceil
    \sim \alpha_i^{-1} \hat \omega^{-0.5}.
  \]
  Then we have
  \[
    M_i \subset \bigcup_{j = 0}^{K_1} g^{-j \lambda_i}(\mathcal V_i).
  \]
  By Lemma~\ref{l:volume} there is $K_2 > 1$ such that for any $\lambda \in [0, \varepsilon^{-1}\Lambda]$ the map $g^{-\lambda}$ expands the volume $m$ by factor at most $K_2$.
  Thus we have the estimate
  \[
    m_i \le K_2 K_1 m(\mathcal V_i) \lesssim
    \alpha_i s_1 \hat \omega^{-1} m_z(B_i).
  \]
  As $\sum_i m_z(B_i) = O(1)$, this gives
  \[
    m(\mathcal E_s) \le \sum_i m_i \lesssim s_1 \hat \omega^{-1} \max_i {\alpha_i} \lesssim s_1 \hat \omega^{-1} \alpha_*.
  \]
\end{proof}

\section{Passing separatrices: proof} \label{s:pass-sep}
\begin{proof}[Proof of Lemma~\ref{l:near-sep-local}]
  Pick $R_Z \sim |\ln^{-1} h_*| \sim |\ln^{-1} \varepsilon|$ such that we have $R_Z < d_Z \omega(h, z)$ for all $z \in B_{c_z}(z_*)$ and $h \in [h_*, 4h_*]$, where the constant $d_Z$ is provided by Lemma~\ref{l:d-z}.
  Build a finite collection of balls $B_i$ with centers
  $z_i \in B_{c_z}(z_*)$
  and equal radii $R_Z$ such that for some $C_B > 0$ we have
  \begin{equation}
    \cup_{i} B_{R_Z / 2}(z_i) \supset B_{c_z}(z_*),
    \qquad \sum_i m_z(B_i) \le C_B
  \end{equation}
  (recall that $m_z$ denotes Lebesgue measure on $\mathbb R^n$).
  Fix some $i$ and set
  \[
    \mathcal L_i = \big\{
      X \in B_{R_Z}(z_i) \times \mathbb R^2_{p, q} \times [0, 2\pi]_{\lambda}
      \; :  h(X) \in [-h_*, 4h_*]
    \big\}.
  \]
  Let us prove the following alternative: for any $X_0$ with $z_0 \in B_{R_Z / 2}(z_i)$ and $h(X_0) = h_*$ (then $X_0 \in \mathcal L_i$)
  the corresponding solution of~\eqref{e:perturbed-pq} either reaches $h = -h_*$ after time at most $\Delta \lambda$ or stays in $\mathcal L_i$ for time $\Delta \lambda$, where
  \[
    \Delta \lambda = \varepsilon^{-1} \sqrt{\varepsilon} |\ln^\gamma \varepsilon|.
  \]

  To prove this alternative, we need to show that this solution cannot leave $\mathcal L_i$ through $|z-z_i| = R_Z$ or through $h = 4h_*$.
  The first part follows from $\dot z = O(\varepsilon)$, thus changing $z$ by $0.5 R_Z$ requires time at least
  $\sim \varepsilon^{-1} / R_Z \gg \Delta \lambda$.

  To prove the second part, we find a non-resonant zone $Z_{r, r+1}$ such that we have $(h, z_i) \in Z_{r, r+1}$ for some $h \in [1.4 h_*, 1.6 h_*]$.
  We can find such non-resonant zone, as different resonant zones with $h \sim h_*$ do not intersect and have width $\sim \sqrt{\varepsilon h_*} \ll h_*$ in $h$.
  By Lemma~\ref{l:d-z} (applied to $s_r$ and $s_{r+1}$) for all $z \in B_{R_Z}(z_i)$ there exists some $h \in [1.2 h_*, 1.8 h_*]$ with $(h, z) \in Z_{r, r+1}$.

  By Lemma~\ref{l:nonres} if at some time we have $X(\lambda) \in Z_{r, r+1}$ (this must happen before reaching $h=4h_*$), we then reach either $\partial \Pi$ (then $h \ll h_*$) or the border of $Z_{r, r+1}$ with $Z_{r+1}$ (then $h$ is again less then in $Z_{r, r+1}$ and the solution must again enter $Z_{r, r+1}$ before reaching $h=4h_*$).
  The third part of Lemma~\ref{l:nonres} gives the estimate $h \le \frac 5 3 \times 1.8 h_* = 3h_*$, so it is impossible to have $h > 3 h_*$ while $z$
  stays in $B_{c_z}(z_*)$.
  The alternative is proved.

  Denote by $g^\lambda$ the flow given by the perturbed system.
  Set
  \[
    K_1 = \lceil \varepsilon^{-1} \Lambda/\Delta \lambda \rceil
    \sim \varepsilon^{-0.5} \ln^{-\gamma} \varepsilon
  \]
  and
  \[
    \mathcal C_i = \bigcup_{j = 0}^{K_1} g^{-j \Delta \lambda}(\mathcal L_i).
  \]
  This set covers all $X_{init}$ such that second part of the alternative can be satisfied for some $\lambda_0$ with $z(X(\lambda_0)) \in B_{R_Z/2}(z_i)$. Indeed, if the second part of the alternative is satisfied, we have $X(\lambda) \in \mathcal L_i$ for all $\lambda \in [\lambda_0, \lambda_0 + \Delta \lambda]$, thus for some natural $k \ge 0$ we have $X(\lambda_{init} + k \Delta \lambda) \in \mathcal L_i$ and $X_{init} \in g^{-k \Delta \lambda}(\mathcal L_i)$.

  By Lemma~\ref{l:volume} there is $K_2 > 1$ such that for any $\lambda \in [0, \varepsilon^{-1}\Lambda]$ the map $g^{-\lambda}$ expands the volume $m$ by factor at most $K_2$.
  As $\Delta h = h_*$ corresponds to $\Delta I \sim h_* \ln \varepsilon$, we have $m(\mathcal L_i)  \sim h_* \ln \varepsilon m_z(B_i)$.
  Thus we have the estimate
  \[
    m(\mathcal C_i) \le K_2 K_1 m(\mathcal L_i) \lesssim
     m_z(B_i) h_* \ln \varepsilon \; \varepsilon^{-0.5} \ln^{-\gamma} \varepsilon
     \lesssim \sqrt{\varepsilon} \ln^{\rho - \gamma + 1} \varepsilon \; m_z(B_i).
  \]
  Set $\mathcal E = \cup \mathcal B_i$.
  As $\sum_i m_z(B_i) = O(1)$, this gives
  \[
    m(\mathcal E) \le \sum_i m(\mathcal B_i)
    \lesssim \sqrt{\varepsilon} \ln^{\rho - \gamma + 1} \varepsilon.
  \]
  Clearly, if $X_0 \not \in \mathcal E_s$, we have $z_0 \in B_{R_Z/2}(z_i), X_0 \not \in \mathcal E_i$ for some $i$. But then the corresponding solution cannot stay in $\mathcal L_i$ for time $\Delta \lambda$, thus it crosses $h = -h_*$.
\end{proof}

\newpage
\section{Probabilities} \label{s:prob}
Denote
\begin{equation}
  \mathcal F(h, z, i) = \frac 1 {2\pi T} \int_0^{2\pi} \oint \Div f dt d\lambda.
\end{equation}
The inner integral above is taken along the closed trajectory of the unperturbed system
given by $h = h, z = z$ that lies inside $\mathcal B_i$, this trajectory is parametrized by the time $t$.
\begin{lemma} \label{l:int-div}
  In the assumptions of Theorem~\ref{t:main},
  there exists a set $\mathcal E_1$ with
  \[
    m(\mathcal E_1) = O(\sqrt{\varepsilon} |\ln^5 \varepsilon|)
  \]
  such that if the initial condition is not in $\mathcal E_1$, for $\lambda \in [\lambda_0, \lambda_0 + \varepsilon^{-1} \Lambda]$ we have
  \begin{equation}
    \varepsilon \int_{\lambda_0}^{\lambda} \Div f(X(\lambda)) d\lambda
    = \varepsilon \int_{\lambda_0}^{\lambda} \mathcal F(\overline X(\lambda)) d\lambda + O(\sqrt{\varepsilon} |\ln \varepsilon|).
  \end{equation}
\end{lemma}
\begin{proof}
  Set $\tilde z = (z, z_\rho)$, $f_{\tilde z} = (f_z, \Div f)$.
  Consider extended perturbed system, where we replace $z$ and $f_z$ by $\tilde z$ and $f_{\tilde z}$.
  Then along solutions of extended perturbed system we have
  $\dot z_\rho = \varepsilon \Div f$
  and along solutions of extended averaged system we have
  $\dot{\overline z}_\rho = \mathcal F$.
  Note that the right-hand side of the extended systems (perturbed and averaged) does not depend on $z_\rho$, the evolution of all variables except $z_\rho$ is given by the ititial systems.
  Apply Theorem~\ref{t:main} to the extended system and denote by $\mathcal E_1$ the excluded set for extended system. Solutions of perturbed and averaged system are $O(\sqrt{\varepsilon} |\ln \varepsilon|)$-close to each other. Taking $z_\rho$-coordinates yields
  \begin{equation}
    \varepsilon \int_{\lambda_0}^{\lambda} \Div f(X(\lambda)) d\lambda
    = z_\rho(\lambda) - z_\rho(\lambda_0)
    = \overline z_\rho(\lambda) - \overline z_\rho(\lambda_0)
    + O(\sqrt{\varepsilon} |\ln \varepsilon|)
    = \varepsilon \int_{\lambda_0}^{\lambda} \mathcal F(\overline X(\lambda)) d\lambda
    + O(\sqrt{\varepsilon} |\ln \varepsilon|).
  \end{equation}
\end{proof}

Given a vector $v_0 = (v_0^1, \dots, v_0^s)$, denote by $C_\delta(v_0)$ the cube with side $2\delta$ and center $v_0$. Set
\begin{equation}
  W^\delta(I_0, z_0) = C_\delta(I_0, z_0) \times [0, 2\pi]^2_{\varphi, \lambda} \subset \mathcal A_3.
\end{equation}
Let $W_i^\delta \subset W^\delta \cap \mathcal E$ be the subset of initial conditions such that the solution is captured in $\mathcal A_i, \; i=1,2$.
\begin{lemma} \label{l:prob-W}
  \begin{equation}
    m(W^\delta_i) = \int_{W^\delta} \frac{\Theta_i(z_*)}{\Theta_3(z_*)} dp_0 dq_0 dz_0
    + O(\sqrt{\varepsilon} |\ln^5 \varepsilon|),
    \qquad i=1,2.
  \end{equation}
  Here $z_*(p_0, q_0, z_0)$ denotes the value of $z$ when the solution of the averaged system with initial data $(h(p_0, q_0, z_0), z_0)$ crosses the separatrices.
\end{lemma}
\begin{proof} This lemma is similar to~\cite[Proposition~2.4]{neishtadt17} and can be proved in the same way. In~\cite{neishtadt17} the proof of Proposition~2.4 is based on Lemma~5.1 and Lemma~5.2.
Lemma~5.1 in~\cite{neishtadt17} is an analogue of Lemma~\ref{l:int-div} above, and we replace the use of Lemma~5.1 in~\cite{neishtadt17} by our Lemma~\ref{l:int-div}. Lemma~5.2 in~\cite{neishtadt17} is a statement about the averaged system (in one-frequency case). One can obtain averaged system in two-frequency case as follows: first, average the perturbation over the time $\lambda$ and obtain a system with one frequency (corresponding to the angle $\varphi$) and then take average over $\varphi$. Thus, averaged system in two-frequency case is the same as averaged system in one-frequency system obtained after averaging over $\lambda$, and the statement of Lemma~5.2 holds for our case. The rest of the proof of Proposition~2.4 can be straightforwardly applied to our case, we just need to change the error terms to adjust for the difference of precision of averaged method and the measure of exceptional set in one-frequency and two-frequency cases.
\end{proof}

Recall the notation
\begin{equation}
  U^\delta(I_0, z_0, \varphi_0, \lambda_0) = C_\delta(I_0, z_0, \varphi_0, \lambda_0).
\end{equation}
\begin{lemma} \label{l:prob-W-to-U}
  We have the following for $i=1,2$.
  Suppose that for $(I_0, z_0)$ in some open set $\mathcal U$ we have
  \begin{equation} \label{e:lim-W}
    \lim_{\delta \to 0} \lim_{\varepsilon \to 0} \frac{m(W^\delta_i)}{m(W^\delta)} = \psi(I_0, z_0),
  \end{equation}
  where the function $\psi$ is continuous in $\mathcal U$ and the limit is uniform for all $(I_0, z_0) \in \mathcal U$.
  Then for any $(I_0, z_0) \in \mathcal U$ and any $\varphi_0, \lambda_0 \in [0, 2\pi]$ we have
  \begin{equation} \label{e:lim-U}
    \lim_{\delta \to 0} \lim_{\varepsilon \to 0} \frac{m(U_i^\delta)}{m(U^\delta)} = \psi(I_0, z_0).
  \end{equation}
\end{lemma}
\begin{proof}
  Assume w.l.o.g. that $i=1$.
  Suppose we are given some tolerance $\kappa$ and we want to show that for each small enough $\delta$ for small enough $\varepsilon$ we have
  $|m(U^\delta_1) / m(U^\delta) - \psi| \le \kappa$ for fixed $I_0, z_0$ and all $\varphi_0, \lambda_0$.
  The value of $\delta$ should be so small that
  in the set $\mathcal U_0 = C_\delta(I_0, z_0)$
  the values of $\psi$ differ from $\psi(I_0, z_0)$ by at most $\kappa/9$.

  Denote $\phi = (\varphi, \lambda)$.
  Let $m_\phi$, $m_\omega$ and $m_{I, z}$ denote the Lebesgue measure on the space where the corresponding variable(s) is defined: $\mathbb R^2$,  $\mathbb R$ and $\mathbb R^{n+1}$, respectively.
  Set $\Omega_0 = \{\omega(I, z), \; (I, z) \in \mathcal U_0 \}$.
  Denote by $\Omega(\nu, \delta_\phi, N) \subset \Omega_0$ the set of all $\omega$ such that the flow $R_\omega^t(\varphi, \lambda) = (\varphi+\omega t, \lambda + t)$ on $[0, 2\pi^2]$ satisfies the following: for any $\phi_1, \phi_2$ there exists $t \in [0, N]$ such that
  \begin{equation} \label{e:triag}
    m_\phi\Big( C_{\delta_\phi}(\phi_1) \triangle R^t_\omega(C_{\delta_\phi}(\phi_2)) \Big) \le \nu \delta_\phi^2.
  \end{equation}
  We have the following
  \begin{itemize}
    \item $\Omega(\nu, \delta_\phi, N)$ is compact.
    \item $\lim_{k \to \infty} m_\omega(\Omega(\nu, \delta_\phi, k))  = m_\omega(\Omega_0)$ for each $\nu > 0$.
  \end{itemize}
  The first property is straightforward (we can write $\Omega(\kappa, \delta_\phi, N)$ as intersection of sets such that the property above is satisfied for each pair $(\phi_1, \phi_2)$, each such set is compact as we can pick convergent subsequence of the values of $t$). To prove the second property, note that if $\omega$ is irrational, we have $\omega \in \Omega(\nu, \delta_\phi, k)$ for large enough $k$, as orbits of the flow $R_\omega^t$ are dense. This gives
  $m_\omega(\cup_k \Omega(\nu, \delta_\phi, k)) = m_\omega(\Omega_0)$, our statement follows from the continuity of Lebesgue measure.

  Now set $\delta_\phi=\delta$, $\nu=\kappa/100$ and
  $\Omega_k = \Omega(\nu, \delta, k)$.
  Take large enough $k$ so that
  \begin{equation} \label{e:tmp-7207}
    m(\{X: \omega(X) \in \Omega_0 \setminus \Omega_k\}) \le \kappa/9 \; m(W^\delta)
  \end{equation}
  (this can be done, as $\pdv{\omega}{X}$ is non-degenerate).
  Each point of $\Omega_k$ has a neighborhood (a segment with center at this point) such that~\eqref{e:triag} holds for all $\omega$ in this segment with the same $t$ but $2\nu$ instead of $\nu$:
  \begin{equation} \label{e:triag-relaxed}
    m_\phi\Big( C_\delta(\phi_1) \triangle R_\omega^t(C_\delta(\phi_2)) \Big) \le 2\nu \delta_\phi^2.
  \end{equation}
  Pick a finite cover $\alpha$ of $\Omega_k$ by such segments. Let $\{I_j\}$ denote all segments between endpoints of the segments from $\alpha$ that intersect $\Omega_k$, the segments $I_j$ are disjoint from each other and their union covers $\Omega_k$. Denote $K_j = \{ (I, z) \in \mathcal U_0 : \omega(I, z) \in I_j \}$.
  We can decompose this set (except a subset of small measure) as disjoint union of cubes $K_{j, l}$ with small side $\delta_w$.
  We pick $\delta_w$ so small that most of volume of $K_j$ is covered (except the proportion at most $\kappa/9$) and
  $|m(W^\delta_1) / m(W^\delta) - \psi| \le \kappa/9$ when $\delta \le \delta_w$ for small enough $\varepsilon$.
  Let us now apply~\eqref{e:lim-W} to the set $W^\delta = W^{j, l} = K_{j, l} \times [0, 2\pi]^2_\phi$.
  We obtain for small enough $\varepsilon$ (in the formula below and thereafter the lower index $1$ denotes that we consider initial data from some set that is captured in $\mathcal A_1$)
  \[
    \frac{m(W_1^{j, l})}{m(W^{j, l})} \in [\psi(I, z) - \kappa/9, \psi(I, z) + \kappa/9]
    \subset [\psi(I_0, z_0) - 2\kappa/9, \psi(I_0, z_0) + 2\kappa/9].
  \]
  Denote by $W^{j, l, \phi_0}$ the subset of $W^{j, l}$ defined by the condition $|\varphi-\varphi_0|, |\lambda-\lambda_0| < \delta_\phi$, where $(\varphi_0, \lambda_0) = \phi_0$.
  We have~\eqref{e:triag-relaxed}. Thus for small enough $\varepsilon$ we have
  \begin{equation}
    m\Big( W^{j, l,\phi_1} \triangle g_\varepsilon^t(W^{j, l,\phi_2}) \Big) \le 3\nu \delta^2 m_{I, z}(W^{j, l})
  \end{equation}
  for any $\phi_1$, $\phi_2$.
  Here $g_\varepsilon^t$ denotes the flow of the perturbed system.
  The set of points captured in $\mathcal A_1$ is invariant under $g_\varepsilon^t$, so the estimate above means that $m(W_1^{j, l, \phi})$ is almost the same for all $\phi$, with difference $\le 3\nu \delta_\phi^2 m_{I, z}(W^{j, l})$. As the average of $m(W_1^{j, l, \phi}) / m(W^{j, l, \phi})$ over $\phi$ is $m(W_1^{j, l}) / m(W^{j, l})$, this means (given $\nu = \kappa/100$) for all $\phi$
  \begin{equation}
    \frac{m(W_1^{j, l, \phi})}{m(W^{j, l, \phi})} \in [\psi(I_0, z_0) - 3\kappa/9, \psi(I_0, z_0) + 3\kappa/9].
  \end{equation}
  Denote $K^{j, \phi} = K_j \times C_\delta(\phi)$.
  Taking sum over $l$ gives
  \[
    \frac{m(K_1^{j, \phi})}{m(K^{j, \phi})} \in [\psi(I_0, z_0) - 4\kappa/9, \psi(I_0, z_0) + 4\kappa/9].
  \]
  Finally, taking union over $j$ and using~\eqref{e:tmp-7207}, we get
  \[
    \frac{m(U_1^\delta)}{m(U^\delta)} \in [\psi(I_0, z_0) - 5\kappa/9, \psi(I_0, z_0) + 5\kappa/9].
  \]
  This estimates holds for any $\kappa$ when $\delta$ and $\varepsilon$ are small enough (and $\varepsilon$ is small compared with $\delta$). This completes the proof.
\end{proof}

\begin{proof}[Proof of Proposition~\ref{p:prob}]
  Fix any open set $\mathcal U \subset \mathcal A_3$ that is separated from the separatrices.
  By Lemma~\ref{l:prob-W} we have
  \begin{equation}
    \lim_{\delta \to 0} \lim_{\varepsilon \to 0} \frac{m(W^\delta_i)}{m(W^\delta)} = \Theta_i(z_*)/\Theta_3(z_*),
  \end{equation}
  where $z_*$ is taken at $(I_0, z_0)$.
  This holds uniformly in $\mathcal U$.
  By Lemma~\ref{l:prob-W-to-U} for any $(I_0, z_0) \in \mathcal U$ and any $\varphi_0, \lambda_0 \in [0, 2\pi]$ we have~\eqref{e:lim-U}.
  Proposition~\ref{p:prob} follows from this statement.
\end{proof}

\section{Systems without capture into resonances} \label{s:no-capture-proof}
In this section we sketch a proof of Remark~\ref{r:no-capture}.
The lemma below is used to estimate the measure of trajectories that come too close to the saddle of perturbed system.
\begin{lemma} \label{l:pass-saddle}
  Take small enough $R > 0$ and define a neighborhood $U$ of the saddle of unperturbed system as
  \begin{equation}
    \{ (p, q, z) : |p-p_C(z)|, |q-q_C(z)| \le R, z \in \mathcal Z_0 \}.
  \end{equation}
  Then for any $r>1$ there exists $C_1 > 0$ such that the measure of the set of initial data $(p, q, z, \lambda) \in U \times [0, 2\pi]$ such that the corresponding solution of perturbed system~\eqref{e:perturbed-pq} does not leave $U$ is at most $\varepsilon^r$.
\end{lemma}
\begin{proof}[Sketch of proof]
  Let us use Moser's normal form near the saddle, it provides us coordinates $x, y$ such that
  unperturbed system rewrites as
  in these coordinates
  \begin{equation}
    \dot x = a(h, z) x, \qquad \dot y = -a(h, z) y, \qquad \dot z = 0, \qquad \dot \lambda = 1
  \end{equation}
  and perturbed system rewrites as
  \begin{equation}
    \dot x = a(h, z) x + \varepsilon f_x, \qquad \dot y = -a(h, z) y + \varepsilon f_y, \qquad \dot z = \varepsilon f_z, \qquad \dot \lambda = 1,
  \end{equation}
  where $f_x, f_y, f_z$ are smooth functions depending on $x, y, z, t$.
  We assume that $U$ is covered by this chart.

  The $x$ direction is expanding for the unperturbed system.
  The cone field
  \begin{equation}
    |dx| \ge \sqrt{dy^2 + dz^2 + d\lambda^2}
  \end{equation}
  is invariant by the flow of the perturbed system for small enough $\varepsilon$ and curves tangent to this cone field are expanded. Take $\Delta t \sim 1$ such that flow over time $\Delta t$ expands such curves by a factor at least $e=2.718...$.

  Take $C > 3 r \Delta t$.
  Cut the phase space into curves $y, z, t = const$ parametrized by $x$. Over time $C |\ln \varepsilon|$ (at least $2r \left \lceil{|\ln \varepsilon|}\right \rceil$ times $\Delta t$) such curves are expanded by a factor at least $e^{2r |\ln \varepsilon|} = \varepsilon^{-2r}$. This means (we use that curve length is equivalent to the measure of projection on $x$ axis, as the curve is tangent to cone field) that the length of the part of the curve that stays in $U$ after such time passes is
  $\lesssim \varepsilon^{2r}$. Integrating over all such curves and taking into account $\varepsilon^{2r} \ll \varepsilon^r$ gives the statement of the lemma.
\end{proof}

Let us say that a trajectory is \emph{captured in $U$}, if (over times $\sim \varepsilon^{-1}$) at some point it enters $U$ and then does not leave $U$ for time at least $2 C_1 |\ln \varepsilon|$, or leaves $U$ through $z \in \partial \mathcal Z_0$.
Let us show that total measure of initial data captured in $U$ is $O(\varepsilon^r)$. We will say that such points form exceptional set of initial data and exclude them from consideration.

Let us deduce this from the lemma above. First, it is possible to show that solutions cannot leave $U$ after time $O(|\ln \varepsilon|)$ via $z \in \partial \mathcal Z_0$, as we assume that along considered solutions of averaged system $z$ stays in $\mathcal Z$ that is $O(1)$-far from $\partial \mathcal Z_0$. Thus it is enough to prove simply that most solutions of perturbed system leave $U$ in due time.

Take $C_1$ determined by Lemma~\ref{l:pass-saddle} so that during one passage through $U$ the measure of initial data in $U$ that stays there for time at least $C_1 |\ln \varepsilon|$ is at most $\varepsilon^{r+1}$. Denote by $V \subset U$ this set. If trajectory starting at some point $X$ (possibly outside $U$) spends time at least $2C_1 |\ln \varepsilon|$ in $U$, then it spends time at least $C_1 |\ln \varepsilon|$ in $V$ and thus lies in a preimage of $V$ under some iterate of the time $C_1 |\ln \varepsilon|$ flow of perturbed system. We consider times $\sim \varepsilon$, thus all initial data such that the corresponding solution stays in $U$ for (continuous) time at least $C_1 |\ln \varepsilon|$ (during times $\sim \varepsilon^{-1}$) are covered by $\sim \varepsilon^{-1} \ln^{-1} \varepsilon$ such preimages. By Lemma~\ref{l:volume} measure of each preimage is $\sim \varepsilon^{r+1}$. This proves the estimate on measure of trajectories captured in $U$.

Consider immediate neighborhood of separatrices with width
$\sim \varepsilon \ln^5 \varepsilon$.
The condition $f_h < 0$ allows to make a stronger estimate for time of passage through this zone than for general case. We assume that for all $z$ the neighborhood $U$ of the saddle $C$ does not conatain a whole separatrix $l_1$ or $l_2$ of saddle $C$.
This means that a trajectory of perturbed system can be split as follows: passage near $l_1$ (outside $U$), passage through $U$, passage near $l_2$, passage through $U$, and so on.
Each passage through $U$ takes time at most $2C_1 |\ln \varepsilon|$ by our definition of exceptional set.
This process can only terminate when the solution leaves the immediate neighborhood of separatrices.
Outside $U$ we have $f_h < -c_1 < 0$ for some $c_1 < 0$ (by compactness). Thus during each passage near $l_1$ or $l_2$ the value of $h$ decreases by $\gtrsim \varepsilon$ and the number of such passages is $\lesssim \ln^5 \varepsilon$.
Passages near separatrices take time $O(1)$ and passage through $U$ takes time $O(\ln \varepsilon)$, thus the whole passage through immediate neighborhood of separatrices takes time $O(\ln^6 \varepsilon)$. As time derivatives of slow variables are $O(\varepsilon)$ for both perturbed and averaged systems, passage through immediate neighborhood of separatrices leads to $O(\varepsilon \ln^6 \varepsilon) \ll \sqrt{\varepsilon}$ deviation between solutions of perturbed and averaged system.

It is also possible to improve the estimate on passage through resonance zones. Recall that width of resonant zones has order $\delta_s$, where $\delta_s$ is given by~\eqref{e:delta}. Outside the immediate neighborhood (in $\mathcal B_3$) of separatrices $h$ solution of perturbed system winds around the separatices, during each wind $h$ decreases by $\sim \varepsilon$ and this wind takes time $\sim \ln h$. Thus on average $h$ decreases with speed $O(\ln^{-1} h)$. As $\pdv{\omega}{h} \sim h^{-1} \ln^{-2} h$, this gives $\dot \omega \sim \varepsilon h^{-1} \ln^{-3} h$ (on average).
Thus crossing resonant zone with width $\sim \delta_r$ takes time
$\sim \varepsilon^{-1}(\sqrt{\varepsilon b_s h \ln^2 h} + \varepsilon \ln^2 \varepsilon)$.
Taking sum over all resonances as in Section~\ref{ss:est-sum} gives that total time of crossing all resonant zones is $O(\varepsilon^{-1/2})$, which gives total change of slow variables in resonant zones $O(\sqrt{\varepsilon})$.

Combining improved estimates for accuracy of averaging method in immediate neighborhood of separatrices and in resonant zones with estimates for nonresonant zones as for general case (as done in Section~\ref{s:approach-proof} for general case) gives the estimate in Remark~\ref{r:no-capture}.

\section{Acknowledgment}
We are greatful to
A.V. Bolsinov for advices on integrable systems and to A.V. Artemyev and V.V. Sidorenko for useful discussion.

\begin{appendices}
\section{Analytic continuation: proofs} \label{a:AC-proofs}

\begin{proof}[Proof of Lemma~\ref{l:period-cont}]
  \begin{figure}[H]
  \centering
  \begin{tikzpicture}[scale=0.8]
    \draw[ultra thick] (0, 2) -- (1.5, 2) node[above right]{$\Gamma_1: y=1$};
    \draw[ultra thick] (2, 0) -- (2, 1.5) node[below right]{$\Gamma_2: x=1$};
    \draw[ultra thick] (0, -2) node[right]{$\Gamma_3: y=-1$} -- (-1.5, -2);
    \draw[ultra thick] (-2, 0) -- (-2, -1.5) node[left]{$\Gamma_4: x=-1$};
    \draw[semithick] (0, 2) -- (0, -2) to[out=270, in=225] (2.5, -2.5) %
      to [out=45, in=0] (2, 0) -- (-2, 0) to[out=180, in=225] (-2.5, 2.5)
      to [out=45, in=90] (0, 2);
    \draw[semithick] (1, 1) to[out=-45, in=180] (2, 0.7)
      to [out=0, in=45] (3, -3) to[out=-135, in=-90] (-0.7, -2)
      to [out=90, in=-45] (-1, -1);
    \draw[semithick] (-1, -1) to [out=135, in=0](-2, -0.7)
      to[out=180, in=-135] (-3, 3) %
      to [out=45, in=90] (0.7, 2)
      to[out=-90, in=135] (1, 1);
    \draw[arrows=->,semithick, dashed](0, -4)--(0,4) node[above right]{$y$};
    \draw[arrows=->,semithick, dashed](-4, 0)--(4, 0) node[above right]{$x$};
    \draw[-{To[length=3mm,width=2mm]}, thick](0, 2)--(0, 1);
    \draw[-{To[length=3mm,width=2mm]}, thick](0, -2)--(0, -1);
    \draw[-{To[length=3mm,width=2mm]}, thick](0, 0)--(1.3, 0);
    \draw[-{To[length=3mm,width=2mm]}, thick](0, 0)--(-1.3, 0);
  \end{tikzpicture}
  \caption{The transversals.}
  \label{f:transversals}
\end{figure}
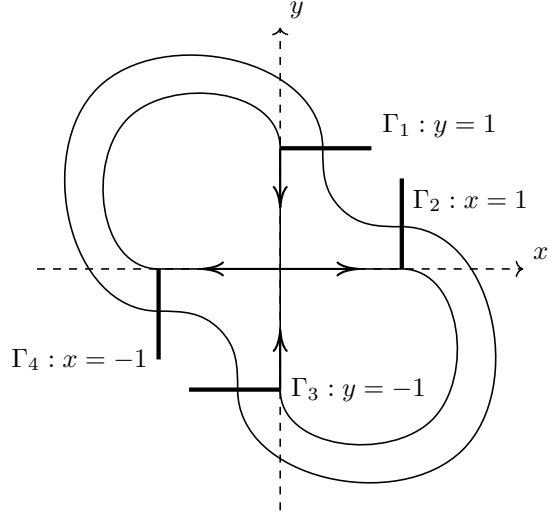

  We will use Moser's normal form~\cite{moser1956analytic} in a neighborhood of $C$.
  There are coordinates $x$ and $y$ such that our system can be written as
  \begin{align} \label{e:Moser}
  \begin{split}
    \dot x &= a(h, z) x, \\
    \dot y &= - a(h, z) y, \\
  \end{split}
  \end{align}
  where $a > 0$ for real $h, z$.
  Rescaling $p, q, x, y$ if necessary, we may assume that new coordinates are defined for all $(x, y) \in \mathbb C^2$ with $|x|, |y| \le 1$.
  Let us consider four transversals $\Gamma_1, \dots, \Gamma_4$ given by $x=\pm 1$ and $y=\pm 1$ as shown in Figure~\ref{f:transversals}.
  We can split $T = T_{12} + T_{23} + T_{34} + T_{41}$, where $T_{ij}$ is the time from $\Gamma_i$ to $\Gamma_j$. Note that as all paths connecting $(h, z)$ and $(h_0, z_0)$ are homotopic to each other (we assume $c_2 < 0.5$, then $h > 0$), the functions $T$ and $T_{ij}$ are continued as single-valued functions.

  The functions $T_{23}$ and $T_{41}$ are holomorphic even for $h=0$. Let us denote $\tilde h = xy$. As it is a first integral of~\eqref{e:Moser}, we have $\tilde h = \tilde h(h, z)$.
  Using~\eqref{e:Moser}, we can compute $T_{12} = T_{34} = - a^{-1} \ln \tilde h$; the branch of the logarithm here is obtained by analytic continuation of the real logarithm.
  As $\tilde h = 0$ for $h=0$, we can write $\tilde h(h, z) = h \tilde h_0(h, z)$ with analytic $h_0$ and thus $T_{12} = T_{34} = -a^{-1} \ln h - a^{-1} \ln \tilde h_0$. Thus for $z$ close to $z_0$ we have $T = A(h, z) \ln h + B(h, z)$, where $A = -2a^{-1}$ and $B$ are holomorphic, $A \ne 0$.

  Let us now prove that $\hat h(z)$ is uniquely defined for any $z$ close to $z_0$. Fix such $z$ and denote $F(h) = T(h) - T_0$ and $F_1(h) = A(h_0, z)\ln h + B(h_0, z) - T_0$. Let $\hat h_1$ be determined by $F_1(\hat h_1) = 0$. Consider the curve $h_\beta = \hat h_1(1 + \alpha e^{i \beta})$, where $\beta \in [0, 2\pi]$ is a parameter and $\alpha > 0$ is a small fixed number.
  We have
  $|F(h_\beta) - F_1(h_\beta)| = O(\alpha |\hat h_1|)(|\ln \hat h_1| + 1)$ and
  $|F_1(h_\beta)| = |A(h_0, z) \ln(1 + \alpha e^{i \beta})| \gtrsim \alpha$.
  As $|\hat h_1|$ is small (we have $\hat h_1 = e^{A^{-1}(T_0 - B)}$, here $T_0$ is large positive number for small $h_0$, $B$ is bounded, and $A$ is close to $A(z_0)$ which is a real negative number), we have $|F_1(h_\beta)| > 2 |F(h_\beta) - F_1(h_\beta)|$.
  Thus, by Rouche's theorem the equation $F(h) = 0$ has a unique solution $\hat h(z)$ in the region bounded by $h_\beta$.

  From $|h-\hat h| < c_2 |\hat h|$ we have
  \[
    \ln h - \ln \hat{h} = \ln(1 + (h - \hat{h})/\hat{h}) = O(c_2)
  \]
  and
  \[
    T(h, z) = O(c_2) + A(h, z) \ln \hat h + B(h, z)
    = A(\hat h, z) \ln \hat h + B(\hat h, z) + O(c_2)(1 + |\ln \hat h|)
    = T_0(1 + O(c_2)),
  \]
  as claimed above.

  Finally, the action $I(h, z)$ is defined by $I = \int p dq$.
\end{proof}

\begin{proof}[Proof of Lemma~\ref{l:cont}]
  Let us first state some helper statements. The local uniqueness and existence theorem (see e.g.~\cite[Theorem~$1.1$]{ilyashenko2008lectures}) claims that in a small neighborhood of any initial condition $y_0, t_0$ the solution of a complex ODE $\dot y = g(y, t)$ exists and depends holomorphically on the initial condition. From this we have the following corollary.
  \begin{corollary}
    For any complex ODE $\dot y = g(y, t)$ for any point $y_0$ there are $\Delta t > 0, \Delta y > 0$ such that for any $y_1$ with $|y_1 - y_0| < \Delta y$ the solution $y(t)$ with $y(0) = y_1$ exists for all complex $t$ with $|t| < \Delta t$.
  \end{corollary}
  \noindent By a compactness argument we get
  \begin{corollary} \label{c:continue}
    Consider a complex ODE $\dot y = g(y, t)$ in some open domain $U$. Consider a compact set $V \subset U$. Then there is $\Delta t > 0$ such that for any $y_0 \in V$ the solution of our ODE with the initial condition $y(0) = y_0$ exists and stays in $U$ for all $t$ with $|t| < \Delta t$.
  \end{corollary}
  \noindent We are now ready to prove Lemma~\ref{l:cont}. We will use the notation of the proof of Lemma~\ref{l:period-cont} such as $\Gamma_i$, $a$ and $T_{ij}$.
  Given $(h, z) \in \pi \mathcal D$ (here $\pi$ is the projection along the $\varphi$ axis),
  let us consider a continuous piecewise-linear path
  $t_s : [0, 1] \mapsto \mathbb C$
  with
  \[
    t_0 = 0, \; t_{1/4} = T_{12}(h, z), \; t_{1/2} = T_{12} + T_{23}, \;
    t_{3/4} = T_{12} + T_{23} + T_{34}, \; t_1 = T
  \]
  that is linearly interpolated between the points $0, 1/4, 1/2, 3/4, 1$.
  Then if we consider the cycle $r(t_s)$ lying in the complex solution of the unperturbed system with given $h, z$ starting with $r(t_0) \in \Gamma_1$, the path $r(t_s)$ then crosses $\Gamma_2, \Gamma_3, \Gamma_4$ for $s=1/4, 1/2, 3/4$ and returns to the same point on $\Gamma_1$ for $s=1$.

  Let us check that the segment $[0, T]$ lies in $c_2$-neighborhood of the path $t_s$ with the constant $c_2 > 0$ that can be made as small as needed by reducing $c$.
  Pick $\alpha \in \mathbb C, \abs{\alpha} = 1$ such that $\alpha T_{12} = \alpha T_{34} \in \mathbb R$. As $T_{23}$ and $T_{41}$ are analytic functions of $h, z$, we have
  $\Im \alpha T_{ij} < O(c)$.
  This means that $\Im \alpha T = O(c)$ and $\Im \alpha t_s = O(c)$ for all $s$, so $[0, T]$ lies in $O(c)$-neighborhood of the path $t_s$.

  Let us choose a neighborhood $V$ (in the space with coordinates $p, q, z$) of the union of the real separatrices of the system for $z = z_0$ such that $r(t_s) \in V$ for all considered $h, z$. We take any $V$ that contains the set given by $x, y \in [-1, 1]$ for all $z$ with $|z-z_0| < c_2$. This implies that $r(t_s)$ stays in $V$ for $s \in [0, 1/4]$ and $s \in [1/2, 3/4]$.
  For other values of $s$ the point $r(t_s)$ also stays $O(c)$-close to the union of the real separatrices of the system for $z = z_0$, so for small enough $c$ it lies in $V$.

  Now we apply Corollary~\ref{c:continue} (with $y = (p, q, z)$, $\dot z = 0$) to the closure of $V$, it gives us some $\Delta t$.
  Denote by $t = \frac{\varphi T}{2\pi}$ the time for the unperturbed system.
  Let us continue $r(h, z, \varphi(t))$ as a function of $t$. Clearly, it is defined for $t \in \{t_s\}$. By Corollary~\ref{c:continue} it can be continued to $\Delta t$-neighborhood of $\{t_s\}$.
  We have proved above that for small enough $c$ we have that $\Delta t/2$-neighborhood of $\{t_s\}$ covers $[0, T]$. This means that $\Delta t$-neighborhood of $\{t_s\}$ covers $\Delta t/2$-neighborhood of $[0, T]$, so $r$ can be continued to the latter neighborhood. Returning to $\varphi$, we get that $r$ can be continued to $\pi \Delta t / \abs{T}$-neighborhood of $[0, 2\pi]$.
  As from Lemma~\ref{l:period-cont} we have $\abs{T} \sim T(h_0, z_0)$,
  this estimate proves the lemma.
\end{proof}

\begin{proof}[Proof of Lemma~\ref{l:int-psi}]
  We can rewrite
  \[
    \omega^{-1} \int_{\varphi=0}^{\varphi_0} \psi(h, z, \varphi, \lambda) d\varphi =
    \int_{t=0}^{t_0} \psi(h, z, \omega t, \lambda) dt,
  \]
  where $t_0 = \omega \varphi_0$.
  Let us again consider the contour $\{ t_s \}$ introduced in the proof of Lemma~\ref{l:cont}. As shown in this lemma, $\delta$-neighborhood of $\{ t_s \}$ covers $[0, T]$ for some $\delta = O(1)$.
  As we have $\abs{\omega^{-1}}\Im(\varphi_0) = O(1)$, we can connect the point $t_0$ with the point $\omega^{-1} \Re \varphi_0 \in [0, T]$ and then with some point $t_{s_0}$ on $\{ t_s \}$ by a path of length $O(1)$ on the complex plane of the values of $t$. As $\psi$ is bounded, we have $\int_{t_{s_0}}^{t_0} \psi dt = O(1)$.
  Thus, it is enough to show that
  \[
    \int_{t=0}^{t_{s_0}} \psi dt = O(1)
  \]
  along the contour $\{ t_s, 0 \le s \le s_0 \}$. This contour is split into up to four parts by the transversals $\Gamma_i$ (cf. Figure~\ref{f:transversals}).
  For the parts when the solution is far from the saddle, the integral along these parts is $O(1)$.
  Let us prove that the integral along the parts near the saddle is also $O(1)$. For defineteness, let us consider the part between $\Gamma_1$ and $\Gamma_2$ (or some fragment of this part if it is included in our contour only partially). During this part the solution is near the saddle, and we can use the $x, y$ chart~\eqref{e:Moser}.
  As $\psi(C) = 0$ ($C$ corresponds to $x=y=0$), we can write $\psi(x, y, z, \lambda) = x \psi_x + y \psi_y$, where $\psi_x, \psi_y$ are analytic. As $\dot x = a x$, we can write (we omit the limits of integration, as they are not important for the $O(1)$ estimate)
  \[
    \int x \psi_x dt = a^{-1} \int \psi_x dx = O(1).
  \]
  Similarly, we have $\int y \psi_y dt = O(1)$.
  Thus, $\int_{t=0}^{t_{s_0}} \psi dt = O(1)$, as required.
\end{proof}

\section{Estimates on Fourier coefficients} \label{a:fourier_proof}
Lemma~\ref{l:cont} can be used to estimate Fourier coefficients near the separatrices.
Let $\psi(p, q, z, \lambda)$ and $\psi_0(p, q, z, \lambda)$ be analytic functions in $\mathcal{\tilde B}$ with $\psi_0 = 0$ at the saddle $C(z)$ for all $\lambda, z$.
\begin{corollary} \label{l:decay-phi}
  There is $K > 0$ such that for any real $(h, z) \in \mathcal B$, $\lambda \in [0, 2\pi]$, $k \in \mathbb Z$ we have
  \begin{align}
  \begin{split}
    \bigg|\int_0^{2\pi} \psi(h, \varphi, z, \lambda) e^{ik\varphi} d\varphi \bigg|
    \le K \exp(-c |k| / T(h, z)), \\
    \bigg|\int_0^{2\pi} \psi_0(h, \varphi, z, \lambda) e^{ik\varphi} d\varphi \bigg|
    \le K T^{-1}(h, z) \exp(-c |k| / T(h, z)).
  \end{split}
  \end{align}
\end{corollary}
\begin{proof}
  The first inequality is proved exactly like the exponential decay of Fourier coefficients of an analytic function.
  We can move the contour of integration up (assuming $k > 0$, otherwise down) by adding $i c / T(h, z)$ to $\varphi$. By the periodicity it will not change the integral. But for the new integral $|e^{ik\varphi}| = \exp(-c k / T(h, z))$ while $|\psi(h, \varphi, z, \lambda)|$ is bounded, as $\psi$ is bounded in the compact set $\tilde B$.

  The second inequality is proved by the same shift of the contour of integration, but we now take into account that $\psi_0$ is small near $C$.
  Denote $\tilde \varphi(t) = 2\pi (t + \frac{ic}{2\pi})/T, \varphi(t) = 2\pi t / T$.
  Rewriting our integral as integral $dt$, we get
  \begin{equation} \label{e:Moser-solutions}
    \omega \int_0^T e^{2\pi ik t/T} e^{-kc/T} \psi_0(h, \tilde \varphi, z, \lambda) dt.
  \end{equation}
  In a neighborhood $\mathcal S$ of $C$ given by $|x|, |y| < 1$ we can use~\eqref{e:Moser}. Solving this system, we get
  \[
    x(t) = e^{a(t-t_0)}x(t_0), \qquad
    y(t) = e^{-a(t-t_0)}y(t_0).
  \]
  Here $a > 0$, as $h$ and $z$ are real.
  Taking $t-t_0 = ic/(2\pi)$, we get
  $|x(\varphi)| = |x(\tilde \varphi)|$
  and
  $|y(\varphi)| = |y(\tilde \varphi)|$.
  Let us now show that the integral of $|\psi_0|$ inside $\mathcal S$ is $O(1)$.
  Indeed, as $\psi_0(C) = 0$ ($C$ corresponds to $x=y=0$), we can write $\psi_0(x, y, z, \lambda) = x \psi_x + y \psi_y$, where $\psi_x, \psi_y$ are analytic.
  Let $[t_1, t_2] \subset [0, T]$ be a segment such that the solution $(x(\varphi(t)), y(\varphi(t)))$ is in $\mathcal S$ for all $t \in [t_1, t_2]$. Then $(x(\tilde \varphi(t)), y(\tilde \varphi(t)))$ is also in $\mathcal S$.
  We may use $x = x(\varphi(t))$ as an independent variable, from $\dot x = a x$ we have
  $dt = dx / (ax)$.
  Thus we have
  \[
    \int_{t_1}^{t_2} |x(\tilde \varphi) \psi_x(\tilde \varphi)| dt
    =
    \int_{t_1}^{t_2} |\psi_x(\tilde \varphi)||x(\varphi)| dt
    =
     a^{-1} \int_{x(t_1)}^{x(t_2)} |\psi_x(\tilde \varphi)| dx = O(1).
  \]
  Similarly, we have $\int_{t_1}^{t_2} |y \psi_y| dt = O(1)$.
  Thus, $\int_{t_1}^{t_2} |\psi_0| dt = O(1)$, as required.

  It is clear that the integral of $|\psi_0|$ outside $\mathcal S$ is also $O(1)$, as the solution only spends time $O(1)$ there. Thus, $\int |\psi_0(\tilde \varphi)| dt = O(1)$ and $\int |\psi_0(\tilde \varphi)| d\varphi = O(T^{-1})$.
  Together with~\eqref{e:Moser-solutions} this yields the second inequality.
\end{proof}

\begin{corollary} \label{l:decay-phi-2}
  There are $C_1, C_2 > 0$ such that for any $(h, z) \in \mathcal B$, $k, l \in \mathbb Z$ we have
  \begin{align} \label{e:fourier-corollary}
  \begin{split}
    \bigg|\int_0^{2\pi} \int_0^{2\pi} \psi(h, \varphi, z, \lambda) e^{ik\varphi} e^{il \lambda} d\varphi d\lambda \bigg|
    \le C_1 \exp\bigg( - C_2 |l| - C_2 \frac{|k|}{T} \bigg), \\
    \bigg|\int_0^{2\pi} \int_0^{2\pi} \psi_0(h, \varphi, z, \lambda) e^{ik\varphi} e^{il \lambda} d\varphi d\lambda \bigg|
    \le C_1 T^{-1} \exp\bigg( - C_2 |l| - C_2 \frac{|k|}{T} \bigg).
  \end{split}
  \end{align}
\end{corollary}
\begin{proof}
  Let us prove the second part, the first one can be obtained similarly.
  Denote by $L$ the left hand side of the second line of~\eqref{e:fourier-corollary}. It is enough to obtain two separate estimates
  \begin{equation}
    L \le C_1 T^{-1} \exp ( - 2 \; C_2 |l| ), \qquad
    L \le C_1 T^{-1} \exp \Big( - 2 \; C_2 \frac{|k|}{T} \bigg).
  \end{equation}
  Shifting the contour of integration in $\lambda$ by $c i$ (as in the usual argument for the exponential decay of Fourier coefficients of an analytic function that we have already used in Corollary~\ref{l:decay-phi}), we get
  \[
    L \lesssim \exp ( - c l ) \int_0^{2\pi} \int_0^{2\pi} |\psi_0(h, \varphi, z, \lambda + c i)| d\varphi d\lambda.
  \]
  Arguing as in the proof of Corollary~\ref{l:decay-phi}, we can show that for all $\lambda$ we have $\int_0^{2\pi} |\psi_0(h, \varphi, z, \lambda + c i)| d\varphi = O(T^{-1})$. This means $L \lesssim T^{-1} \exp(-cl)$ and thus the first required estimate holds.

  Multiplying the estimate from Lemma~\ref{l:decay-phi} by $e^{il \lambda}$ and integrating by $\lambda$, we obtain the second required estimate.
\end{proof}

\begin{proof}[Proof of Lemma~\ref{l:fourier}]
  As $\pdv{I}{h} = \omega^{-1}$, the Fourier coefficients of $f_I$ can be expressed via the Fourier coefficients of $f_h$: $f_{I, m} = \omega^{-1} f_{h, m}$. Thus the first part follows from Corollary~\ref{l:decay-phi-2}.

  Let us prove the second part. Let us first estimate $\pdv{f_{h, m}}{h}$.
  We have
  \begin{equation}
    \norm{\pdv{f_{h, m}}{h}}
    \sim
    \norm{
      \int_{0}^{2\pi} \int_{0}^{2\pi}
      \pdv{f_h}{h}\/ e^{-i(m_1 \varphi + m_2 \lambda)}
      d\lambda d \varphi
    }
    \lesssim
    \int_{0}^{2\pi} \norm{
      \int_{0}^{2\pi}
      \pdv{f_h}{h}\/ e^{-i m_2 \lambda}
      d\lambda
    } d \varphi
  \end{equation}
  By~\eqref{e:est-real} we have $\pdv{p}{h}$, $\pdv{q}{h} = O_*(h^{-1} \ln^{-1} h)$. Note that these expressions do not depend on $\lambda$ and $\int_0^{2\pi} O_*(h^{-1} \ln^{-1} h) d \varphi = O(h^{-1} \ln^{-2} h)$ by~\eqref{e:int-O-star}.
  As $\pdv{f_h}{h} = \pdv{f_h}{p} \pdv{p}{h} + \pdv{f_h}{q} \pdv{q}{h}$, we may continue the estimate above as follows.
  \begin{equation}
    \norm{\pdv{f_{h, m}}{h}}
    \lesssim
    h^{-1} \ln^{-2} h \;
    \max_{\varphi} \int_{0}^{2\pi}
    \norm{\pdv{f_h}{p}\/ e^{-i m_2 \lambda}}
    +
    \norm{\pdv{f_h}{q}\/ e^{-i m_2 \lambda}}
    d\lambda
    \lesssim
    h^{-1} \ln^{-2} h \; \exp(- C_F |m_2|),
  \end{equation}
  as the Fourier coefficients of smooth functions $\pdv{f_h}{p}$ and $\pdv{f_h}{q}$ decrease exponentially.
  As $f_{I, m} = \omega^{-1} f_{h, m}$, this gives
  \[
    \norm{\pdv{f_{I, m}}{h}} \lesssim \omega^{-1} \norm{\pdv{f_{h, m}}{h}} + \norm{f_{h, m}} \; \pdv{}{h}\/ (\omega^{-1})
    \lesssim |h^{-1} \ln^{-1} h| \; \exp(- C_F |m_2|).
  \]
  The estimate $\norm{\pdv{f_{h, m}}{z}} \lesssim h^{-1} \ln^{-2} h \; \exp(- C_F |m_2|)$ is obtained in the same way as the estimate for $\norm{\pdv{f_{h, m}}{h}}$.
\end{proof}

\section{Proof of estimates on $u$} \label{a:appendix}
The following three lemmas will be used to prove Lemma~\ref{l:est-u}.
Fix a resonance $s_2/s_1$.
In all three lemmas we will assume that $s_2/s_1$ is the nearest to $\omega$ resonance.
We will use the notation $m = (m_1, m_2) \in \mathbb Z^2$.
In summation over $m$ we will often need to skip the vectors $m$ that are equal to $(\nu s_1, -\nu s_2)$ for some $\nu \in \mathbb Z$. This will be denoted by an upper index $s$ in the summation symbol.
\begin{lemma} \label{l:big-sum}
  \begin{equation} \label{e:big-sum}
    \sideset{}{^{(s)}}\sum_{1 \le |m| \le N}
    \norm{f_m} |m_1 \omega + m_2|^{-1}(|m_2| + 1)
    \lesssim
    |s_1| + |\ln h| \ln |\ln \varepsilon|.
  \end{equation}
  Moreover, this holds not only for the Fourier coefficients $f_m$, but for any non-negative numbers $\norm{f_m}$ with $\norm{f_m} \lesssim e^{-C_F |m_2|}$.
\end{lemma}
\begin{proof}
  Denote by $S_{m_2}$ the part of the left-hand side of~\eqref{e:big-sum} with this $m_2$.
  Note that we have $m_1 \omega + m_2 \ne 0$.
  For fixed $m_2$ let $m^+_1$ ($m^-_1$) be the values of $m_1$ corresponding to the smallest in absolute value positive (negative) value of $m_1 \omega + m_2$.
  Here we consider \emph{all} integer values of $m_1$, including the ones with $|m_1| > N$, so we may have $|m^\pm_1| > N$.
  Let $A^+_{m_2}$ ($A^-_{m_2}$) denote the corresponding term in $S_{m_2}$ if it exists (i.e. $1 \le |m| \le N$), or $0$ otherwise.
  Denote by $B_{m_2}$ the sum of all other terms, $S_{m_2} = A^+_{m_2} + A^-_{m_2} + B_{m_2}$.

  As in~\cite[Proof of Lemma $7.1$]{neishtadt2014}, for all $1 \le |m| \le N$ we have
  $|m_1 \omega + m_2| \ge (4s_1)^{-1}$.
  This means $A^+_{m_2} \lesssim s_1 (|m_2| + 1)e^{-C_F |m_2|}$ and $\sum_{m_2} A^+_{m_2} \lesssim s_1$.
  Similarly, $\sum_{m_2} A^-_{m_2} \lesssim s_1$.

  For fixed $m_2$ set $k(m_1)= m_1 - m_1^+$ for $m_1 > m_1^+$ and $k(m_1)= m^-_1 - m_1$ for $m_1 < m^-_1$. We have $|m_1 \omega + m_2| \ge \omega k$, so
  \[
    B_{m_2} \lesssim 2 \omega^{-1} (|m_2| + 1) e^{-C_F |m_2|}\sum_{k = 1}^{2N} k^{-1}
    \lesssim \omega^{-1} \ln |\ln \varepsilon| \; (|m_2| + 1) e^{-C_F |m_2|}
  \]
  and $\sum_{m_2} B_{m_2} \lesssim \omega^{-1} \ln |\ln \varepsilon|$.
\end{proof}
\begin{lemma} \label{l:big-sum-2}
  For the Fourier coefficients $f_m$ of $f$ we have
  \begin{equation} \label{e:big-sum-2}
    \sideset{}{^{(s)}}\sum_{1 \le |m| \le N,}
    \norm{f_m} |m_1 \omega + m_2|^{-1}|m_1| \lesssim |s_1 \ln h| + \ln^2 h \ln |\ln \varepsilon|.
  \end{equation}
\end{lemma}
\begin{proof}
  Let us argue as above, adapting the notation $A^\pm_{m_2}$ and $B_{m_2}$ to the corresponding terms in~\eqref{e:big-sum-2}. We now have
  (we use $|m^+_1\omega + m_2| \le \omega$, so $|m^+_1 + \omega^{-1} m_2| \le 1$)
  \[
    A^+_{m_2} \lesssim s_1 |m^+_1| e^{-C_F |m_2|} \lesssim s_1 (\omega^{-1} |m_2|+1) e^{-C_F |m_2|}
  \]
  and $\sum_{m_2} A^+_{m_2} \lesssim |s_1 \ln h|$.
  Similarly, $\sum_{m_2} A^-_{m_2} \lesssim |s_1 \ln h|$.

  We also have
  \[
    \abs{\frac{m_1}{m_1 \omega + m_2}}
    = \abs{\omega^{-1} - \frac{m_2 \omega^{-1}}{m_1 \omega + m_2}} \lesssim |\omega^{-1}| + |\omega^{-1}||m_2||m_1 \omega + m_2|^{-1}.
  \]
  Thus we can split $B_{m_2} \lesssim B_{1, m_2} + B_{2, m_2}$, where
  the two summands correspond to the summands $\omega^{-1}$ and $\omega^{-1}|m_2||m_1 \omega + m_2|^{-1}$ above (in this order).
  By the estimate~\eqref{e:fourier} for $f_m$ we have $\sum_{|m_1| \le N} \norm{f_{(m_1, m_2)}} \lesssim e^{-C_F|m_2|} \omega^{-1}$, thus
  $\sum_{m_2} B_{1, m_2} \lesssim \omega^{-2}$.
  Denote by $\tilde B$ the terms $B$ in~\eqref{e:big-sum}.
  The estimate $\sum_{m_2} \tilde B_{m_2} \lesssim \omega^{-1} \ln |\ln \varepsilon|$ from the proof of~\eqref{e:big-sum} gives
  $\sum_{m_2} B_{2, m_2} \lesssim \omega^{-2} \ln |\ln \varepsilon|$.
\end{proof}

\begin{lemma}
  For the Fourier coefficients $f_m$ of $f$ we have
  \begin{equation} \label{e:big-sum-3}
    \sideset{}{^{(s)}}\sum_{1 \le |m| \le N}
    \norm{f_m} |m_1 \omega + m_2|^{-2}|m_1| \lesssim |\ln h| \ln^4 \varepsilon.
  \end{equation}
\end{lemma}
\begin{proof}
  Let us once again argue as in the proof of Lemma~\ref{l:big-sum}, adapting the notation $A^\pm_{m_2}$ and $B_{m_2}$ to the corresponding terms in~\eqref{e:big-sum-3}.
  As in the proof of Lemma~\ref{l:big-sum-2}, we have
  \[
    A^+_{m_2} \lesssim s_1^2 |m_1| e^{-C_F |m_2|} \lesssim s_1^2 \omega^{-1} |m_2| e^{-C_F |m_2|}.
  \]
  Thus, $\sum_{m_2} A^+_{m_2} \lesssim |\ln h| \ln^4 \varepsilon$ and, similarly, $\sum_{m_2} A^-_{m_2} \lesssim |\ln h| \ln^4 \varepsilon$.

  To estimate $B$, let us reuse the notation $k$ from the proof of Lemma~\ref{l:big-sum}. Recall that $|m_1 \omega + m_2| \ge \omega k$. We have
  \[
    B_{m_2} \lesssim N e^{-C_F |m_2|} \sum_{m_1} |m_1 \omega + m_2|^{-2} \lesssim e^{-C_F |m_2|} N \omega^{-2} \sum_{k} k^{-2} \lesssim e^{-C_F |m_2|} N \omega^{-2}.
  \]
  This implies $\sum_{m_2} B_{m_2} \lesssim \ln^2 h \ln^2 \varepsilon$.
\end{proof}

\begin{proof}[Proof of Lemma~\ref{l:est-u}]
  We have $u = A + B$ with
  \begin{equation}
    A = \sideset{}{^{(s)}}\sum_{1 \le |m| \le N, \; m \in \mathbb Z^2} \frac{f_m e^{i(m_1\varphi + m_2 \lambda)}}{i(m_1\omega + m_2)},
    \qquad
    B = \sum_{1 \le |\nu s| \le N, \; \nu \in \mathbb Z}
    \frac{f_{\nu s} e^{i\nu (s_1\varphi - s_2 \lambda)}}{i\nu(s_1\omega - s_2)}.
  \end{equation}
  We have $s_1 \le N \lesssim \ln^2 \varepsilon$ and $\omega \sim \ln h \lesssim \ln \varepsilon$.
  By~\eqref{e:big-sum} we get $\norm{A} \lesssim \ln^2 \varepsilon$.
  As $\Delta = |\omega - \xi_s|$, we have $\abs{s_1 \omega - s_2} = s_1 \Delta$.
  We also have from~\eqref{e:fourier} that $\norm{f_{\nu s}} \lesssim a_s^{\abs{\nu}}, a_s \le e^{-C_F} < 1$. Hence,
  $\norm{B} \lesssim a_s s_1^{-1} \Delta^{-1}$. This gives the estimate on $\norm{u}$.

  Taking $\pdv{}{\lambda}$ or $\pdv{}{\varphi}$ of the terms in $A$ multiplies them by $m_2$ or $m_1$, respectively.
  The estimates on $\pdv{A}{\lambda}$ and $\pdv{A}{\varphi}$ follow from~\eqref{e:big-sum} and~\eqref{e:big-sum-2}, respectively.
  We also have
  $\norm{\pdv{B}{\lambda}} \lesssim a_s s_2 s_1^{-1} \Delta^{-1}$ and
  $\norm{\pdv{B}{\varphi}} \lesssim a_s \Delta^{-1}$.

  Taking $\pdv{}{h}$ of $A$ or $B$ creates two terms. In one of them $f_m$ are replaced by $\pdv{f_m}{h}$, denote this term by $A_{h, 1}$ or $B_{h, 1}$.
  By Lemma~\ref{l:fourier} and Lemma~\ref{l:big-sum} we have
  $A_{h, 1} \lesssim h^{-1} \ln^{-1} h \ln^2 \varepsilon$.
  By Lemma~\ref{l:fourier} we also have $B_{h, 1} \lesssim b_s s_1^{-1} h^{-1} \ln^{-1} h \Delta^{-1}$.

  For the second term we have
  \begin{equation}
    A_{h, 2} = - \pdv{\omega}{h} \quad \sideset{}{^{(s)}}\sum_{1 \le |m| \le N}
    \frac{f_m m_1 e^{i (m_1\varphi + m_2 \lambda)}}
    {i(m_1\omega + m_2)^2},
    \qquad
    B_{h, 2} = - \pdv{\omega}{h} \sum_{1 \le |\nu s| \le N}
    \frac{f_{\nu s} s_1 e^{i\nu (s_1\varphi - s_2 \lambda)}}
    {i \nu (s_1\omega - s_2)^2}.
  \end{equation}
  We have $\pdv{\omega}{h} \sim h^{-1} \ln^{-2} h$.
  By~\eqref{e:big-sum-3} we have
  $\norm{A_{h, 2}} \lesssim |\pdv{\omega}{h}| \ln h \ln^4 \varepsilon \lesssim h^{-1} \ln^{-1} h \ln^4 \varepsilon$.
  As $\norm{f_{\nu s}} \lesssim b_s^{|\nu|}$, we have
  $\norm{B_{h, 2}} \lesssim b_s \pdv{\omega}{h} s_1^{-1} \Delta^{-2}$.

  The estimate on $\pdv{u}{z}$ is obtained in the same way.
\end{proof}

\begin{proof}[Proof of Lemma~\ref{l:d-f-0-I}]
  From the Hamiltonian equations we have $\pdv{h}{I} = \omega$.
  By~\cite[Corollary~3.2]{neishtadt17} we have
  \[
    \pdv{I}{z} = O(1), \; \pdv{I}{zh} = O(\ln h), \; \pdv[2]{I}{z} = O(1).
  \]
  As $\pdv{I}{h} = \omega^{-1}$, the first estimate implies $\norm{\pdv{w}{v}} = O(\ln h)$.
  We have
  $\pdv{I}{h} (\pdv{h}{z})_{I=const} + (\pdv{I}{z})_{h=const} = 0$, this gives $(\pdv{h}{z})_{I=const} = O(\ln^{-1} h)$ and $\norm{\pdv{h}{w}} = O(\ln^{-1} h)$.

  We have
  \begin{equation}
    f_{I, 0} = \pdv{I}{h} f_{h, 0} + \pdv{I}{z} f_{z, 0}.
  \end{equation}
  This rewrites as
  \begin{equation} \label{e:local-f-I-0}
    f_{I, 0} = (2\pi)^{-1} \oint f_h dt + \pdv{I}{z} f_{z, 0}.
  \end{equation}
  The contour integral is taken along the closed trajectory of the unperturbed system given by the values of $h, z$ and this trajectory is parametrized by the time $t$.
  By~\cite[Lemma~3.2]{neishtadt17} we have
  \begin{equation}
    \oint f_h dt = O(1), \qquad
    \pdv{}{h} \oint f_h dt = O(\ln h), \qquad
    \pdv{}{z} \oint f_h dt = O(1).
  \end{equation}
  Plugging the first estimate in~\eqref{e:local-f-I-0} gives $f_{I, 0} = O(1)$.
  From~\cite[Lemma~3.2]{neishtadt17} we also have
  \begin{equation}
    \pdv{f_{h, 0}}{h}, \pdv{f_{z, 0}}{h} = O(h^{-1}\ln^{-2}h), \qquad
    \norm{\pdv{f_{h, 0}}{z}} = O(\ln^{-1}h),
    \qquad \norm{\pdv{f_{z, 0}}{z}} = O(1).
  \end{equation}
  Plugging this in the derivatives of~\eqref{e:local-f-I-0} gives
  \begin{equation}
    \pdv{f_{I, 0}}{h} = O(h^{-1} \ln^{-2} h), \qquad
    \norm{\pdv{f_{I, 0}}{z}} = O(1).
  \end{equation}
  As $\pdv{h}{w} = O(\ln^{-1} h)$, we have
  $\norm{\pdv{f_{I, 0}}{w}} = \norm{\pdv{f_{I, 0}}{h} \pdv{h}{w} + \pdv{f_{I, 0}}{z}\pdv{z}{w}} = O(h^{-1} \ln^{-3} h)$.
  The estimate for $\norm{\pdv{f_{z, 0}}{w}}$ is obtained in the same way.
\end{proof}

\section{Proof of auxiliary lemma} \label{a:proof-aux}
\begin{proof}[Proof of Lemma~\ref{l:volume}]
  Recall that the divergence of a vector field $v$ with respect to a volume form $\alpha$ is a function $\Div_\alpha (v)$
  such that
  $\mathcal L_v(\alpha) = \Div_\alpha (v) \cdot \alpha$ (here $\mathcal L$ denotes the Lie derivative).
  Let $\alpha = dp \wedge dq \wedge dz \wedge d\lambda$ be the volume form.
  Set
  \[
    b(X, \lambda) = ((g^\lambda)^* \alpha)_X / \alpha_X,
  \]
  i.e. we take pullback of $\alpha$ by the flow and divide it by $\alpha$ at the point $X$. This gives a number, as the space of $n+3$-forms on $n+3$-manifold in the given point is one-dimensional.
  By definition of Lie derivative the number $b(X, \lambda)$ satisfies
  \begin{equation}
    \dv{b(X, \lambda)}{\lambda} =  \Div v(g^\lambda(X)) \; b(X, \lambda),
  \end{equation}
  where $v$ is the right-hand side of~\eqref{e:perturbed-pq}.
  As the Hamiltonian terms have zero divergence and $\dot \lambda = 1$, we have $\Div v = \varepsilon \Div f = O(\varepsilon)$. This shows that for $\lambda$ with
  $|\lambda| < \varepsilon^{-1} \Lambda$
  we have $b(X, \lambda) \in (0, C]$ with
  $C = \exp(\Lambda \max|\Div f|)$. Integrating over all $X \in A$ gives the required estimate.
\end{proof}

\section{Reduction of two-frequency systems to periodically perturbed one-frequency system} \label{a:two-freq-reduction}
In this appendix we present the proof of Lemma~\ref{l:two-freq-one-freq}. This proof was kindly communicated to us by A.V. Bolsinov. Then we show how this lemma can be used to reduce perturbations of two-frequency integrable systems to time-periodic perturbations of one-frequency systems.
\begin{proof}[Proof of Lemma~\ref{l:two-freq-one-freq}]
  We will consider the case without the parameter $z$, as with $z$ one can simply construct the new coordinates separately for each $z$ as described below.
  The proof is based on results presented in the book~\cite{bolsinov2004integrable}.

  Let $L$ denote the singular leaf of the Liouville foliation (i.e., the foliation of the phase space into \emph{Liouville tori} given by $H=H_0$, $F=F_0$; each Lioville torus is parametrized by the values $H_0$, $F_0$ of the two first integrals) and let $Q^3$ denote the isoenergy level that contains $L$.
  By~\cite[Theorem 3.2]{bolsinov2004integrable} there exists a \emph{periodic integral} $s_1$ defined in a four-dimensional neighborhood $V(L)$, i.e. a function $s_1$ that is smooth even on separatices and is such that the flow of the vector field $\sgrad s_1$ is $2\pi$-periodic (we use the notation $\sgrad U$ for the \emph{Hamiltonian vector field} of the function $U$, it is determined by $\omega(v, \sgrad U)=dU(v)$, where $\omega$ is the symplectic structure, $v$ is arbitrary tangent vector and $dU(v)$ denotes the derivative of $U$ in the direction $v$).
  The periodic integral $s_1$ allows to define the structure of a Seifert fibration in a neighborhood $U(L) \subset Q^3$ of the singular leaf $L$, the fibration of $U(L)$ by the orbits of the flow of $\sgrad s_1$ (\cite[Theorem 3.3]{bolsinov2004integrable}).

  As shown in \cite[Chapter 3]{bolsinov2004integrable}, there are two cases:
  \begin{enumerate}
    \item one can take a two-dimensional surface $P \subset U(L)$ that intersects each leaf of the Seifert fibration once;
    \item one can take a two-dimensional surface $\hat P \subset U(L)$ that intersects each regular leaf of the Seifert fibration twice and each singular leaf once.
  \end{enumerate}
  By~\cite[Proposition 5.4]{bolsinov2004integrable} topological stability of the isoenergy level $Q^3$ implies that $P$ and $\hat P$ can be taken transversal to $\sgrad H$.
  On $Q^3$ the vector fields $\sgrad s_1$ and $\sgrad H$ are tangent to $Q^3$. Thus we can continue $P$ and $\hat P$ to $3$-dimensional transversals $P^3, \hat P^3 \subset V(L)$.

  Let us consider the first case and construct a phase variable $\varphi_1$ conjugate to $s_1$. To do so, we set $\varphi_1=0$ on $P^3$ and propagate it along the leaves of Seifert fibration. Indeed, we want $\{ \varphi_1, s_1 \} = 1$ (here $\{\cdot, \cdot\}$ denotes the Poisson bracket), this condition can be rewritten in the following way: the derivative of $\varphi_1$ along $\sgrad s_1$ is $1$ and used to define $\varphi_1$.
  As trajectories of $\sgrad s_1$ are $2\pi$-periodic, this correctly defines the angle variable $\varphi_1$.

  In the second case we can construct $\varphi_1$ in the same way, the difference will be that $\varphi_1$ will be defined on a double cover. We will consider the lift of the unperturbed system to the covering space instead of the original unperturbed system in the rest of the proof.

  Let us now define variables $p, q$ so that $s_1, p, \varphi_1, q$ are canonical variables. Fix variables $p, q$ on some two-dimensional section $\{ \varphi_1 = 0, s_1 = \tilde s_1 \}$ so that these variables are canonical with respect to the restriction of the symplectic structure to this section. Spread these coordinates on the whole $V(L)$ by the flows of $\sgrad s_1$ and $\sgrad \varphi_1$ (these flows commute, as $\{ \varphi_1, s_1 \} = 1$). As these flows are symplectic, we have $\{ q, p \} = 1$ on $V(L)$. By construction we have $\{ a, b \} = 0$, where $a=\varphi_1, s_1$ and $b=p, q$. Thus $s_1, p, \varphi_1, q$ are canonical variables.

  In these new variables the dynamics of the unperturbed system rewrites as
  \begin{equation}
    \dot \varphi_1 = \pdv{H}{s_1}, \qquad \dot s_1 = 0, \qquad \dot p = -\pdv{H}{q}, \qquad \dot q = \pdv{H}{p}.
  \end{equation}

  Take $\varphi_1$ as a new independent variable (new time). Denote $\psi'=\frac{d\psi}{d\varphi_1}$.
  Denote $h=H$ and take $h$ as a new variable that replaces $s_1$.

  According to general formulas of isoenergetic reduction~\cite[\S9.45.B]{arnold1989mathematical}, the dynamics of $p$ and $q$ with respect to the time $\varphi_1$  is given by the Hamiltonian $S(p, q, h) = - s_1(p,q,h)$. Finally, denote $s=\varphi_1$. We have transformed the unperturbed system to the form

  \begin{equation}
    s' = 1, \qquad h' = 0, \qquad p' = -\pdv{S}{q}, \qquad q' = \pdv{S}{p}.
  \end{equation}
\end{proof}

In the coordinates of Lemma~\ref{l:two-freq-one-freq} the perturbed system rewrites as
\begin{equation}
  s' = 1 + \varepsilon f_s, \qquad
  h' = \varepsilon f_h, \qquad
  p' = - \pdv{S}{q} + \varepsilon f_p, \qquad
  q' = \pdv{S}{p} + \varepsilon f_q, \qquad
  z' = \varepsilon f_z
\end{equation}
where $(f_s, f_h, f_p, f_q, z)$ is the lift of the perturbation $f$ under the cover \[
  (s, h, p, q, z) \mapsto (p_1, p_2, q_1, q_2, z).
\]
This vector field is smooth. Taking $s$ as new time gives (denoting $\dot a = \dv{a}{s}$).
\begin{equation}
  \dot s  = 1, \qquad
  \dot h = \varepsilon g_h, \qquad
  \dot p = - \pdv{S}{q} + \varepsilon g_p, \qquad
  \dot q = \pdv{S}{p} + \varepsilon g_q, \qquad
  \dot z = \varepsilon g_z
\end{equation}
where
\begin{align}
\begin{split}
  g_h &= f_h/(1 + \varepsilon f_s), \\
  g_p &= (1 + \varepsilon f_s)^{-1} \Big( f_p + \pdv{S}{q} f_s \Big), \\
  g_q &= (1 + \varepsilon f_s)^{-1} \Big( f_q - \pdv{S}{p} f_s \Big), \\
  g_z &= f_z/(1 + \varepsilon f_s).
\end{split}
\end{align}
Thus perturbed two-frequency system is reduced to time-periodic perturbation of one-frequency system with Hamiltonian depending on additional parameter $h$ (that can be included in the vector $z$).

\end{appendices}

\printbibliography

\vskip 15mm

\noindent Anatoly Neishtadt,

\noindent {\small Department of Mathematical Sciences,}

\noindent {\small Loughborough University, Loughborough LE11 3TU, United Kingdom;}

\noindent {\small Space Research Institute, Moscow 117997, Russia}

\noindent {\footnotesize{E-mail : a.neishtadt@lboro.ac.uk}}

\vskip 5mm

\noindent Alexey Okunev,

\noindent {\small Department of Mathematical Sciences,}

\noindent {\small Loughborough University, Loughborough LE11 3TU, United Kingdom}

\noindent {\footnotesize{E-mail : a.okunev@lboro.ac.uk}}

\end{document}